\documentclass[11pt]{amsart}

\usepackage{lineno}

\usepackage{amssymb}
\usepackage{graphicx}
\usepackage{subfig}
\usepackage{tabularx}

\usepackage{bm}
\usepackage{mathrsfs}     
\usepackage{helvet}         
\usepackage{courier}        
\usepackage{type1cm}      
\usepackage[usenames, dvipsnames]{xcolor}
\usepackage{url}

\usepackage[percent]{overpic}

\usepackage{geometry,calc,color}
\usepackage{mathtools}
\usepackage[nointegrals]{wasysym}

\usepackage{enumerate}
\usepackage{arydshln}
\usepackage{siunitx}
\usepackage{booktabs}
\usepackage{multirow}
\usepackage{tikz}
\usepackage{pgfplots}
\usepackage[shortlabels]{enumitem}

\usepackage{appendix}
\pgfplotsset{compat=1.9}

\newtheorem{theorem}{Theorem}[section]

\newtheorem{problem}[theorem]{Problem}

\newtheorem{comment}[theorem]{Comment}
\newtheorem{remark}[theorem]{Remark}

\numberwithin{equation}{section}

\newcommand{\roundPrecision}{2}

\definecolor{ForestGreen} {cmyk}{0.91,0,0.88,0.12}
\definecolor{LemonChiffon}{rgb}{1.,0.98,0.8}
\definecolor{LightCoral}  {rgb}{0.94,0.5,0.5}
\definecolor{MediumYellow}{rgb}{1.,1.,0.50}
\definecolor{PaleGreen}   {rgb}{0.6,0.98,0.6}
\definecolor{PineGreen}   {cmyk}{0.92,0,0.59,0.25}
\definecolor{SeaGreen}    {cmyk}{0.69,0,0.50,0}

\newcommand{\blu}[1]{{\color{black} #1}} 
\newcommand{\rosso}[1]{{\color{black} #1}} 

\newcommand{\rev}[1]{{\color{black} #1}} 

\long\def\COMMENT#1{{\color{red} \begin{comment} #1 \end{comment}}} 
\long\def\NOTE#1{{\par\small\color{blue}  {\noindent{\bf Nota}.~ #1 }\par\medskip}} 

\long\def\COMMENT#1{}
\long\def\NOTE#1{}

\long\def\TESI#1{}

\usepackage{amssymb}
\usepackage{algorithm}
\usepackage[algo2e]{algorithm2e}

\newcommand{\numberset}{\mathbb}

\newcommand{\R}{\numberset{R}}

\newcommand{\PP}{\numberset{P}}

\newcommand{\B}{\numberset{B}}

\newcommand{\xK}{\x}

\renewcommand{\aa}{\Kmatr}

\newcommand{\abilinear}{a}
\newcommand{\ah}{a_h}
\newcommand{\aloc}{a^\elementVEM}
\newcommand{\ahloc}{a_h^\elementVEM}
\newcommand{\shloc}{S_h^\elementVEM}

\newcommand{\Ker}{Ker}

\newcommand{\LdueO}{L^2(\Omega)}

\newcommand{\HunozeroO}{H^1_0(\Omega)}
\newcommand{\elementVEM}{K}
\newcommand{\refelementVEM}{\widehat{\elementVEM}}
\newcommand{\mesh}{\mathcal{T}_h}
\newcommand{\meshfamily}{\{\mathcal{T}_h\}_h}
\newcommand{\cspace}{V}

\newcommand{\Puno}{\numberset{P}_1}
\newcommand{\grad}{\nabla}

\newcommand{\dweight}{\omega}

\newcommand{\bK}{{\partial\elementVEM}}
\newcommand{\Bkloc}{\B(\partial\elementVEM)}
\newcommand{\vemkloc}{V_1(\elementVEM)}
\newcommand{\vemk}{V_h}

\newcommand{\PiNabla}{\Pi^\nabla}

\newcommand{\dof}{dof}
\newcommand{\Span}{span}
\newcommand{\param}{\mu}
\newcommand{\uN}{u_\delta}
\newcommand{\VN}{\mathcal{V}_\delta}
\newcommand{\manifold}{\mathcal{M}}
\newcommand{\rbspace}{\mathcal{W}_M}
\newcommand{\funrb}{\xi}

\newcommand{\vv}{\mathbf{v}}
\newcommand{\x}{\mathbf{x}}
\newcommand{\Bcal}{\mathscr{B}}
\newcommand{\Acal}{\mathcal{A}}

\newcommand{\xiHat}{\widehat{\xi}}

\newcommand{\xhat}{\widehat{\x}}
\newcommand{\X}{\mathbf{X}}

\newcommand{\tolalfa}{t_\alpha}
\newcommand{\tolrho}{t_\rho}

\newcommand{\That}{\widehat{T}}

\newcommand{\vemunoloc}{\cspace_1(\elementVEM)}

\newcommand{\coef}{\gamma}

\newcommand{\dRB}{\sigma_M^{\text{rb}}}
\newcommand{\uhRB}{u_M^{\text{rb}}}

\newcommand{\Umatr}{\mathsf{U}}
\newcommand{\Bmatr}{\mathsf{B}}
\newcommand{\Cmatr}{\mathsf{C}}
\newcommand{\Amatr}{\mathsf{A}}
\newcommand{\Fmatr}{\mathsf{F}}
\newcommand{\Smatr}{\mathsf{S}}
\newcommand{\Kmatr}{\mathcal{K}}

\newcommand{\matrixB}{\mathsf{B}}

\newcommand{\xvec}{\mathsf{x}}

\newcommand{\norm}[1]{\left\lVert#1\right\rVert}

\newcommand{\W}{\mathcal{V}}
\newcommand{\parameters}{\mathscr{P}}

\newcommand{\aRB}{A}
\newcommand{\hTheta}{\widehat{\Lambda}}
\newcommand{\hu}{\widehat{u}}
\newcommand{\hv}{\widehat{v}}
\newcommand{\meshfine}{\mathcal{T}_\delta}

	\newcommand{\dd}{\delta}
\newcommand{\meshK}{\mathcal{T}^K_{\delta}}

\newcommand{\hatvv}{\widehat{\vv}}

\newcommand{\hatxiK}[1]{\widehat{d}_{#1}[K]}
\newcommand{\hatxijK}{\hatxiK{j,\delta}}

\newcommand{\stiff}{\Smatr}

\newcommand{\errstar}[1]{\mathtt{err^\star}(#1)}
\newcommand{\errzero}[1]{\mathtt{err^0}(#1)}
\newcommand{\erruno}[1]{\mathtt{err^1}(#1)}
\newcommand{\erren}[1]{\mathtt{err^\Kmatr}(#1)}
\newcommand{\errinf}[1]{\mathtt{err^\infty}(#1)}

\newcommand\parr{\mathtt{z}}
\newcommand\paru{\mathtt{z_1}}
\newcommand\paruu{\mathtt{z_2}}
\begin{document}
	
	\title[Reduced Basis method for VEM]{
		Reduced basis stabilization and post-processing\\ for the virtual element method
	}%
	
	\author[F. Credali]{Fabio Credali}
	\address{Dipartimento di Matematica ``F. Casorati'', University of Pavia (Italy), CEMSE division, King Abdullah University of Science and Technology, Thuwal (Saudi Arabia)}%
	\email{fabio.credali@kaust.edu.sa}%
	
	\author[S. Bertoluzza]{Silvia Bertoluzza}
	\address{IMATI ``E. Magenes'', CNR, Pavia (Italy)}%
	\email{silvia.bertoluzza@imati.cnr.it}%
	
	\author[D. Prada]{Daniele Prada}
	\address{IMATI ``E. Magenes'', CNR, Pavia (Italy)}%
	\email{daniele.prada@imati.cnr.it}%
	
	\subjclass{}%
	\keywords{}%

	\begin{abstract} 
		We present a reduced basis method for cheaply constructing (possibly rough) approximations to the nodal basis functions 
		of the virtual element space, and propose to use such approximations for the design of the stabilization term in the virtual element method and for the post-processing of the solution.
	\end{abstract}
	\maketitle
	
	\section{Introduction}\label{sec:intro}
	The recent years have seen an increasing interest in the design of numerical methods for the solution of partial differential equations based on the decomposition of the physical space domain with polytopal meshes. With respect to standard finite element meshes, these allow for a much greater flexibility that can be exploited e.g. for better dealing with complex geometries  or for more efficiently handling refinement and coarsening in adaptivity. Different approaches can be found in the literature, ranging  from the design of conforming discretizations based on the explicit definition and evaluation of local spaces by means of different types of barycentric coordinates \cite{Hormann:2017:GBC,sukumar2006recent,floater2015generalized}, to different non conforming approaches,  such as discontinuous Galerkin methods \cite{cangiani2017hp} and hybrid high order methods \cite{di2020hybrid},  where the solution is looked for in discontinuous polynomial spaces and where interelement continuity is only weakly imposed, by different strategies, up  to methods that do not rely on an underlying discretization space, such as  the mimetic finite difference method \cite{da2014mimetic}.
	
	In this landscape, the Virtual Element Method (VEM), which was introduced in \cite{beirao2013basic}  and have since already obtained considerable attention, aims at getting the best of both worlds. In such a method, the discretization space is a conforming one, that is, it is a subspace of the natural space for the continuous problem (e.g. $H^1$ for second order elliptic equations). However, the method systematically avoids the computationally expensive evaluation and handling of the basis functions. To this aim, leveraging a splitting of  the local function space into the direct sum of a polynomial component and a residual space, the VEM exactly evaluates only the polynomial component, while dealing with the residual component by means of suitable stabilization terms. The name ``virtual'' reflects the core idea at the basis of the method: the conforming discretization space, which underlies the method, is virtually there, providing the foundations for the design of the method and  its theoretical analysis, but is never fully evaluated in the actual numerical computations. 
	
	The resulting algorithms turn out to be surprisingly robust: optimal convergence, which can be proven under somewhat restrictive assumptions on the shape regularity of the polytopal mesh, can be observed in practice also for very badly shaped meshes, at least for isotropic problems. Also thanks to this robustness, the method has quickly gained the attention of the scientific community, as demonstrated by the ever increasing literature focused on different aspects:  theoretical analysis \rosso{\cite{beirao2017stability, cangiani2017posteriori, brenner2017some, chen2018some,  bertoluzza2022weakly, bertoluzza2022interior}}, implementation techniques \cite{sutton2017virtual,mengolini2019engineering,yu2022mvem}, preconditioning \rosso{\cite{bertoluzza2017bddc,prada2019feti,bertoluzza2020feti,dassi2020parallel,dassi2020parallel2,dassi2022robust}}, different generalizations and extensions \rosso{\cite{brezzi2014basic,beir2016basic,de2016nonconforming,da2016serendipity,berrone2021lowest}}, as well as applications in several  fields, mainly in engineering:  linear elasticity problems \cite{da2013virtual,gain2014virtual,nguyen2018virtual}, plate bending \cite{brezzi2013virtual}, to reaction--diffusion problems with variable coefficients \cite{beirao2016virtual}, \rosso{ Helmholtz equation \cite{mascotto2019nonconforming,mascotto2019nonconforming2,mascotto2022nonconforming}, fluid dynamics and porous media \cite{antonietti2014stream,da2017divergence,vacca2018h,DASSI20221,berrone2022virtual}, contact and deformation problems \cite{wriggers2017low,wriggers2017efficient,cihan2022virtual}, as well as geophysical applications and discrete fracture networks problems \cite{benedetto2016hybrid,berrone2022efficient,berrone2023virtual}.}  
	
	\
	
	The  virtual approach, characteristic of the VEM, has, however, also some downsides. First of all, while  the numerical solution belongs to a conforming discretization space (i.e., it is continuous), it cannot generally be accessed directly without solving a PDE in each polytopal element (we recall that the elemental VEM functions are themselves solution of suitable partial differential equations). Only a projection onto a discontinuous polynomial space can be actually accessed, \rev{resulting, for most practical purposes, in a discontinuous approximation of the solution. Among other things, this turns out to be an obstacle to  the design of geometric multigrid preconditioners for VEM~\cite{multi,antonietti2023agglomeration}, where the most natural prolongation operators would require the evaluation of the discrete functions at points interior to the elements.}	Moreover, it has been observed that the introduction of the (inherently isotropic) stabilization term, aimed at handling the non polynomial residual component of the solution, injects a measure of isotropy in the discrete problem that, when dealing with anisotropic problems, somehow pollutes the results. To overcome these limitation, and fully exploit the potential of the method, it would be desirable to have at our disposal a possibly cheap method to locally evaluate (a more or less accurate approximation of) the local basis functions, if, when and where needed. What we propose here, is to resort to the reduced basis  method to that aim.

\

Indeed, the Reduced Basis (RB) method \cite{hesthaven2016certified} aims at the efficient numerical solution of para\-metrized partial differential equations, and it is specifically designed for addressing the need of repeatedly solving the same  equation for a large number of different values of the parameters. In the RB method, the discretization space is spanned by a small number of (linear combinations of) solutions of the equation with suitable values of the parameter (typically, these are randomly generated according to some probability law), which are evaluated in a computationally intensive offline phase, to be carried out once and for all. In our framework, the idea (see \cite[Section 6.2]{hesthaven2016certified}) is to treat the shape of the polygonal elements as parameter, by rewriting, by means of a change of variable, the elemental PDE defining a VEM discrete function as a parametrized equation on a fixed reference element. This approach will provide us with a tool allowing to efficiently  reconstruct the non polynomial part of virtual element functions more or less accurately, and compute different quantities of interest, starting from the stiffness matrices.

\

In the design and implementation of the VEM, this tool, that thanks to the local nature of the VEM is fully parallelizable, can be used in different ways, with different aims, and with different accuracies and computational costs (we can choose different values of the dimensions $M$ of the RB space use to construct the VEM functions). 
It can  be used in the postprocessing phase to obtain, in output, an actual $H^1$ (for second order problem) function, that can be handled for performing tasks such as visualization, pointwise evaluation, full $H^1$ error evaluation when benchmarking the method (we recall that the standard error evaluation in  the VEM normally relies on broken norms for the discontinuous polynomial part). 
It can also be exploited in the design of the discrete bilinear form.  
More precisely, true to the core idea of virtual elements, we  propose to exploit the RB reconstruction of the VEM basis functions only for the design of a better stabilization term, in those cases where the stabilization recipes available in the literature underperform (as it happens for anisotropic problems). In such a case, we will not need to reconstruct the basis functions  with particular accuracy (we will see that choosing $M = 1$ will be sufficient to improve the performance of the method for anisotropic PDEs!). This approach is similar to the one proposed in  \cite{wriggers2017low}, where a coarse finite element discretization is leveraged for  dealing with the non polynomial component in nonlinear virtual elements for finite deformation. However, it would be possible to use the RB reconstruction of the VEM basis functions to design a fully conforming discretization method based on the VEM function space (which, for the lowest order VEM, coincides with the polygonal FEM space with harmonic generalized barycentric coordinates \cite{DeRose2006HarmonicC}). In such case, the VEM machinery could still be exploited for the implementation of the method: similarly to the approach of \cite{manzini2014new}, the polynomial component can be exactly evaluated and handled as in the virtual element method, while the non polynomial component is efficiently handled by the RB method.

\

The paper is organized as follows: in Section \ref{sec:vem} we will recall the definition and useful properties of the lowest order Virtual Element Method, on which the paper is focused. In Section \ref{sec:thebasis} we will focus on the partial differential equations defining the nodal basis for the elemental virtual element space, and on its interpretation as a parametric equation on 
fixed reference elements. In Section \ref{sec:rb} we will briefly recall the basic ideas underlying the Reduced Basis method, which we will apply for the computation of the VEM basis functions in Section \ref{sec:vemrb}. In Section \ref{sec:validation} we will present some numerical tests aimed at the validation of the RB VEM function reconstruction tool. In Sections \ref{sec:stab} and \ref{sec:postproc} we will show some applications to the VEM  stabilization in the  anisotropic case and to the postprocessing of the VEM solution. Section~\ref{sec:conclusion} presents some conclusions and discusses future perspectives.


	\section{The lowest virtual element method}\label{sec:vem}
	Hereon, we will use the following standard notation for functional spaces \cite{lions2012non}. If $D$ is an open bounded domain, $L^2(D)$ denotes the space of square integrable functions endowed with scalar product $(\cdot,\cdot)_D$. Moreover, $H^1(D)$ denotes the Sobolev space of square integrable functions with square integrable gradient, while $H^1_0(D)$ is the subspace of $H^1(D)$ of functions with zero trace on the boundary $\partial D$. The space of continuous functions over $D$ is denoted by $C^0(D)$, while $\Puno(D)$ is the space of polynomials of degree $\le 1$ in $D$.
	
	\

	Let us consider a simple model problem, namely a second order diffusion equation with homogeneous Dirichlet boundary conditions in a polygonal domain $\Omega \subset \R^2$, that we assume to be open, bounded and connected. In strong form, the equation reads as 
	\begin{equation}\label{strong}
		\left\{\begin{aligned}
			-\nabla \cdot \aa \nabla u &= f, \mbox{\quad in }\Omega,\\
			u&=0, \mbox{\quad on }\partial \Omega,\\
		\end{aligned}\right.
	\end{equation}
where $f \in L^2(\Omega)$ denotes the given right hand side, and where $\aa$ is a possibly non symmetric positive definite matrix, which, for the sake of simplicity, we assume to be constant. Letting
\[
\abilinear(u,v) = \int_\Omega \aa \grad u\cdot\grad v\,d\x,
\]
the variational counterpart of \eqref{strong} is given by the following problem, which is well known to be well posed, $\abilinear$ being continuous and coercive \cite{lions2012non}.
	\begin{problem}\label{pro:cont_poisson}
		Find $u\in\cspace=\HunozeroO$ such that
		\begin{equation}
			a(u,v) = \int_\Omega fv\,d\x, \quad \forall v\in\cspace.
		\end{equation}
	\end{problem}
	\TESI{We define the bilinear form $\abilinear:\cspace\times\cspace\longrightarrow\R$ and the linear functional of the right hand side $F:\cspace\longrightarrow\R$ as
	\begin{equation}
		\begin{aligned}
		&a(u,v) = \int_\Omega\aa \grad u\cdot\grad v\,d\x \\
		&F(v) = \int_\Omega fv\,d\x
		\end{aligned}
	\end{equation}
	assuming that $f\in\LdueO$.}
	We consider, in this paper, the lowest order virtual element method, of which, hereafter, we recall the definition and main properties. 	Let $\meshfamily$ denote a family of decompositions of $\Omega$, each made up of a finite number of polygons $\elementVEM$, which we assume to be star shaped.

Following \cite{beirao2013basic}, as for the finite element method, the virtual element space is first defined locally on each element $\elementVEM$.  The local spaces are then glued together to form the global space.
		Let us then consider a generic polygonal element $\elementVEM\in\mesh$. 
	The local lowest order VEM space is defined as 
\begin{equation}
	\label{eq:local_vem_space}
	\vemkloc = \{ v\in H^1(\elementVEM):v_{|\partial K}\in\Bkloc,\,\Delta v_{|\elementVEM}=0\}\mbox{,}
\end{equation}		
	with
	\begin{equation}
		\label{eq:bolle}
		\Bkloc = \{ v\in C^0(\partial\elementVEM): v_{|e}\in\Puno(e)\quad \text{for all edge } e\in\partial\elementVEM\}.
	\end{equation}
Each function $v\in\vemkloc$ can be identified via the set of degrees of freedom given by the values of $v$ at the vertices of $\elementVEM$, which is unisolvent for $\vemkloc$.  Notice that $\Puno(K)\subseteq \vemkloc$.
Glueing by continuity all the local spaces, we finally define the global lowest order virtual element space as
\begin{equation}
	\label{eq:global_vem_space}
	\vemk = \{ v\in V:  v_{|\elementVEM}\in \vemkloc\quad \text{for all }\elementVEM\in\mesh\}\subset\cspace\mbox{.}
\end{equation}
Such a space is	associated to the set of global degrees of freedom given by the values of $v\in\vemk$ at  all the internal vertices of the decomposition $\mesh$.



\

	Having defined the space, we next introduce a discrete version of the bilinear form. We start by splitting the bilinear form $\abilinear$ into the sum of local contributions $\aloc$, that is
	\begin{equation}\label{continuous_weak}
	\abilinear(u,v) = \sum_{\elementVEM\in\mesh} \aloc(u,v),
	\end{equation}
	where $\aloc:H^1(\elementVEM)\times H^1(\elementVEM)\rightarrow\R$ is defined as
	\begin{equation}
		\aloc(u,v) = \int_\elementVEM \aa \grad u\cdot\grad v\,d\x.
	\end{equation}

\TESI{With the same reasoning, we can introduce the discrete bilinear form $\ah:\vemk\times\vemk\longrightarrow\R$ such that
	\begin{equation}
		\ah(u_h,v_h) = \sum_{\elementVEM\in\mesh} \ahloc(u_h,v_h) \quad\forall u_h,v_h\in\vemk.
	\end{equation}
	where $\ahloc:\vemkloc\times \vemkloc\longrightarrow\R$ is its local counterpart.
	}
The definition of the virtual element method stems from the observation that, while evaluating local forms $a^K(u,v)$ given the values of the degrees of freedom for $u$ and $v$ would generally require solving PDEs in the elements $K$,  no PDE needs to be solved to  compute exactly $\aloc(q,v_h)$ for all $q\in\Puno(\elementVEM)$ and $v_h\in\vemkloc$ (in the virtual elements terminology we say that this quantity is ``computable''). 
		In order to build a computable approximate bilinear form $\ah$, we introduce the local projection operator
	\begin{equation}
		\PiNabla:\vemkloc\rightarrow\Puno(\elementVEM)\subset\vemkloc
	\end{equation}
defined as 
		\begin{equation}\label{eq:PiELk_def}			
						\int_{\elementVEM} \nabla \PiNabla v_h\cdot \nabla q \,d\x= \int_{\elementVEM} \nabla v_h \cdot \nabla q\,d\x \quad \forall q\in\Puno(\elementVEM),\qquad
						\int_{\bK} \PiNabla v_h \,ds= \int_{\bK} v_h\,ds.
			\end{equation}
		We recall that the action of $\PiNabla$ on the elements of $\vemkloc$ is computable (see, e.g., \cite{beirao2014hitchhiker}).
	We then consider a decomposition of $\vemkloc$ as
	\begin{equation}\label{splitting}
		\vemkloc = \Puno(\elementVEM)\oplus V_\perp(\elementVEM), \qquad \text{ with }V_\perp(\elementVEM) = \ker \PiNabla \subset \vemkloc,
	\end{equation}
	where $V_\perp(\elementVEM)$ is 
the kernel of the projector $\PiNabla$. We observe that,	
	 splitting $u_h,v_h\in\vemkloc$ into a polynomial part plus a non polynomial contribution according  to \eqref{splitting}, we can write
	\begin{multline}\label{splitting:2}
		\aloc(u_h,v_h) = \aloc(\PiNabla u_h,\PiNabla v_h) + \aloc(\PiNabla u_h, (I-\PiNabla) v_h)\\ + \aloc((I-\PiNabla) u_h ,\PiNabla v_h) + \aloc((I-\PiNabla)u_h,(I-\PiNabla)v_h) \\= \aloc(\PiNabla u_h,\PiNabla v_h) + \aloc((I-\PiNabla)u_h,(I-\PiNabla)v_h),
	\end{multline}
where the last equality descends from the fact that, if $\aa$ is constant, by definition, $\nabla (I-\PiNabla) v_h$ is ortogonal to  $\aa \nabla q$ for all order one polynomials $q$ (indeed it is easy to see that there exists an order one polynomial $q'$  such that $\aa \nabla q = \nabla q'$).
The first  term on the right hand side of \eqref{splitting:2} can be computed exactly directly from the value of the degrees of freedom, while the second term is not computable. In defining the VEM it is then replaced by a computable \textit{stabilization term} endowed with appropriate properties. More precisely, we  replace the last term in the sum with \begin{equation}
		\shloc((I-\PiNabla)u_h,(I-\PiNabla)v_h)
	\end{equation}
	where $\shloc$ is any computable semi scalar product, inducing an equivalent $H^1(K)$ semi norm on $V_\perp$(K). In other words,  $\shloc$ is a
	 symmetric, positive definite bilinear form  defined on $\vemkloc$ such that there exist two constants $c_\star$ and $c^\star$ so that	for all $v_h \in { V_\perp(\elementVEM)}$ we have
	\begin{equation}\label{spectral_equivalence}
		c_\star  a^K(v_h,v_h) \le \shloc(v_h,v_h) \le c^\star  a^K(v_h,v_h).
	\end{equation}
	The local virtual element bilinear form  is defined as
	\[
	\ahloc(u,v) = \aloc(\PiNabla u, \PiNabla v) + \shloc((I-\PiNabla)u,(I- \PiNabla)v).
	\]
	Since also the linear operator at the right hand side of equation \eqref{continuous_weak} is not computable, we need as well to introduce an approximate linear operator $F_h: V_h \to \mathbb{R}$, which we define as
\[
F_h(v_h) = \sum_{K\in \mesh} F_h^K (v_h) = \sum_K | \bK |^{-1}  \int_{K} f \,d\x \int_{\bK}  v_h \,ds.
\]

The discrete problem then reads:

\begin{problem}\label{pro:disc_poisson}
	Find $u_h\in \vemk$ such that
	\begin{equation}
		a_h(u_h,v_h) = F_h(v_h) \quad \forall v_h \in\vemk.
	\end{equation}
\end{problem}

\newcommand{\vPio}{\Pi^0}

\begin{remark}
An alternative, often preferred, way of defining the approximate local bilinear form $\aloc$ is (see \cite{beirao2016virtual})
\[
\aloc(u_h,v_h) = \int_K \aa \vPio(\nabla u_h) \cdot \vPio(\nabla v_h)\,d\x + \shloc((I-\PiNabla)u,(I- \PiNabla)v),
\]
where $\Pi^0: L^2(K)^2 \to \Puno(K)^2$ is the $L^2(K)$ projection, which is also computable. We remark however that, while for general order $k$ virtual element spaces, the two definitions yield distinct bilinear forms, in the case of the lowest order VEM we have that
\(
\vPio \nabla u_h = \nabla \PiNabla u_h
\), so that the methods resulting from the two definitions coincide.
\end{remark}

	There are several possible choices for the stabilization term, depending also on the problem under consideration \cite{beirao2013basic,cangiani2015hourglass,mascotto2018ill}, see also \cite{bertoluzza2022stabilization}. In our implementation, we consider the simplest, so called \textit{dofi--dofi} stabilization, which is defined, for the lowest order VEM, as
	\begin{equation}
		\shloc(v_h,w_h) = \sum_{i=1}^{N} \dof_i(v_h)\,\dof_i(w_h) = \sum_{i=1}^{N} v_h(\vv_i)\,w_h(\vv_i) \quad\forall v_h,w_h\in\vemkloc,
	\end{equation}
where $\vv_1,\cdots,\vv_N$ are the $N$ vertices of $K$, and	where $\dof_i$ is the operator that associates to each smooth function $\varphi$ the $i$--th local degree of freedom on $\elementVEM$, that is, the value of $\varphi$ at the vertex $\vv_i$. We also consider the so called {\em D-recipe} stabilization \rosso{(\cite{mascotto2018ill}),} that is defined as
		\begin{equation}
		\shloc(v_h,w_h)  = \sum_{i=1}^{N} \dweight_i \, v_h(\vv_i)\,w_h(\vv_i) \quad\forall v_h,w_h\in\vemkloc,\quad \dweight_i = \max\{1, a^K(\PiNabla \rosso{e_i}, \PiNabla \rosso{e_i}) \},
	\end{equation}
	where \rosso{$e_i \in \vemkloc$} is the basis function for the local VE space corresponding to the node $\vv_i$.

	\

Provided \eqref{spectral_equivalence} holds, we have the following error bound (see \cite{beirao2013basic,beirao2016virtual}): $u \in H^2(\Omega)$ implies
\[
\| u - u_h \|_{1,\Omega} \lesssim h | u |_{2,\Omega} + \mathfrak{E}(f),
\]
with $\mathfrak{E}(f) = \sup_{v\in H^1(\Omega)} | \langle f - f_h , v \rangle |/\| v \|_{1,\Omega}$.

\

\NOTE{
The error estimate for the general case is exactely the same as for the laplacian. Indeed we have 
\begin{gather}
	\| u_h - I_h u \|_{1,\Omega}^2 \lesssim a_h(u_h - I_h u, e_h)  = 
	\langle f_h , e_h \rangle - a_h(I_h u - u_\pi,e_h) + a(u_\pi - u,e_h) - \langle f , e_h \rangle \\
	\lesssim \| f_h - f \|_{-1,\Omega} \| e_h \|_{1,\Omega} 
	+ \| I_h u - u_\pi \|_{1,\Omega} \| e_h \|_{1,\Omega} + \| u_\pi - u \|_{1,\mathcal{T}_h} \| e_h \|_{1,\Omega}.
	\end{gather}
which gives us 
\[
	\| u_h - I_h u \|_{1,\Omega} \lesssim \| f_h - f \|_{-1,\Omega} 
+ \| I_h u - u_\pi \|_{1,\Omega}  + \| u_\pi - u \|_{1,\mathcal{T}_h}.
\]
}

\section{The VE nodal basis functions as solutions to  parametric equations}\label{sec:thebasis}

As it happens for the finite element method, the local VE space $\vemkloc$ is endowed with a nodal basis $\{e_j, \ j=1,\cdots ,N\}$, $N$ denoting the number of vertices of the  polygonal element $\elementVEM$, such that all functions $w_h \in \vemkloc$ can be written as
\[
w_h = \sum_{j=1}^N w_h(\vv_j) e_j.
\]
The basis functions $e_1,\dots,e_N$ are the solutions of the following problem, where $\delta$ is the Kronecker delta.
	\begin{problem}\label{pro:basis1}
		For $j=1,\dots,N$, find $e_j$ such that
		\begin{equation}\label{truebasisfunctions}
			\left\{\begin{aligned}
				-\Delta e_j&= 0 &&\mbox{ in }\elementVEM\\
				e_j &= g_j &&\mbox{ on }\partial\elementVEM
			\end{aligned}\right.
		\end{equation}
		with $g_j \in \Bkloc$   such that $g_j(\vv_i)= \delta_{i,j}$ for $i=1,\dots,N$.
	\end{problem}
In the finite element method, the explicit knowledge of the nodal basis functions is exploited for different goals, particularly for the evaluation of the entries of the local stiffness matrix, but also, in the framework of the post-processing of the solution, for the evaluation of quantities of interest for the final user, such as the value of the solution at a given point or along a given line, just to make an example. Conversely, the idea at the core of the virtual element method is to never solve the above equations, and carry on without the explicit knowledge of the nodal basis functions. The stiffness matrix is replaced by the approximate stiffness matrix, directly computed in terms of degrees of freedom via polynomial projections plus stabilization, and, once the discrete problem is solved, the degrees of freedom of the discrete solution, which are the quantities that are actually computed, give the final user access to one or more projections of the solution onto spaces of discontinuous piecewise polynomials. 

Depending on the problem and on the final goal of the computation, this purely virtual approach that completely avoids constructing (some approximation of) the non polynomial component of the local basis functions might not be enough. For example, when dealing with anisotropic problems, the inherently isotropic stabilization term might hinder the performance of the method, or  the final user might wish for a continuous discrete solution.  

In such cases, a cheap way of (locally) reconstructing, more or less accurately, virtual functions, given the values of the degrees of freedom, would be desirable. For instance, even a very rough approximation of the non polynomial component of the virtual basis functions might allow to design anisotropic stabilization terms, while a more accurate reconstruction can be leveraged to retrieve point values of the continuous discrete solution within the elements. We propose here to carry out such a task by locally resorting to a model order reduction, in the spirit of the reduced basis method. To this aim, for each $N$,  we can look at the collection of Problems of the form \ref{pro:basis1} for different $N$ vertices polygons as an equation parametrized by the geometry. More precisely, we reformulate Problem \ref{pro:basis1}, in which the solution space depends on the geometry, as a parameter dependent problem on a reference element, set in a geometry independent space,  the parameter being the element $K$ itself (see \cite[Chapter 6.2]{hesthaven2016certified}).

We introduce the $N$ vertices regular polygon $\refelementVEM$ with  unit diameter and centered in the origin. We let $\hatvv_1,\cdots,\hatvv_N$ denote the $N$ vertices,  ordered counterclockwise,  and $\xhat_\elementVEM=(0,0)$ the barycenter. As the behaviour of the solution of Problem \ref{pro:basis1} with respect to translations and rescaling of the domain is well understood, we restrict ourselves to polygons {having circumscribed circle with unit diameter}, centered in the origin. We then introduce the following parameter space:
\begin{equation}\label{defpolyparameters}
\parameters = \{K: \ K\text{ polygon with $N$ \rev{vertices} with }\Ker(K) \not=\emptyset,\  h_K = 1, \ \xK_K  = (0,0)\},
\end{equation}
where $\Ker(K)$ denotes the {\em kernel} of $\elementVEM$ (that is the  set of points with respect to which $\elementVEM$ is star shaped), $\x_\elementVEM$ denotes the barycenter of $\Ker(K)$, and $h_K$ the diameter of the circumscribed circle.
Remark that, depending on the characteristics of the tessellations considered, one might further restrict the parameter space. For example, when handling Voronoi tessellations, one might add a convexity condition to the definition of the parameter space $\parameters$.

\begin{figure}
	\begin{overpic}[trim = 100 200 100 200,width=15cm]{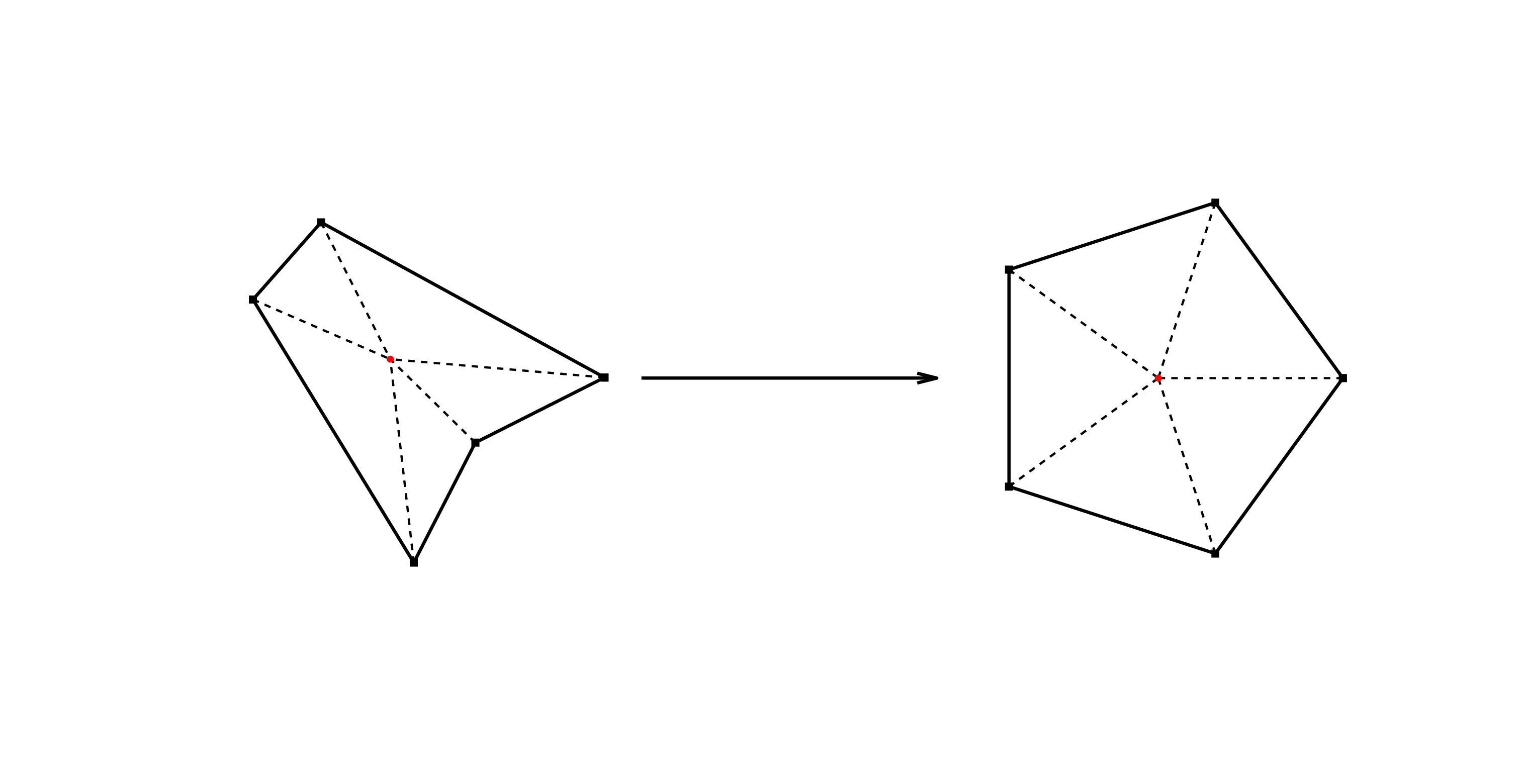}
		\put (71,14) {$\That_i$}
		\put (19,12) {$T_i$}
		\put (22,17.5) {$\xK_\elementVEM$}
		\put (82,16) {$\xhat$}
		\put (50,17) {$\Bcal_\elementVEM$}
		\put (9,25) {$\elementVEM$}
		\put (93,25) {$\refelementVEM$}
	\end{overpic}
	\caption{An example of affine mapping $\Bcal_\elementVEM$  between a random \rev{star shaped} pentagon and its regular counterpart.}
	\label{fig:affine_mapping}
\end{figure}
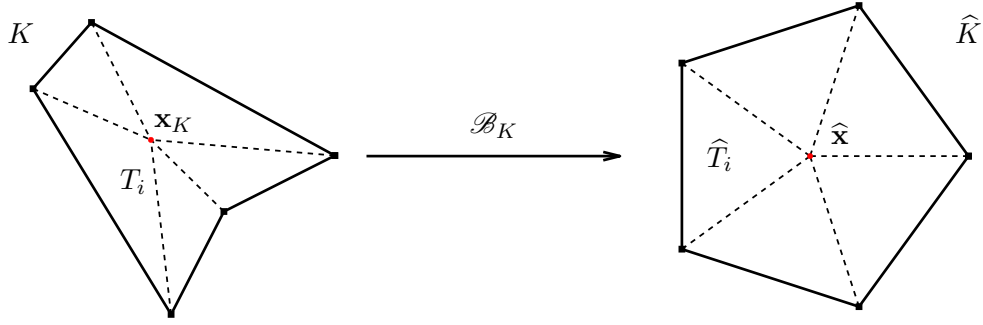

\rev{We can construct a piecewise affine transformation between the star shaped polygons $\elementVEM\in\parameters$ and the reference polygon $\refelementVEM$}. To this aim, as shown in Figure \ref{fig:affine_mapping}, we partition both $\elementVEM$ and $\refelementVEM$ in as many triangles as there are edges. More precisely, assuming that the vertices $\vv_1$, $\cdots$,$\vv_N$ of $K$  are also ordered counterclockwise, we let $T_i$ and $\That_i$ denote respectively the triangles with vertices 
$\vv_i$, $\vv_{i+1}$, \rosso{$\xK_K$}  and $\widehat{\vv}_i$, $\widehat{\vv}_{i+1}$, $\xhat=(0,0)$, with the convention that $\vv_{N+1}=\vv_1$ (resp. $\hatvv_{N+1} = \hatvv_1$),
and  we have the decomposition
\begin{equation}\label{eq:partition}
	\elementVEM = \bigcup_{i=1}^N T_i \quad \text{and} \quad \refelementVEM = \bigcup_{i=1}^N \That_i.
\end{equation}
We next introduce the continuous piecewise affine transformation
\begin{equation}
	\Bcal_\elementVEM:\elementVEM\longrightarrow\refelementVEM, \qquad \Bcal_\elementVEM(\x) = \Bmatr_{\elementVEM,i}\x \quad \text{on } T_i
\end{equation}
with $\Bmatr_{\elementVEM,i}$ invertible $2 \times 2$ matrix,	defined in such a way that $\Bcal_\elementVEM(\vv_i)=\widehat{\vv}_i$ for all the vertices $\vv_i$ of $K$, and $\Bcal_\elementVEM(\xK ) = \xhat$.  By construction, $\That_i=\Bcal_\elementVEM T_i$ for $i=1,\dots,N$. An easy calculation yields the following expression for the matrix $\Bmatr_{\elementVEM,i}$:
\begin{equation}
	\Bmatr_{\elementVEM,i}
	=
	\begin{bmatrix}
		x_i & x_{i+1} \\
		y_i & y_{i+1}
	\end{bmatrix} \begin{bmatrix}
		\widehat{x}_i & \widehat{x}_{i+1} \\
		\widehat{y}_i & \widehat{y}_{i+1}
	\end{bmatrix}^{-1},
\end{equation}
where $(\widehat{x}_i,\widehat{y}_i)$ and $(x_i,y_i)$ denote the coordinates of the $i$--th vertex $\widehat{\vv}_i$ and $\vv_i$ of $\refelementVEM$ and $\elementVEM$ respectively.

\

We can write Problem \ref{pro:basis1} in weak form as: find $e_j \in H^1(K)$ such that
\begin{equation}\label{eq:postprocessing}
e_j = g_j \quad\text{ on }\partial K, \qquad \text{ and } \qquad	\int_K \nabla e_j \cdot \nabla v = 0, \quad  \text{ for all }v \in H^{1}_0(K). 
	\end{equation}

By performing a change of variables  in the integrals in equation \eqref{eq:postprocessing}, we can transfer the  Poisson Problem \ref{pro:basis1} to a problem in $H^1(\refelementVEM)$.  Indeed we have
\begin{equation*}
	\int_\elementVEM\grad u\cdot\grad v\,d\x = \sum_{i=1}^{N} \int_{T_i}\grad u\cdot\grad v\,d\x= \sum_{i=1}^{N} \int_{\That_i} \rosso{| \det(\Bmatr^{-1}_{\elementVEM,i})|}\, \Bmatr_{\elementVEM,i}^\top \Bmatr_{\elementVEM,i}\,\grad\widehat{u}\cdot\grad\widehat{v}\,d\xhat.
\end{equation*}

\TESI{where we point out that the matrix $\Bmatr_{\elementVEM,i} \Bmatr_{\elementVEM,i}^\top$. }
Then, introducing the parameter dependent bilinear form
\begin{equation}
	\aRB(\widehat u,\widehat v;\elementVEM) = \sum_{i=1}^{N} \int_{\That_i} \rosso{| \det(\Bmatr^{-1}_{\elementVEM,i})|}\, \Bmatr_{\elementVEM,i}^\top \Bmatr_{\elementVEM,i}\,\grad\widehat{u}\cdot\grad\widehat{v}\,d\xhat,
\end{equation}
we have the following equivalent formulation for Problem \ref{pro:basis1}, in the form of a parametric equation on the reference element:
find $\widehat{e}^K_j \in H^1(\refelementVEM)$ such that for all $\widehat{v} \in  H^1_0(\refelementVEM)$ 
\begin{equation}\label{eq:postprocessinghat}	
\widehat e^K_j = \widehat g_j \quad\text{ on }\partial \widehat K, \qquad \text{ and }\qquad	\aRB(\widehat e^K_j,\widehat v;\elementVEM)		 = 0, \quad  \text{ for all } \widehat v \in H^1_0(\widehat{K}), 
\end{equation}
where  $\widehat{g}_j$ is the piecewise linear function on $\partial \refelementVEM$ such that $\widehat g_j(\widehat \vv_i)= \delta_{i,j}$ for $i=1,\dots,N$.

\

As already anticipated, the idea is now to solve such a problem for all the $N$ edges elements in the mesh, by resorting to the reduced basis method, which was devised as an efficient way to solve  parameter dependent partial differential equations for a large number of different instances of the parameter.


	\section{The reduced basis method}\label{sec:rb}
We devote this section to  recalling the main features of the reduced basis method, and we refer the reader to \cite{hesthaven2016certified} for further information. 

\

Letting $\W_\delta$ denote an Hilbert space, we let $\aRB(\cdot,\cdot;\mu): \W \times \W \to \mathbb{R}$  denote a parameter dependent bilinear form, which we assume to be continuous and coercive for all values $\mu\in \parameters$ of a vector parameter $\mu$ in the given parameter set $\parameters \subseteq \mathbb{R}^L$. Analogously, let $F(\cdot;\mu): \W \to \mathbb{R}$ denote a parameter dependent bounded linear operator. We consider the following class of parameter dependent problems.
\begin{problem}\label{pro:model_rb}
	Find $u[\param]\in \W$ such that 
	\begin{equation}\label{eq:model_rb}
		\aRB(u[\param],v;\param) = F(v;\param) \quad \forall v\in \W.
	\end{equation}
\end{problem}

\newcommand{\NN}{\mathcal{N}}

Given a 
 finite dimensional approximation space $\VN$, with
\begin{equation}
	\VN = \Span\{\varphi_1,\dots,\varphi_\mathcal{N}\} \subset \W,
\end{equation}  for each value of the parameter $\mu \in \parameters$, let   $\uN[\param]\in\VN$ denote an approximation to $u[\param]$, computed by the preferred method (this is usually, but not necessarily,  a Galerkin method).
The assumption underlying the reduced basis method is that the (approximate) solution manifold, that is
the set 
\begin{equation}
	\manifold = \{ u_\delta[\param] : \param \in \parameters \} \subset \W_\delta,
\end{equation} 
of all the discrete solutions of the problem as the parameter varies,
can be approximated by a lower dimensional space
\begin{equation}
	\rbspace = \Span\{\funrb_1,\dots,\funrb_M\} \subset \VN,
\end{equation}
with $M\ll \mathcal{N}$, where, for $\ell =1,\cdots, M$, the functions $\funrb_\ell$ are linear combinations of discrete snapshots $u_\delta[\param_k]$, computed offline for some suitably chosen values $\param_k\in \parameters$.

\

Assuming that the functions $\funrb_\ell$, $\ell = 1,\cdots,M$,   have been computed in an offline phase, for new values $\param$ of the parameter we look for the corresponding approximate solutions of Problem~\ref{pro:model_rb} in $\rbspace$ 
by a Galerkin method: the solution of the form $u^\text{rb}[\param] = \sum_{\ell=1}^M x_\ell \xi_\ell$ can be computed by solving a (small) linear system with unknown $\xvec = (x_\ell)_\ell$, of the form
\begin{equation}\label{linsysrb}
\Amatr[\mu]\, \xvec = \Fmatr[\mu], 
\end{equation}
where the matrix $\Amatr[\mu]$ and the vector  $\Fmatr[\mu]$ are defined as
\[
\Amatr[\mu] = (\aRB(\xi_{\ell'},\xi_{\ell};\mu))_{\ell,\ell'}, \qquad \Fmatr[\mu] = (F(\xi_\ell;\mu))_{\ell},
\]
where we use the notation $(X(\ell,\ell'))_{\ell,\ell'}$ to denote the matrix whose $\ell$-th row, $\ell'$-th column entry is $X(\ell,\ell')$, and use an  analogous notation for vectors. While one would think that assembling the linear system \eqref{linsysrb}, which requires evaluating the bilinear  and linear forms  on elements of the  space $\VN$, can a priori be quite expensive if $\NN$ is very large,
under suitable assumptions, the matrix $\Amatr[\mu]$ and the right hand side vector $\Fmatr[\mu]$ can instead be assembled very cheaply, by relying on the pre-computation of a number of quantities, carried out offline at the time of the construction of the basis \{$\xi_\ell\}$. More precisely, assuming that  the bilinear form $\aRB(\cdot ,\cdot; \mu)$ and right hand side $F(\cdot;\mu)$ allow an affine decomposition, that is, that
there exist parameter independent bilinear forms $\aRB^q$, $q = 1, \cdots, Q_A$ and linear operators $F^q$, $q = 1, \cdots, Q_F$ such that $A$ and $F$ can be decomposed as
\begin{equation}\label{eq:affine}
	\aRB(u,v;\param) = \sum_{q=1}^{Q_A} \alpha^q_A[\param]\aRB^q(u,v),
	\qquad
	F(v;\param) = \sum_{q=1}^{Q_F} \alpha^q_F[\param]F^q(v),
\end{equation}
where $\alpha^q_A: \parameters \to \R$ and $\alpha^q_F: \parameters \to \R$ are given functions,
we can precompute and store the matrices $\Amatr^q = (\aRB^q(\funrb_{\ell'},\funrb_\ell))_{\ell,\ell'}$ and $\Fmatr^q = (F^q (\funrb_\ell))_{\ell}$.
The matrix $\Amatr[\mu]$ and the right hand side $\Fmatr[\param]$ can then be obtained as
\[
\Amatr[\mu] = \sum_{q=1}^{Q_A} \alpha^q_A[\param] \Amatr^q, \qquad \Fmatr[\mu] = \sum_{q=1}^{Q_F} \alpha^q_{F}[\param] \Fmatr^q,
\]
where only the coefficients $\alpha^q_A[\param] $ and $\alpha^q_{F}[\param]$ have to be computed for each new value of the parameter.

\


%
The computational intensive phase of the procedure is the offline phase, where the initial snapshots of the solution are computed by numerically solving a number of instances of Problem \ref{pro:model_rb} in the large discrete space $\VN$.
To start,  a set of trial parameters
\begin{equation}
	S = \{\param_1,\dots,\param_P\} \subset \parameters
\end{equation}
is selected. 
A common strategy is to randomly choose the elements of $S$ according to a certain probability distribution. Hopefully, if $S$ is \rosso{large} enough, the subspace spanned by the set $\{ u[\param] : \param\in S \}$  allows for a good representation of $\manifold$.  For all $\param_\ell\in S$, we then compute the solution $\uN[\param]$ by solving an $\NN \times \NN$ linear system.

\

Once the snapshots $\uN[\param_1],\dots,\uN[\param_P]$ have been computed, we need to select a suitable $M$ dimensional subspace of the space they span. There are different strategies for carrying out such a task, one of which is to perform a proper orthogonal decomposition (POD, \cite{hesthaven2016certified}). 
Assuming that each $\uN[\param_k]$ takes the form of a column vector of length $\NN$,  we assemble the $\NN \times P$ matrix $\Umatr$ containing all the snapshots:
\begin{equation}\label{eq:snapshots_matrix}
	\Umatr = \begin{pmatrix}
		\uN[\param_1] &\vline\quad \dots \quad\vline& \uN[\param_{P}]
	\end{pmatrix}.
\end{equation}
We then construct the correlation matrix $\Cmatr = P^{-1} \Umatr^\top \rosso{\stiff} \Umatr$, with $\stiff$ denoting the stiffness matrix for a scalar product in $\VN$, and compute its eigenvalues and eigenvectors $(\lambda_\ell,\mathbf{z}_\ell)$, $\ell = 1,\cdots, P$, which we assume to be ordered in such a way that the sequence $\{\lambda_\ell\}$ is non increasing. Setting, for $\ell \leq P$, 
\rosso{\begin{equation}
	\funrb_\ell = \frac{1}{\sqrt{P}}\sum_{k=1}^{P} z_\ell^k \uN(\param_k), 
\end{equation}
where $z_\ell^k$ denotes the $k^{\text{th}}$ entry of the eigenvector $\mathbf{z}_\ell$ (that is $\mathbf{z}_\ell = (z_\ell^k)_k$),
}
 we  obtain a new ordered basis $\{\funrb_1,\dots,\funrb_P\}$ for the span of the snapshots. The reduced basis is then obtained by truncating the new basis to the first $M$ elements:
\[
\rbspace = \Span\{ \funrb_\ell, \ \ell=1,\cdots, M\}.
\]

\TESI{\begin{algorithm}
	\caption{Proper orthogonal decomposition}
	\begin{flushleft}
		Snapshots: $\uN(\param_i)$ for $i=1,\dots,P$\\
		Set $\Umatr = \begin{bmatrix}
			\uN(\param_1) &\vline\quad \dots \quad\vline& \uN(\param_{P}) 
		\end{bmatrix}$\\
		\vspace{2mm}
		Build $\Cmatr = \frac{1}{P} \Umatr^\top \Umatr$\\
		\vspace{2mm}
		Solve the eigenproblem $\Cmatr \mathbf{z}_n = \lambda_n \mathbf{z}_n$ for $n=1,\dots,P$\\
		Consider the eigenvectors associated with the $M$ largest eigenvalues\\
		\vspace{2mm}
		Compute $\funrb_i = \frac{1}{\sqrt{P}}\sum_{m=1}^{P} {z}^i_m \uN(\param_m) \quad\text{for }i=1,\dots,M$
	\end{flushleft}
	\label{alg:pod}
\end{algorithm}}


Once the  $\xi_\ell$, $\ell = 1,\cdots,M$ are selected,  the building blocks
\begin{equation}
	\Amatr^{q} = (\aRB^q(\funrb_{\ell'},\funrb_\ell))_{\ell,\ell'},\quad \text{and} \quad \rosso{\Fmatr^{q} =(F^q(\funrb_\ell))_\ell,}
\end{equation}
 for the affine decomposition of the forms $A$ and $F$ are precomputed once and for all, and stored.\TESI{ The steps of the offline procedure are reported in Algorithm~\ref{alg:offline}.
 }
\

	\section{Computation of the virtual basis functions by reduced basis method}\label{sec:vemrb}

The idea  is now to resort to the reduced basis method for solving, for elements $K$ in the mesh $\mesh$, the elemental Laplace equation that defines the virtual functions in the space $\vemk$.  We then fix the number of polygon edges $N$, and, as described in Section \ref{sec:thebasis}, we write the collection of Laplace problems on $N$-edges polygons as a single parametrized PDE of the form  \eqref{eq:postprocessinghat}, set on the reference $N$ edges regular polygon $\refelementVEM$. We start by introducing a fine mesh $\meshfine$ for  $\refelementVEM$, and we let $\VN$ denote the corresponding finite element space
\[
\VN = \{ \hu_\delta \in H^1(\refelementVEM): \ \hu_h|_\tau \in \mathbb{P}_1(\tau), \ \forall \tau \in \meshfine \}.
\] 

	
	 In the following, it will be convenient to reformulate  equation \eqref{eq:postprocessinghat} as a problem with homogeneous boundary conditions. To this aim we 
	  introduce the discrete harmonic lifting $\hTheta_j \in \VN$ of $\widehat g_j$, defined as
	\begin{equation}\label{eq:defTheta}
	\int_{\refelementVEM} \nabla \hTheta_j  \cdot \nabla v_\delta\,d\xhat = 0, \text{ for all }v_\delta \in \VN \cap H^1_0(\refelementVEM), \qquad \hTheta_j = \widehat{g}_j, \text{ on }\partial\refelementVEM.
	\end{equation}
The functions $\widehat e^K_j$ can be written as $\widehat e^K_j = \hTheta_j + \hatxiK{j,\delta}$
 where $\hatxiK{j,\delta} \in H^1_0(\refelementVEM)$ is solution to the following parametrized problem:
	\begin{problem}\label{pro:refK} For all $j = 1, \cdots, N$, 
		find $\hatxijK \in H^1_0(\refelementVEM)$ such that 
		\begin{equation}
			A(\hatxiK{j,\delta} ,\widehat v; K) = - A( \hTheta_j,\widehat v; K)\qquad \forall \widehat v \in  H^1_0(\refelementVEM).
		\end{equation}
	\end{problem}
	
	\

\subsection{The offline phase: computing the snapshots}\label{sec:Koffline}	For all $K \in S$, where $S$ is a randomly generated collection of polygons in $\parameters$, 
		we need to compute snapshots $\hatxijK$,  $j = 1,\cdots,N$, of solutions of Problem \ref{pro:refK}. As, in the presence of badly shaped polygons, {the mapping $\Bcal_K$ may present strong gradients and, consequently, the bilinear form $A(\cdot,\cdot; K)$ may be ill conditioned}, instead of solving such a problem directly by the finite element method on the mesh $\meshfine$, we rather solve the equivalent Problem \ref{pro:basis1} on the polygon $K$. This is done by resorting to  finite elements on a  shape regular fine mesh $\meshK$, of meshsize $\delta^K$, defined directly on $K$, independently of the fine mesh on $\refelementVEM$. The resulting finite element function $e^K_{j,\dd} \in H^1(K)$ is the pull back of a function $\widetilde e^K_j$ in $H^1(\refelementVEM)$.
	The snapshot $\hatxijK  \in \VN$ is then obtained by interpolation as
	\[
	\hatxijK= I_\delta \widetilde e^K_j -\hTheta_j, 
	\]
	where $I_\delta: C^0(\refelementVEM) \to \VN$ is the standard Lagrangian interpolation operator.

	\
	
	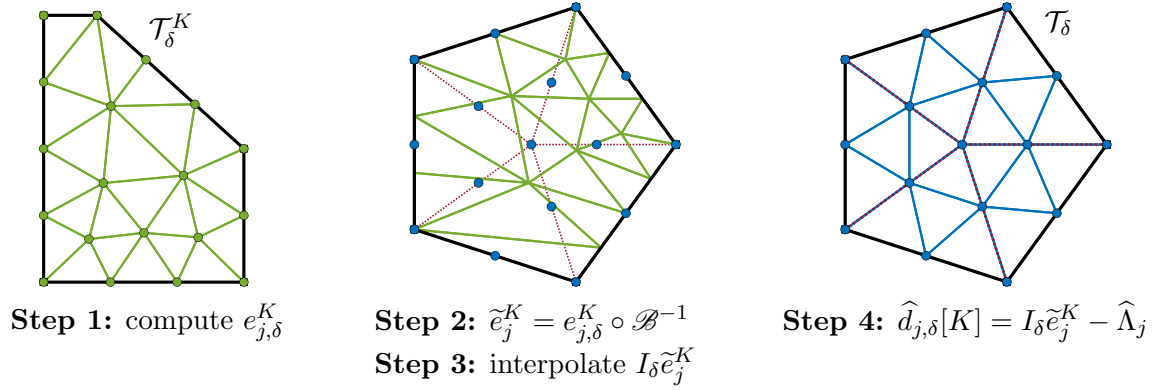
\begin{figure}
		\centering
		\begin{overpic}[width=0.17\linewidth]{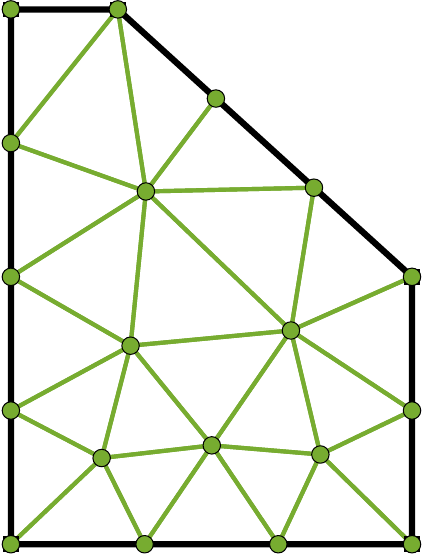}
			\put (40,90) {$\meshK$}
			\put (-10,-15) {\textbf{Step 1:} compute $e^K_{j,\dd}$}
		\end{overpic}\hspace{2cm}
		\begin{overpic}[width=0.22\linewidth]{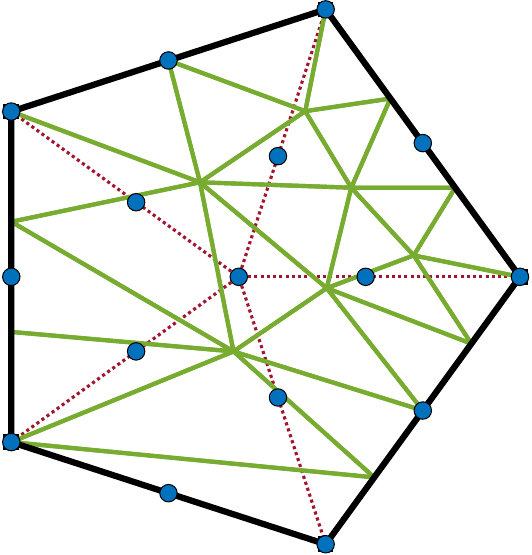}
			\put (-12,-15) {\textbf{Step 2:} $\widetilde e^K_j = e^K_{j,\dd} \circ \Bcal^{-1}$}
			\put (-12,-30) {\textbf{Step 3:} interpolate $I_\delta \widetilde e^K_j$}
		\end{overpic}\hspace{2cm}
		\begin{overpic}[width=0.22\linewidth]{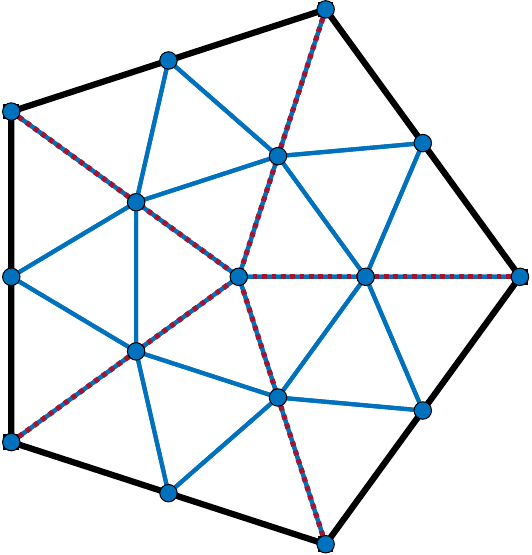}
			\put (72,90) {$\meshfine$}
			\put (-20,-15) {\textbf{Step 4:} $\hatxijK= I_\delta \widetilde e^K_j -\hTheta_j$}
		\end{overpic}
		\vspace{1.3cm}
		\caption{Sketch of snapshots computation. For the sake of clarity we use a depiction with very coarse meshes (computations are of course carried out on much finer meshes).}
		\label{fig:offline_sketch}
	\end{figure}
	
	Once the snapshots are computed, we need to construct the reduced basis functions.  Among other things, we want to leverage the newly computed basis for constructing stabilization terms for non isotropic problems, and we believe that, for such a task, it is important to take into consideration the relation between the different basis functions. Therefore, we propose to perform the POD on the  basis taken as a whole, rather than on the individual basis functions taken one by one. 
	More precisely, letting $K^k$ be the $k^\text{th}$ polygon in $S$ and letting $\widehat d^k_{j,\delta} = \widehat d_{j,\delta}[K^k]$, the $k^\text{th}$ column of the matrix $\Umatr$ defined in \ref{eq:snapshots_matrix} is obtained by stacking on top of each other the finite element coordinates of the snapshots $\widehat d^k_{1,\delta},\dots,\widehat{d}_{N,\delta}^k$  so that each column is associated to the whole basis for the corresponding polygon:
	\[
	\Umatr = \left(
	\begin{array}{ccc}
		\widehat d^1_{1,\delta} & \cdots & \widehat d^P_{1,\delta} \\
		\vdots & & \vdots \\
		\widehat d^1_{N,\delta} & \cdots & \widehat d^P_{N,\delta}
	\end{array}
	\right)
	\]
	(by abuse of notation we use the same symbol for  functions in the finite element space $\VN$ and for the relative vector of degrees of freedom). Applying the POD to this matrix finally results 
	in an ordered sequence of N-tuples $(\xiHat_1^\ell,\cdots , \xiHat_N^\ell)^\top$, $\ell = 1,\cdots,P$, with
	\begin{equation}\label{defhatxi}
	\left( 
	\begin{array}{c}
		\xiHat_1^\ell\\
		\vdots\\
		\xiHat_N^\ell 
	\end{array} 	
	\right) = \frac{1}{\sqrt{P}}\sum_{k=1}^{P} z_\ell^k
	\left( 
	\begin{array}{c}
		\widehat d_{1,\delta}^k\\
		\vdots\\
		\widehat d_{N,\delta}^k 
	\end{array} 	
	\right), 
	\end{equation}
	where $(\lambda_\ell,\mathbf z_\ell)$, $\ell = 1, \cdots, P$, with $\mathbf{z}_\ell = (z_\ell^k)_{k=1}^P$, are the eigenvalue -- eigenvector pairs of the correlation matrix $\Cmatr = P ^{-1} \Umatr^\top\, \stiff\, \Umatr$, ordered in such a way that the sequence $(\lambda_\ell)_\ell$ is not increasing.

	\TESI{The dimension $M$ reduced basis in which we will look for an approximation of the basis function corresponding to the $j$-th node of a new polygon $\elementVEM$ is then 
		defined as
		\[
		\rbspace^j = \Span \{\xiHat_j^1,\cdots,\xiHat_j^M
		\} \subset \VN \cap H^1_0(\refelementVEM).
		\]}

	\subsection{Existence of the affine decomposition} 
In order to implement the online phase efficiently we need to provide an affine decomposition  of the form \eqref{eq:affine} for the bilinear form $A(\cdot,\cdot; K)$ and of the right hand side $A(\hTheta_j,\cdot;K)$. To this aim we start by observing that
\[
A(\hu,\hv;K) = \sum_{i=1}^N A(\hu,\hv;{T}_i)
\qquad \text{ with } A(\hu,\hv;T_i) = \int_{\widehat{T}_i} 
\rosso{| \det(\matrixB_{K,i}^{-1}) |}\,\matrixB_{K,i}^\top\, 
\matrixB_{K,i} \nabla \hu\cdot \nabla \hv\,d\xhat.
\]

We can now expand the symmetric matrix $\rosso{| \det(\matrixB_{K,i}^{-1}) |}\,\matrixB_{K,i} \matrixB_{K,i}^\top$, using a basis for the space of  symmetric $2 \times 2$ matrices, such as the basis $\{\Acal^1,\Acal^2,\Acal^3\}$, where
	\begin{equation}\label{def:symmatbasis}
		\Acal^1= \begin{bmatrix}
				1 & 0 \\
				0 & 0
			\end{bmatrix},\qquad 
			\Acal^2 = \begin{bmatrix}
				0 & 0 \\
				0 & 1
			\end{bmatrix},\qquad 
			\Acal^3 = \begin{bmatrix}
				0 & 1 \\
				1 & 0
			\end{bmatrix},
	\end{equation}
and write
	\begin{equation}\label{eq:matr_dec}
	\rosso{| \det(\matrixB_{K,i}^{-1}) |}\,\Bmatr_{K,i}^\top \,	\Bmatr_{K,i}= \sum_{\nu=1}^3 c_\nu^i[K] \Acal^\nu.
	\end{equation}
We thus have the affine decomposition
\[
A(\hu,\hv;K) = \sum_i \sum_{\nu = 1}^3 c_\nu^i[K] A^{i,\nu}(\hu,\hv), \quad \text{ with }\quad A^{i,\nu}(\hu,\hv) = \int_{\widehat{T}_i}\Acal^\nu \nabla \hu \cdot \nabla \hv\,d\xhat.
\]
Analogously, as far as the right hand side is concerned, we have that
\[
 A( \hTheta_j,\widehat v; K) = \sum_{i=1}^N  \sum_{\nu=1}^3 c^i_\nu[K] F^{i,\nu}_j(\hv) 
 \quad \text{ with }\quad F^{i,\nu}_j(\hv) = \int_{\widehat{T}_i} \Acal^\nu \nabla \hTheta_j\cdot \nabla \hv\,d\xhat.
 \]


 The assembly, in the online phase, of the linear system resulting from solving Problem~\ref{pro:refK} by the Galerkin method in the reduced basis space can then be carried out efficiently, as described in Section \ref{sec:rb}. In the offline phase, once the snapshots are computed and the reduced bases selected, we compute and store the affine decomposition bricks related to the bilinear form and right hand side, as well as some further bricks we will need in the online phase. More precisely, we compute and store the following quantities 
	\begin{equation}\label{defAinu}
		\Amatr_i^\nu(j,j^\prime,\ell,\ell^\prime) = \int_{\That_i} \Acal_\nu \nabla\xiHat_j^\ell\cdot\nabla\xiHat_{j^\prime}^{\ell^\prime}\,d\xhat
		\qquad\text{and}\qquad
		\Fmatr_i^\nu(j,j^\prime,\ell) = \int_{\That_i} \Acal_\nu \nabla\xiHat_j^\ell\cdot\nabla\hTheta_{j^\prime}\,d\xhat.
	\end{equation}
for $\ell, \ell' = 1,\cdots, M$, $j,j' = 1,\cdots,N$, and for $\nu = 1,\cdots,3$.

\begin{algorithm}
	\caption{The offline phase}
	\begin{flushleft}
		\textbf{Initialization:}\\
		\vspace{1mm}
		\begin{itemize}[leftmargin=15pt]
		\item Select $S$: sample of $P$ polygons with $N$ vertices\\
		\item Construct $\refelementVEM$: regular polygon with $N$ vertices\\
		\item Compute $\hTheta_j\in\VN$ according to \eqref{eq:defTheta}\\
		\end{itemize}
		\vspace{2.5mm}
		\textbf{Snapshots computation:}\\
		\vspace{1mm}
		\For{$\elementVEM\in S$}{
			Compute $e^K_{j,\dd}$ solving Problem~\ref{pro:basis1} with FEM on triangulation $\meshK$ of $\elementVEM$\\
			Set $\widetilde e^K_j = e^K_{j,\dd} \circ \Bcal^{-1}$ defined on $\Bcal(\meshK)$ in $\refelementVEM$\\
			Compute $I_\delta \widetilde e^K_j$ interpolating $\widetilde e^K_j$ on $\meshfine$\\
			Compute snapshots $\hatxijK= I_\delta \widetilde e^K_j -\hTheta_j\in\VN$\\
		}
		\vspace{2.5mm}
		\textbf{Proper Orthogonal Decomposition:}\\
		\vspace{1mm}
		Build snapshots matrix $\Umatr$ and correlation matrix $\Cmatr = P ^{-1} \Umatr^\top\, \stiff\, \Umatr$\\
		Compute ordered sequence $(\lambda_\ell,\mathbf z_\ell),$  $\ell=1,\dots,P$, solutions of eigenvalue problem $\Cmatr\mathbf z = \lambda \mathbf z$\\
		Build $(\xiHat_1^\ell,\cdots , \xiHat_N^\ell)^\top$, $\ell = 1,\cdots,P$, according to \eqref{defhatxi}\\
		Compute and store affine decomposition building blocks $\Amatr_i^\nu$, $\Fmatr_i^\nu$
		
	\end{flushleft}
	\label{alg:rb_offline}
\end{algorithm}
\newcommand{\wljK}{w_\ell^{K,j}}
\newcommand{\hejKrb}{\widehat e_{M,j}^{\,\text{rb}}(K)}

	\subsection{Online phase: reconstruction of the basis functions} We can now use the selected reduced basis to contruct (approximations of) the basis functions for the local virtual element space $\vemkloc$. Given a new $N$ vertices polygon $K$, we look for $\hatxijK$ solution to Problem~\ref{pro:refK} in the form:
	\[
	\hatxijK = \sum_{\ell=1}^M \wljK  \xiHat_j^\ell.
	\]

Thanks to the affine decomposition the corresponding linear system is very cheaply assembled, as described in Section \ref{sec:rb}, and, if $M$ is small, it can be cheaply solved. The basis functions $e_j^K$ can then be constructed as the pull back $\widehat e_{M,j}^{\,\text{rb}}[K] \circ \Bcal_\elementVEM$ of the function
	\begin{equation}\label{rbsolution}
\widehat e_{M,j}^{\,\text{rb}}[K]= \hTheta_j +  \sum_{\ell=1}^M \wljK \xiHat_j^\ell \in \VN.
	\end{equation}
	These can be then used to reconstruct the solution and evaluate the desired quantities.

	We would like to point out once more that, depending on our goal, we might not necessarily need a good approximation of the basis functions. If our aim is the design of better stabilization terms for non isotropic problems (see Section \ref{sec:stab}), a very rough approximation might be sufficient. Moreover, as the value of the polynomial part of virtual functions can be computed exactly, similarly to what suggested in \cite{manzini2014new}, we can use the approximated basis functions only to handle the non polynomial part, and this reduces the impact of the error committed in their evaluation. Then, we believe that, for many purposes, very small values of the size $M$ of the reduced basis will be enough.
	
	\begin{algorithm}
	\caption{The online phase}
	\begin{flushleft}
		\textbf{Data:}\\
		\vspace{1mm}
		\begin{itemize}[leftmargin=15pt]
			\item $\hTheta_j\in\VN$: harmonic lifting of $\widehat g_j$\\
			\item $\elementVEM$: polygon with $N$ vertices\\
		\end{itemize}
		\vspace{2.5mm}
		Set $M\le P$ and consider reduced basis $(\xiHat_1^\ell,\cdots , \xiHat_N^\ell)^\top$, $\ell = 1,\cdots,M$\\
		\vspace{2.5mm}
		\textbf{Approximating basis functions for $\vemunoloc$:}\\
		\vspace{1mm}
		\For{j=1,\dots,N}{
		Assemble reduced linear system $\Amatr[\elementVEM]\,\mathsf{w}^{K,j} = \Fmatr[\elementVEM], $ by affine decomposition bricks $\Amatr_i^\nu$, $\Fmatr_i^\nu$\\
		Solve for $\mathsf{w}^{K,j}$ and construct $\widehat e_{M,j}^{\,\text{rb}}[K]= \hTheta_j +  \sum_{\ell=1}^M \wljK \xiHat_j^\ell \in \VN$\\
	}
	\end{flushleft}
	\label{alg:rb_online}
\end{algorithm}

\begin{remark}
	We remark that a straightforward application of the reduced basis method, as described in Section \ref{sec:rb}, would consist looking for a coefficient vector $(w^K_\ell)_\ell$ such that the  functions 
		\[
	\hatxijK = \sum_{\ell=1}^M w^K_\ell  \xiHat_j^\ell,
	\]
where the coefficient $w_\ell^K$ does not depend on $j$, satisfy
			\begin{equation}
	\sum_{j=1}^N	A(\hatxiK{j} , \xiHat_j^\ell; K) = - \sum_{j=1}^N A( \hTheta_j, \xiHat_j^\ell; K)\qquad \forall \ell = 1,\cdots,M.
	\end{equation}
	Here we prefer to build each	$\hatxijK$  function independently of the other, as this does, a priori, likely imply a better approximation of the exact VEM basis functions.
\end{remark}

	\begin{remark}
		We underline that thanks to the affine decomposition, the online phase manages to avoid direct manipulation of elements of the fine space $\VN$. Moreover,  as the computation on each element is completely independent of the other elements, the whole online phase is fully parallelizable. A high degree of parallelization can be easily obtained also for the offline phase, as, on the one hand, the snapshots can be evaluated in parallel, and on the other hand
		the whole offline phases for different values of  $N$ (number of polygon edges) can be performed independently of each other.
	\end{remark}

	\section{Numerical validation}\label{sec:validation}
	In this section, we present some numerical tests, carried out with the aim of assessing the performance of the reduced basis method for reconstructing the virtual element basis functions, both in terms of accuracy and of computational efficiency. As we plan on using our approach in the context of the virtual elements on  Voronoi meshes, where the elements are convex, we will restrict our parameter set by including a convexity condition.

	\subsection{Dataset generation and reduced basis construction}
\begin{figure}
		\includegraphics[trim = 100 10 100 10,width=7.3cm]{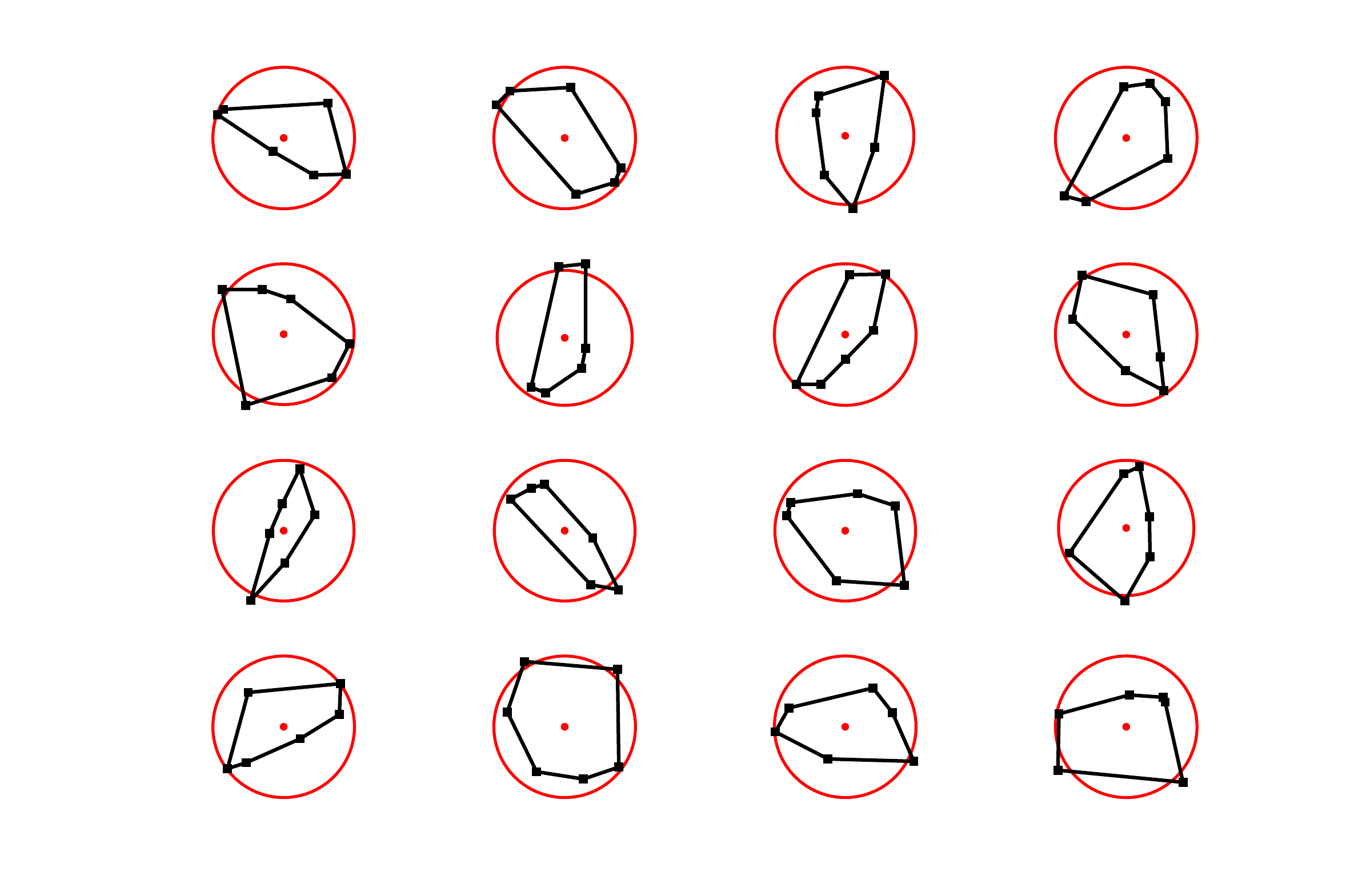}\qquad\vline\qquad
		\includegraphics[trim = 100 10 100 10,width=7.3cm]{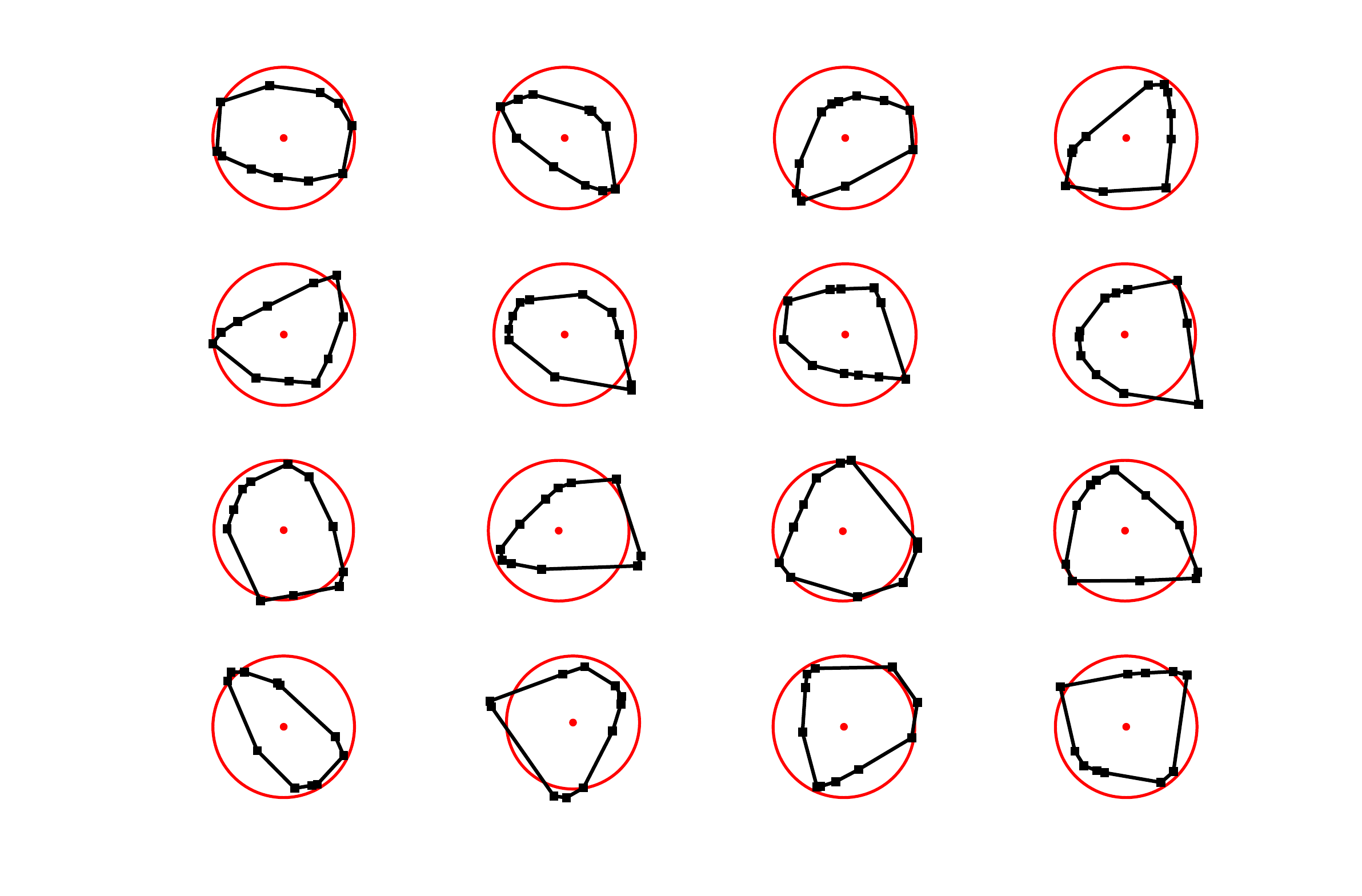}
		\caption{Some random polygons generated for building the reduced basis. On the left  $N=6$, on the right $N=11$.}
		\label{fig:dataset}
	\end{figure}

 For each $N=4, \cdots,14$ we randomly generate a dataset of $5000$ convex polygons with $N$ edges, by using a rejection sampling algorithm, which is   based on \cite{valtr1994probability}, as described by S. Vanderschot in \cite{Vendorschot}. Remark that for  $N=3$ the lowest order VEM basis functions coincide with the order one FEM basis functions on triangles, which are known in their closed form, so that such a case is of no interest in our framework. In Figure~\ref{fig:dataset}, we show, for $N=6$ and $N=11$, some of the randomly generated polygons in the datasets.  
 
 \TESI{Nella tesi aggiungerei una sottosezione che descrive meglio l'algoritmo (in pratica una ripetizione di quello che c'\`e sulla pagina web dalla quale hai preso l'algoritmo). Quello che farei è una Sezione ``Dataset generation'' con due sottosezioni: ``General star shaped polygons'' e ``Convex polygons''.}
 
 	
We then construct a reduced basis: we randomly select $P = 300$ trial polygons $\{\elementVEM^\ell, \ell=1, \cdots,300\}$ out of the dataset, and,
following the procedure described in Section \ref{sec:vemrb}, we  construct and store the first $60$ elements of the ordered sequence of  $N$-tuples $(\xiHat_1^\ell,\cdots , \xiHat_N^\ell)^\top$, $\ell = 1,\cdots,300$, out of which reduced bases of different length $M \leq 60$ can be immediately obtained by simple truncation. Together with the basis functions, we construct and store the different precomputed quantities needed in the online phase.

In our tests the maximum mesh size for both the reference element mesh $\meshfine$ and for the mesh $\meshK$ on the physical polygons that we use in the computation of the snapshots are set to $\dd = \delta^K  = 0.01$.
 In order to map onto $\refelementVEM$ the bases computed on the physical element $K$, and to compute the corresponding function in $\mathcal{V}_\delta$, we implement an interpolation rule based on barycentric coordinates on the triangular elements of the mesh $\meshK$ (see Figure \ref{fig:offline_sketch}).
 
 For different values of $M$, the reduced basis of size $M$ is tested over a set of $500$  polygons, also randomly selected out of the database. The reduced basis thus constructed will also be used later on for performing the tests in Sections \ref{sec:stab} and \ref{sec:postproc}.
 
%

	
\subsection{Accuracy} At first, we test the accuracy of the reconstruction of virtual element functions $u_h \in \vemkloc$. We consider two different situation
\begin{enumerate}[a)]
	\item a function $u_h \in \vemkloc$ obtained by interpolating a smooth function $u$ (in our tests we take $u = x^5 + y^5$);
	\item a function $u_h \in \vemkloc$ whose node values are randomly generated in $(0,1)$ with a normal distribution. 
	\end{enumerate}
As $\PiNabla u_h$ can be computed exactly, in order to reconstruct $u_h$ we only need to resort to the reduced basis method to compute $u_h - \PiNabla u_h$. Then, for each test polygon we carry out the following steps:
\begin{enumerate}[i)]
	\item we compute $\PiNabla u_h$, as usual in the VEM literature;
	\item we compute  $\widehat e_{M,j}^{\,\text{rb}}[K]$ by the reduced basis method with $M$ basis functions;
	\item we compute $\dRB \simeq (I-\PiNabla)u_h$
	 by pulling back a linear combination of the $\widehat e_{M,j}^{\,\text{rb}}[K]$;
	\item we compute $u^\text{rb}_M = \PiNabla u_h + \dRB$;
	\item we compute, for comparison,  the ``exact'' reconstruction $u^\text{fe}_h$, by actually solving the Laplace equation by a finite element method on a triangulation $\meshK$  on the test polygon $K$ of mesh size $\dd$.
\end{enumerate}

\COMMENT{DP: what about moving list i) - iii) to the introduction to better explain what we are going to do in the following pages?}

\begin{figure}
	\includegraphics[width=8cm]{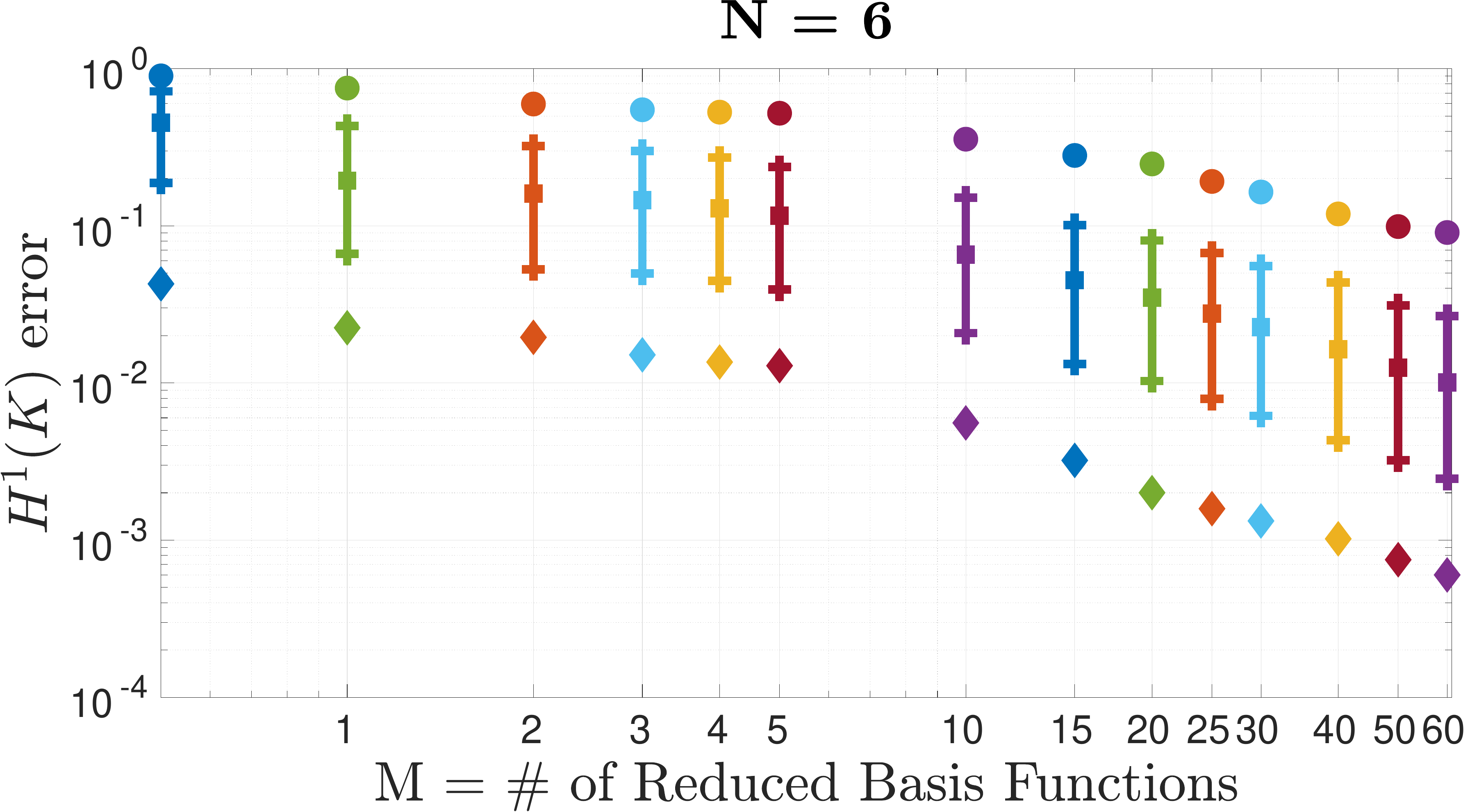}\quad
	\includegraphics[width=8cm]{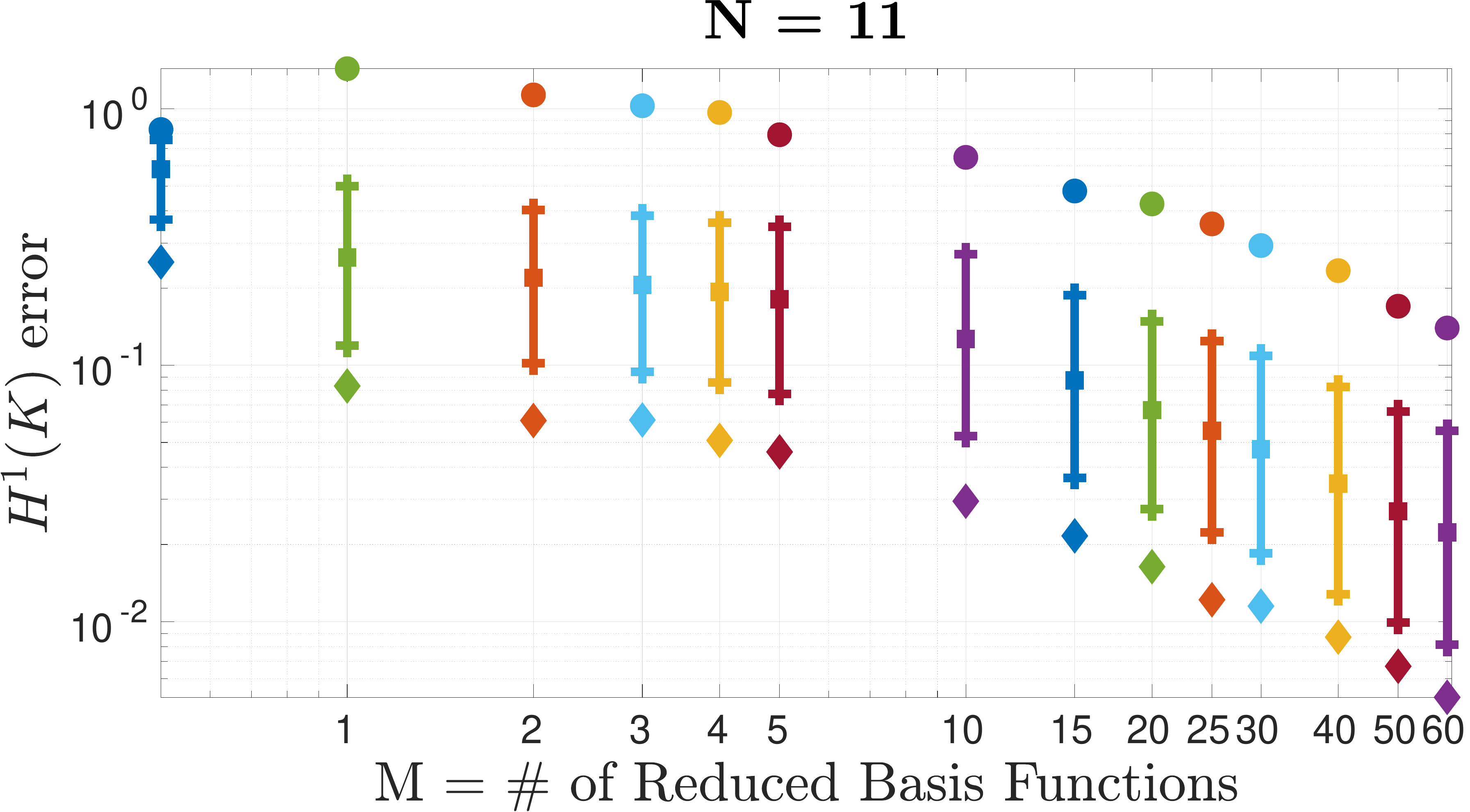}
	\caption{Statistical plots of the errors for $N=6,11$ varying the number of reduced basis $M$.  Test case a) (the degrees of freedom are imposed evaluating $p(x)=x^5+y^5$ at the vertices). The first data, in correspondence with $M=0$, refers to the error between $u_h^\text{fe}$ and $\PiNabla u_h$. Circles represent the maximum values, whereas diamonds represent the minimum values. The averages are depicted by a square. The vertical lines are drawn by connecting the 95\textsuperscript{th} and 5\textsuperscript{th} percentiles.}
	\label{fig:stat_poly}
\end{figure}

\begin{figure}
	\includegraphics[width=8cm]{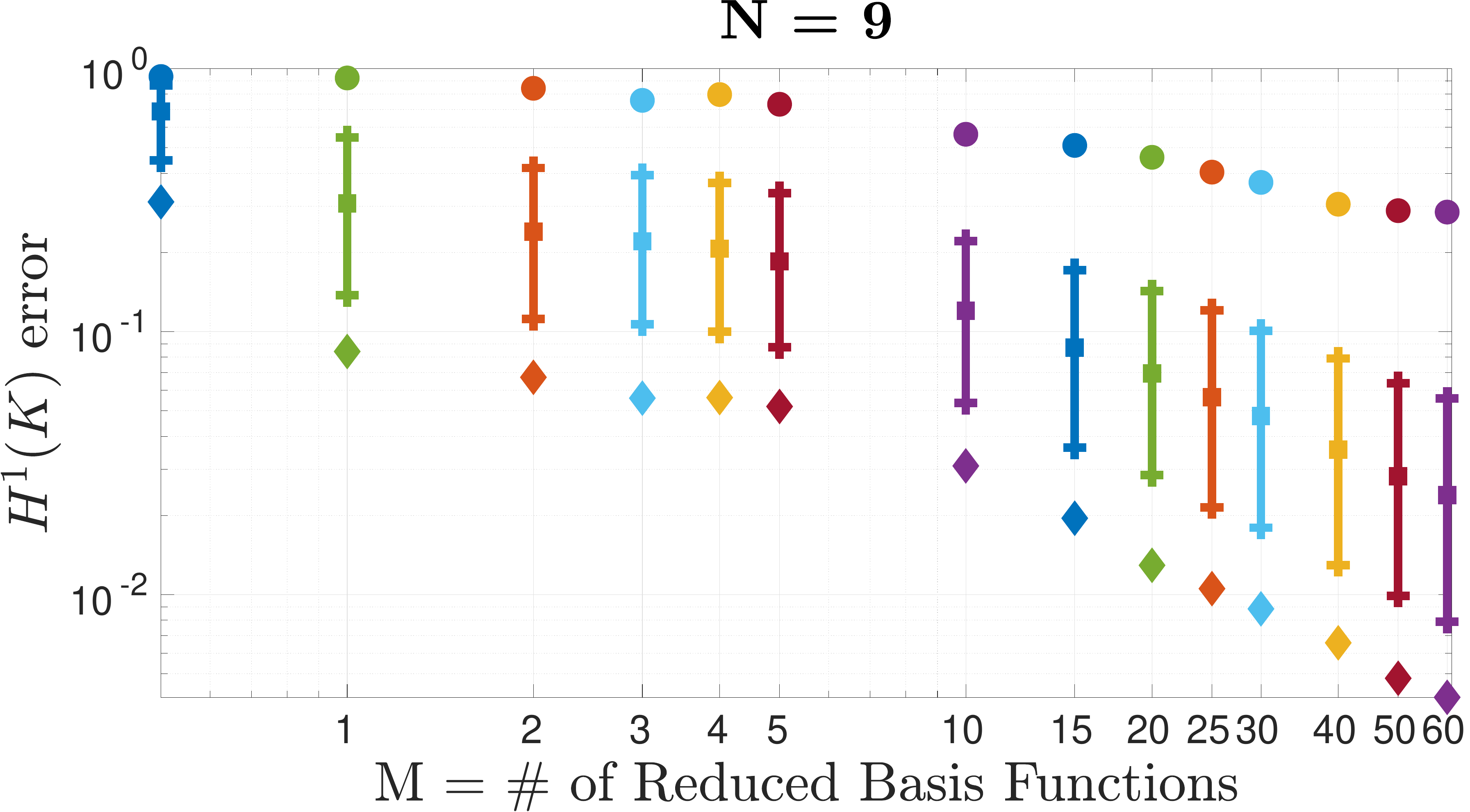}\quad
	\includegraphics[width=8cm]{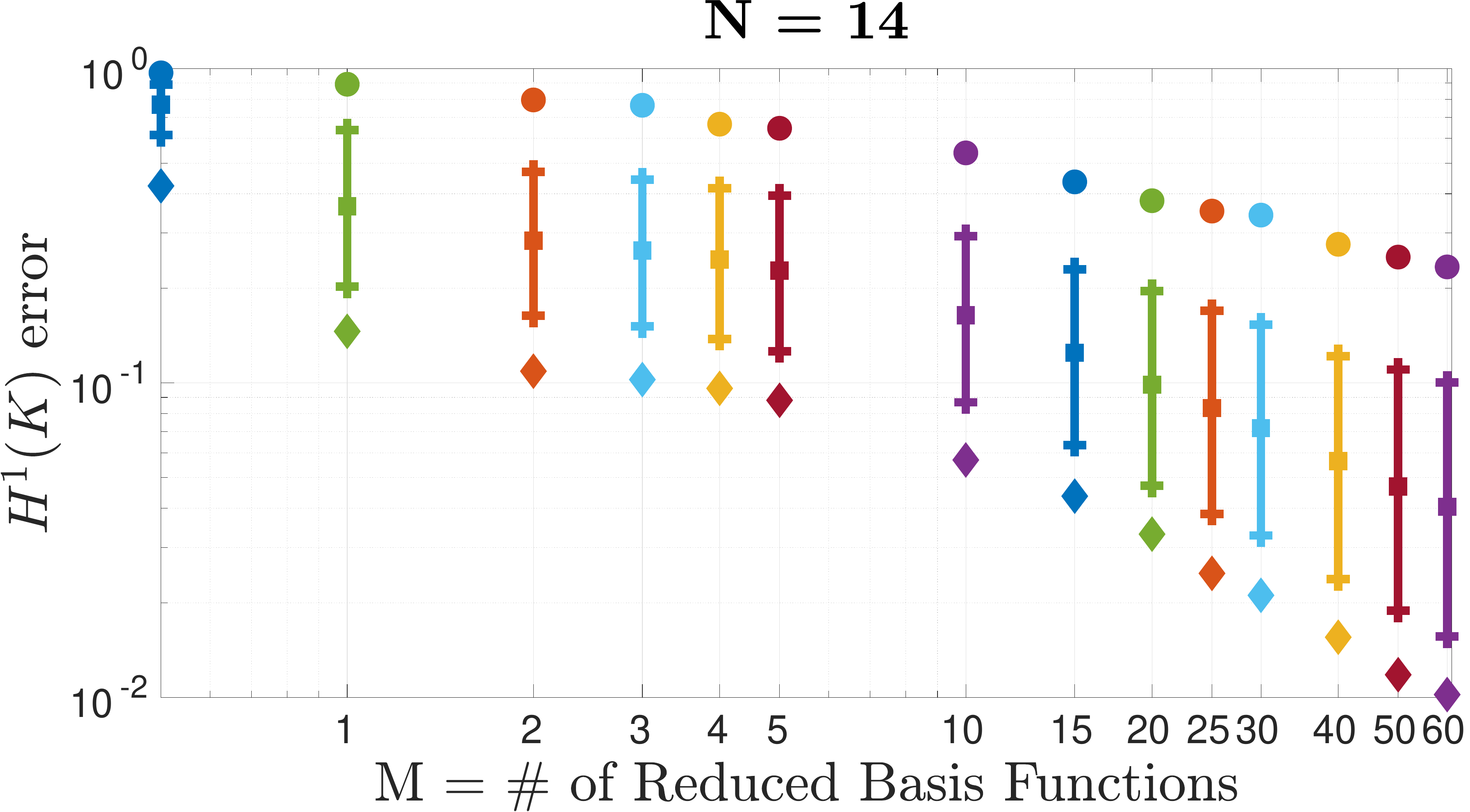}
	\caption{Statistical plots of the errors for $N=9,14$ varying the number of reduced basis $M$. Test case b) (the degrees of freedom are randomly generated in $(0,1)$). Same format as in Figure~\ref{fig:stat_poly}.}
	\label{fig:stat_rand}
\end{figure}

To asses the accuracy of the method, we  compute \rosso{$\norm{u^\text{fe}_h-\uhRB}_{1,\elementVEM}/\norm{u^\text{fe}_h}_{1,\elementVEM}$} for different values of the number $M \geq 1$ of reduced basis functions.
\rosso{In Figure \ref{fig:stat_poly}} we report statistical plots for $N=6,11$ case a), while in Figure \ref{fig:stat_rand}, we report similar plots for $N=9,14$ case b). 
For each $M$ we depict the maximum  (circle) and  minimum (diamond) value. The 5\textsuperscript{th} and the 95\textsuperscript{th} percentiles are connected by a straight line on which the average value is marked by a square.
We can see that in the 95\% of cases, the use of just one reduced basis function already improves the results obtained by projecting onto polynomials. If we look at test case~a), where the function $u_h$  is smoother (Figure \ref{fig:stat_poly}), we see that there are instances in which $\PiNabla u_h$ approximates $u^\text{fe}_h$ better that the reduced basis reconstruction $u^\text{rb}_h$. However, also in such a case, the error is of the same order of magnitude.  
On the other hand, for test case~b), with just two reduced basis functions ($M=2$) we already obtain  \rosso{some} improvement also in the worst case scenario.  This makes sense since for non smooth functions the non polynomial part plays a relevant role. If we look at the best case scenario, one order of magnitude is gained with just $M=1$. \rosso{Finally, looking at the best case scenario for $M=60$, we notice an improvement of almost two orders of magnitude with respect to the best case scenario obtained by $\PiNabla u_h$. This phenomenon is observed in both case a) and case b).}

\begin{figure}
	\centering
	\vspace{2.5mm}
	\includegraphics[trim = 40 40 40 40, width=0.23\linewidth]{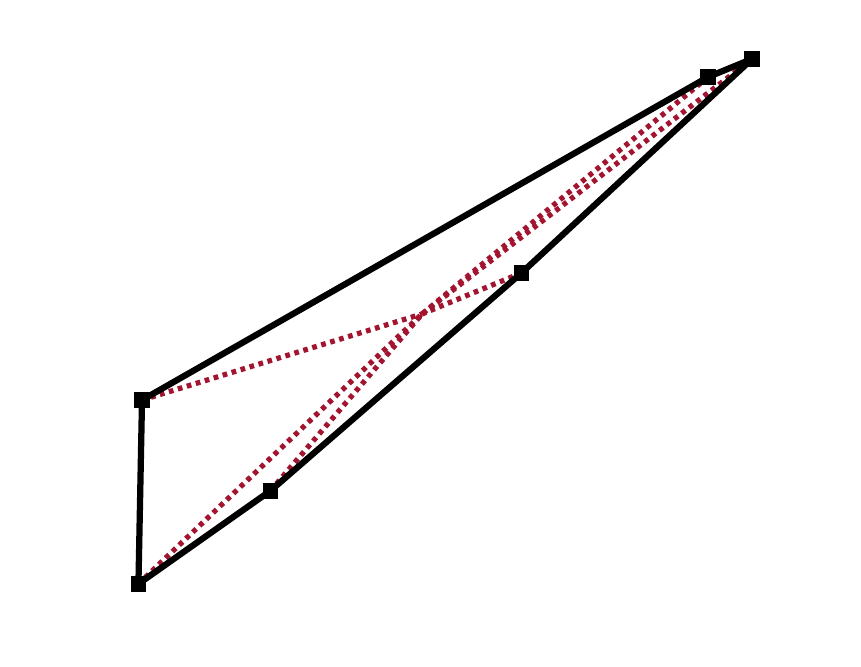}\quad
	\includegraphics[trim = 40 40 40 40, width=0.23\linewidth]{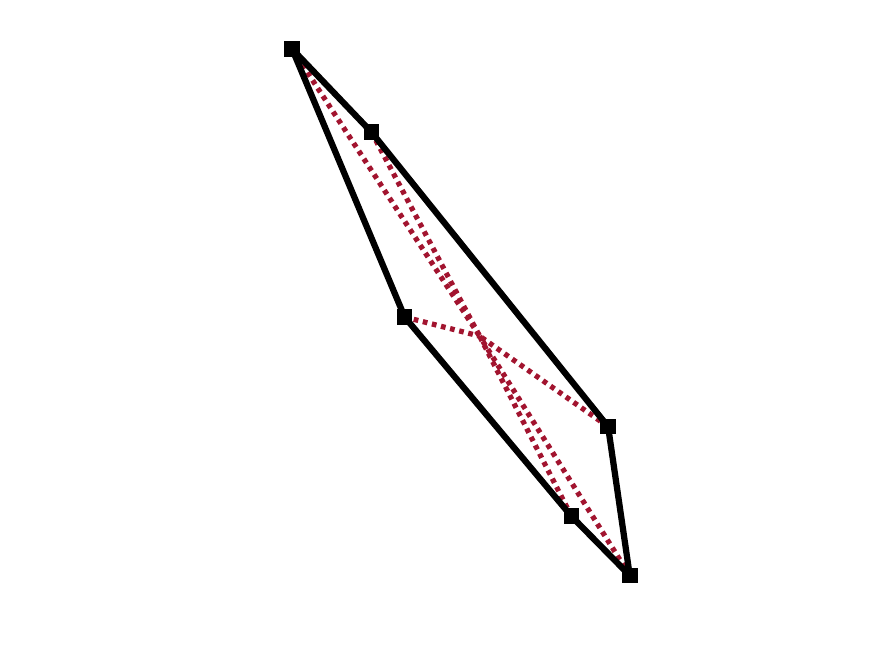}\quad
	\includegraphics[trim = 40 40 40 40, width=0.23\linewidth]{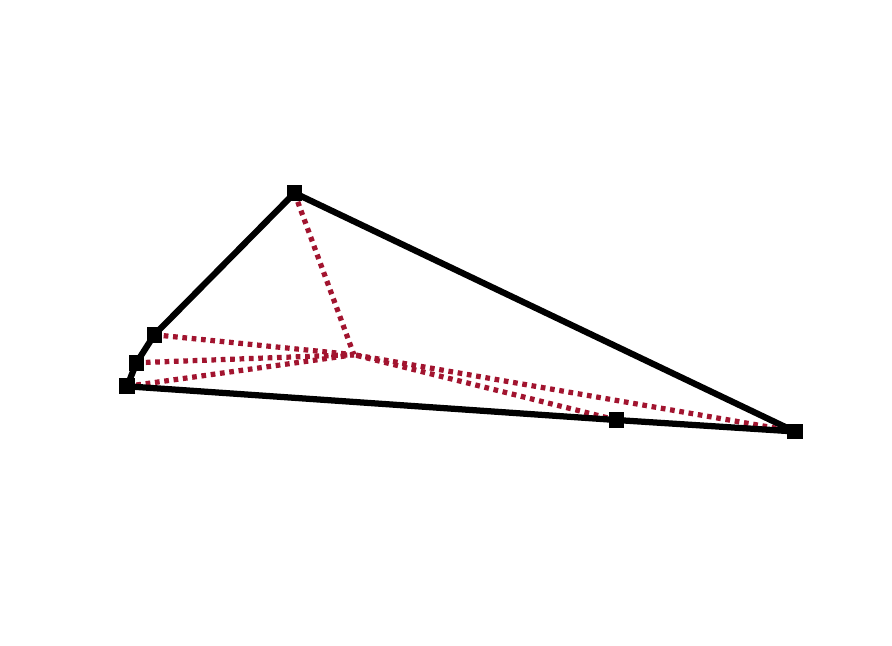}\quad
	\includegraphics[trim = 40 40 40 40, width=0.23\linewidth]{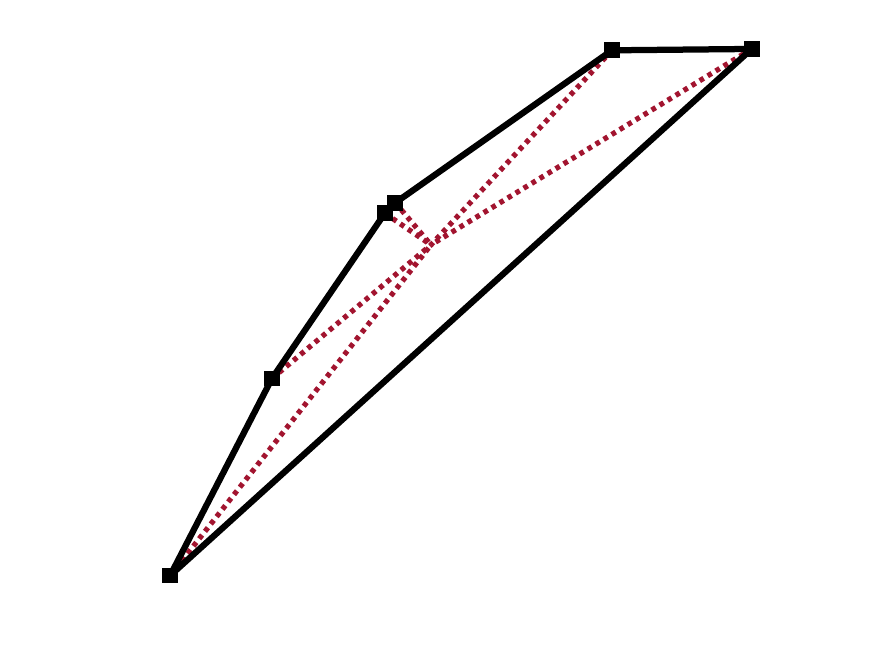}\\
	\vspace{7mm}
	\includegraphics[trim = 40 40 40 40, width=0.23\linewidth]{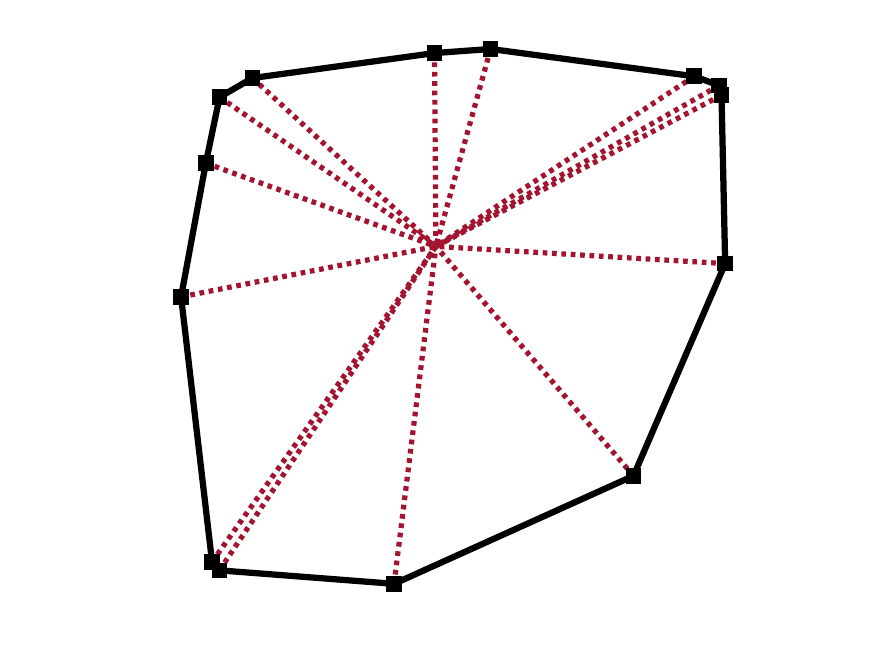}\quad
	\includegraphics[trim = 40 40 40 40, width=0.23\linewidth]{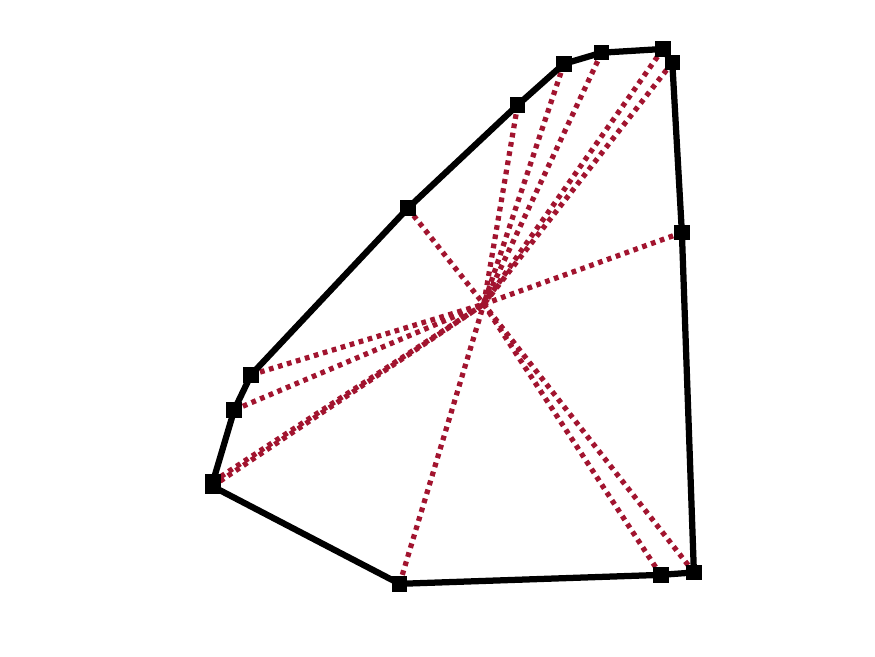}\quad
	\includegraphics[trim = 40 40 40 40, width=0.23\linewidth]{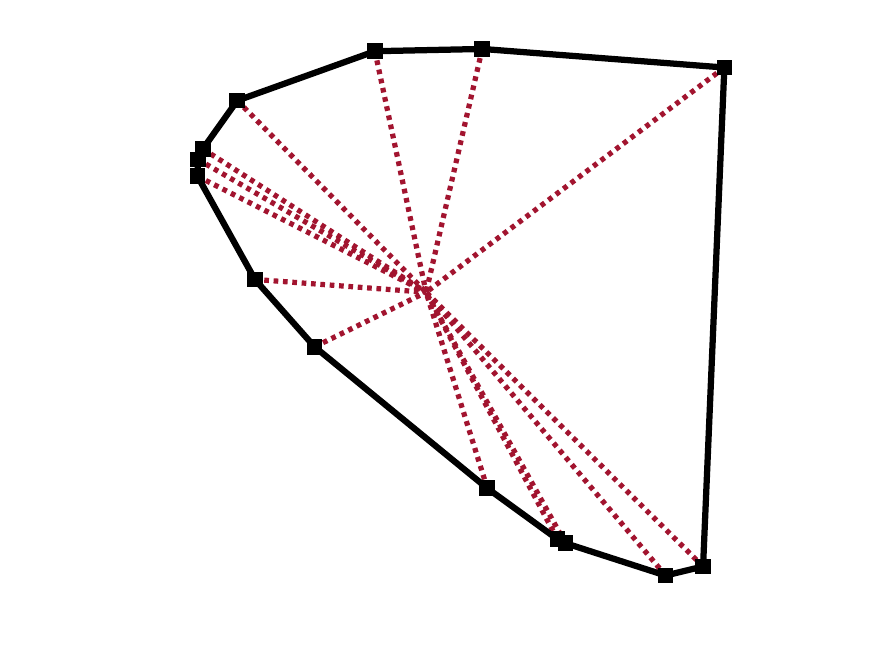}\quad
	\includegraphics[trim = 40 40 40 40, width=0.23\linewidth]{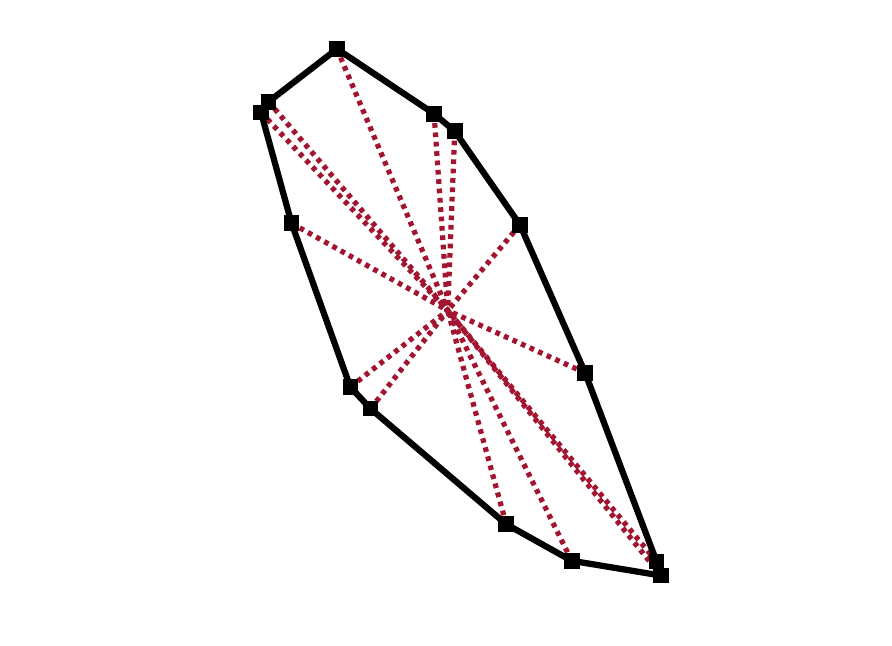}
	\vspace{2mm}
	\caption{Some polygons realizing the maximum values in Figures~\ref{fig:stat_poly}~and~\ref{fig:stat_rand} for $N=6$ (first line) and $N=14$ (second line). It is clear that some triangles of the decomposition are extremely badly shaped.}
	\label{fig:bad_polygons}
\end{figure}

Some polygons realizing the maximum values of the error in Figures~\ref{fig:stat_poly}~and~\ref{fig:stat_rand} for $N=6,14$ are visualized in Figure~\ref{fig:bad_polygons}. 
In all cases we notice some extremely badly shaped triangles  $T_i$, which we believe could be at the source of the issue.
 Further studies will be focused on finding \emph{a priori} criteria to single out polygons that need special treatment to reduce the error.


%


\subsection{Computational efficiency} To assess the computational efficiency of the method, we compare the computational cost of the evaluation of $u^\text{rb}_h$ for different values of $M$, with both the cost of the evaluation of the finite element approximation $u^\text{fe}_h$ (this, of course, is not a practically viable method, and is reported only for the sake of comparison), and the cost  of the evaluation of $\PiNabla u_h$.

In Table~\ref{tab:times1} we report the CPU times in seconds for $N=6,9,11,14$, for the evaluation of $\PiNabla u_h$, $u^\text{fe}_h$  and \rosso{$\uhRB$ for $M=1,5,30,60$}. For the evaluation of  $\PiNabla u_h$, we measure the time $T_{build}$ needed to build the associated matrix \cite{beirao2014hitchhiker} and the time $T_{apply}$ needed to evaluate $\PiNabla u_h$ in the triangulation nodes of each polygon. For the finite element evaluation $u^\text{fe}_h$, we measure the time $T_{assemble}$ required to generate the triangulation and assemble the finite element stiffness matrix and $T_{solve}$ to solve the linear system associated to the post--processing problem. For the reduced basis approximations, $T_{assemble}$ and $T_{solve}$ represent the time required to assemble and solve respectively the linear system in the online phase.  We observe that, though it is increasing as $M$ increases, the time we need to compute $u_h$ by the the reduced basis reconstruction is, for all values of $M$, comparable with the time needed for evaluating $\PiNabla$ and may therefore be considered  an acceptable overhead to the overall cost of a Virtual Element method.

We emphasize once again that the actual
reconstruction of the virtual element functions by the RB method
is highly parallelizable, as each element can be handled completely independently from the others. Moreover, as demonstrated by the numerical tests presented in Table~\ref{tab:times1}, the computation on each element is quite cheap. Indeed, thanks to the affine decomposition, all the relevant quantities can be computed directly from precomputed quantities without the need of constructing or referring to any finite element mesh on the physical or on the reference element.


\renewcommand{\arraystretch}{1.4}
\begin{table}[h]
	\begin{center}
		\begin{tabular}{r|rr|rr|rr}
			\multicolumn{7}{c}{\textbf{Comparison in time of post--processing techniques}}\\
			\hline
			& \multicolumn{2}{c|}{$\PiNabla u_h$} & \multicolumn{2}{c|}{$u^\text{fe}_h$} & \multicolumn{2}{c}{$\uhRB$, $M=60$}\\
			\hline
			$N$ & $T_{build}(s)$& $T_{apply}(s)$ & $T_{assemble}(s)$ & $T_{solve}(s)$& $T_{assemble}(s)$ & $T_{solve}(s)$\\
			\hline
			6 & 6.34e-4 & 1.04e-3 & 2.58e-1 & 3.08e-1 &
			 1.30e-3 & 4.26e-4\\
			9 & 7.22e-4 & 1.46e-3 & 4.27e-1 & 5.97e-1 &
			 2.77e-3 & 6.02e-4\\
			11& 7.72e-4 & 1.80e-3 & 4.42e-1 & 6.32e-1 &
			3.94e-3 & 6.87e-4\\
			14& 7.69e-4 & 2.31e-3 & 5.72e-1 & 9.07e-1 &
			6.04e-3 & 8.13e-3\\
			\hline
			\hline
			& \multicolumn{2}{c|}{$\uhRB$, $M=1$}& \multicolumn{2}{c|}{$\uhRB$, $M=5$}& \multicolumn{2}{c}{$\uhRB$, $M=30$}\\
			\hline
			$N$ &  $T_{assemble}(s)$ & $T_{solve}(s)$& $T_{assemble}(s)$ & $T_{solve}(s)$& $T_{assemble}(s)$ & $T_{solve}(s)$\\
			\hline
			6 & 5.00e-4	 & 3.25e-5  & 5.47e-4 & 1.22e-4 & 8.10e-4 & 2.39e-4 \\
			9 &	9.78e-4  & 3.84e-5  & 1.09e-3 & 1.55e-4 & 1.72e-3 & 3.20e-4 \\
			11&	1.28e-3  & 3.79e-5  & 1.45e-3 & 1.68e-4 & 2.35e-3 & 3.58e-4 \\
			14&	1.92e-3  & 4.11e-5  & 2.18e-3 & 1.71e-4 & 3.60e-3 & 3.92e-4 \\
			\hline
		\end{tabular}
	\end{center}
	\caption{\rosso{Average CPU times required to evaluate $\PiNabla u_h$, $u^\text{fe}_h$ and $\uhRB$ for $M=1,5,30,60$} on $500$ test polygons. The average is computed by considering each polygon twice: indeed, we take into account both \rosso{cases a) and b)} for imposing the degrees of freedom. $T_{build}$ = CPU time to build the projection matrix associated to $\PiNabla$; $T_{apply}$ = CPU time to evaluate $\PiNabla u_h$ on mesh nodes. For the exact reconstruction of $u^\text{fe}_h$: $T_{assemble}$ = CPU time to assemble the FEM linear system; $T_{solve}$ = CPU time to solve it. For the reduced basis approximations: $T_{assemble}$ = CPU time to assemble the linear system of the online phase; $T_{solve}$ = CPU time to solve it. the linear systems are solved with the $\backslash$ command provided by Matlab. \rosso{The code was ran on an Intel Xeon Gold 6230R core running at 2.10GHz.}}
	\label{tab:times1}
\end{table}
\renewcommand{\arraystretch}{1}

\TESI{Nella tesi mettere i test per tutti gli N}

	

\newcommand{\haloc}{\widehat a^K}

\section{Application I: Stabilization}\label{sec:stab}

As mentioned in Section \ref{sec:intro}, we can exploit the possibility of cheaply constructing approximations to the basis functions of the local virtual element space, to design of ad hoc local stabilization bilinear forms for strongly anisotropic problems.

The idea is to define $\shloc$ as
\[
\shloc(
e_i,e_j
) = \haloc(\widehat e_{M,i}^{\,\text{rb}}[K],\widehat e_{M,j}^{\,\text{rb}}[K]), 
\]
where
\[
\haloc(u,w) = \sum_{i=1}^{N} \int_{\That_i} \rosso{| \det(\Bmatr^{-1}_{\elementVEM,i})|}\, \Bmatr_{\elementVEM,i}^\top\, \aa\, \Bmatr_{\elementVEM,i}\grad\widehat{u}\cdot\grad\widehat{v}\,d\xhat
\]
is the bilinear form on the reference element obtained by change of variable from the bilinear form~$\aloc$.

As far as the choice of the size $M$ of the reduced basis is concerned, we observe that 
taking a large value of $M$ would ideally result in accurately reconstructing the local shape functions and, consequently, in taking $\shloc \sim \aloc$. In other words, for $M$ large, rather than a virtual element method we would have a polygonal finite element method, based on the virtual elements discrete space, which, in the lowest order case considered here, coincides with the polygonal finite element space with harmonic barycentric coordinates \cite{martin2008polyhedral,DeRose2006HarmonicC}. 
In the spirit of the virtual element method, we rather take $M$ small/very small, so that $\shloc$ only roughly approximates $\aloc$. 
We observe that $\shloc(
e_i,e_j
)$ can be computed efficiently thanks to the affine decomposition. Indeed, we have
\begin{equation}
	\shloc(
	e_i,e_j
	)  =
	\haloc(  \hTheta_i, 
	\hTheta_j 
	) + \sum_{\ell=1}^M w^{K,i}_\ell
	\haloc(   \xiHat_i^\ell, 
	\hTheta_j
	) +  \sum_{\ell=1}^M \wljK
	\haloc(  \hTheta_i , 
	\xiHat_j^\ell
	)
	+  \sum_{\ell,\ell'=1}^M 
	\wljK
	\haloc(  \xiHat_i^{\ell'}, 
	\xiHat_j^\ell
	),\label{affinestab}
\end{equation}
where $\wljK$ are the coefficient in the expansion \eqref{rbsolution}.
\blu{As we already did in \eqref{eq:matr_dec},} the matrix 
$\rosso{| \det(\Bmatr^{-1}_{\elementVEM,i})|}\, \Bmatr_{\elementVEM,i}^\top\, \Kmatr\, \Bmatr_{\elementVEM,i}
$ can be easily decomposed as 
\[
\rosso{| \det(\Bmatr^{-1}_{\elementVEM,i})|}\, \Bmatr_{\elementVEM,i}^\top\, \aa\, \Bmatr_{\elementVEM,i} = \sum_{\nu = 1}^4 \coef_\nu^i[K] \Acal^\nu,
\]
where, for $\nu = 1,2,3$, $\Acal^\nu$ is defined in \eqref{def:symmatbasis} and 
\[
\Acal^4= \begin{bmatrix}
	0 & 1 \\
	-1 & 0
\end{bmatrix},
\]
so that we have
\[
\haloc(\xiHat_{j}^\ell,\xiHat_{j^\prime}^{\ell'}) = \sum_{i=1}^N \sum_{\nu = 1}^4 \coef_\nu^i[K] \Amatr_i^\nu(j,{j^\prime},\ell,\ell'). 
\]
Then, having precomputed $\Amatr^{\nu}_i$, $i=1,\cdots,N$, $\nu = 1,\cdots,4$, the right hand side of \eqref{affinestab} is easily  assembled.

%

\begin{algorithm}
	\caption{Reduced basis stabilization in $\vemunoloc$}
	\begin{flushleft}
		\textbf{Data:}\\
		\vspace{1mm}
		\begin{itemize}[leftmargin=15pt]
			\item $K$: element of $\mesh$\\
		\end{itemize}
		\vspace{2.5mm}
		
		\textbf{Compute stiffness matrix with RB stabilization:}\\{
			\vspace{2.5mm}	
			Project virtual basis functions onto $\PP_1(\elementVEM)$: $\PiNabla e_j$, $j=1,\dots,N$\\
			Build $\mathsf{\Pi}$ such that $\mathsf{\Pi}_{i,j}=\aloc(\PiNabla e_i,\PiNabla e_j)$ for $i,j=1,\dots,N$\\
			
			Compute $\mathsf{R}=\mathsf{I}-\mathsf{\Pi}$, i.e. $\mathsf{R}_{i,j} = e_i(\vv_j) - \PiNabla e_i(\vv_j)$ for $i,j=1,\cdots,N$\\
			\vspace{2.5mm}
			
			Go to \textbf{Online phase} (see Algorithm~\ref{alg:rb_online}):\\
			\quad \textbf{Input:} $\elementVEM$\\
			\quad \textbf{Output:} RB approximation of VEM basis functions, $\widehat e_{M,j}^{\,\text{rb}}[K], j=1,\cdots,N$\\
			\vspace{2.5mm}
			
			Compute the affine decomposition coefficients $\coef_\nu^i[K]$, $i=1,\cdots,N$, $\nu = 1,\cdots,4$\\
			Construct $\haloc(\xiHat_{j}^\ell,\xiHat_{j^\prime}^{\ell'}) = \sum_{i=1}^N \sum_{\nu = 1}^4 \coef_\nu^i[K] \Amatr_i^\nu(j,{j^\prime},\ell,\ell')$ by affine decomposition\\
			
			Build approximate stiffness $\mathsf{K}^\text{rb}$, i.e. $\mathsf{K}^\text{rb}_{i,j}=\haloc(\widehat e_{M,i}^{\,\text{rb}}[K],\widehat e_{M,j}^{\,\text{rb}}[K])$ for $i,j=1,\cdots,N$\\
			Compute stabilization term $\mathsf{S}^\text{rb} =\mathsf{R}^\top\,\mathsf{K}^\text{rb}\,\mathsf{R}$\\
			Compute VEM stiffness $\mathsf{K} = \mathsf{\Pi} + \mathsf{S}^\text{rb}$\\
		}
	\end{flushleft}
	\label{alg:rb_stab}
\end{algorithm}

\rev{ The proposed stabilization is of course more expensive than standard diagonal stabilizations such as \emph{dofi--dofi} or \emph{D--recipe}. For the reduced basis stabilization, the bulk of the computational cost of the construction of $S^ K_h$ lies in the evaluation of the coefficients $\wljK$, for which we refer to Table \ref{tab:times1}. Its overall cost is of the same order of magnitude of the cost of the evaluation of the consistency component of the stiffness matrix. Conversely, standard diagonal stabilization terms such as \emph{dofi--dofi} or \emph{D--recipe} are  extremely cheap, their cost being two to three orders of magnitude lower, and therefore, overall negligible. On the other hand, as we will see below, depending on the problem, the reduced basis stabilization might lead to an improvement in the convergence of the method that, in our opinion, makes such computational overhead acceptable. }

\

\begin{figure}
	\centering
	\includegraphics[width=0.325\linewidth]{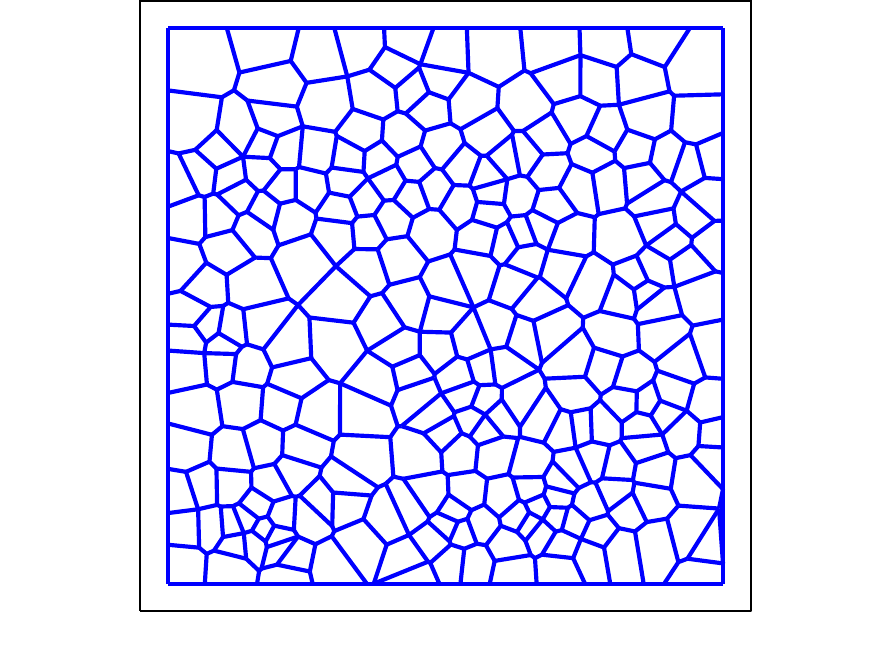}
	\includegraphics[width=0.32\linewidth]{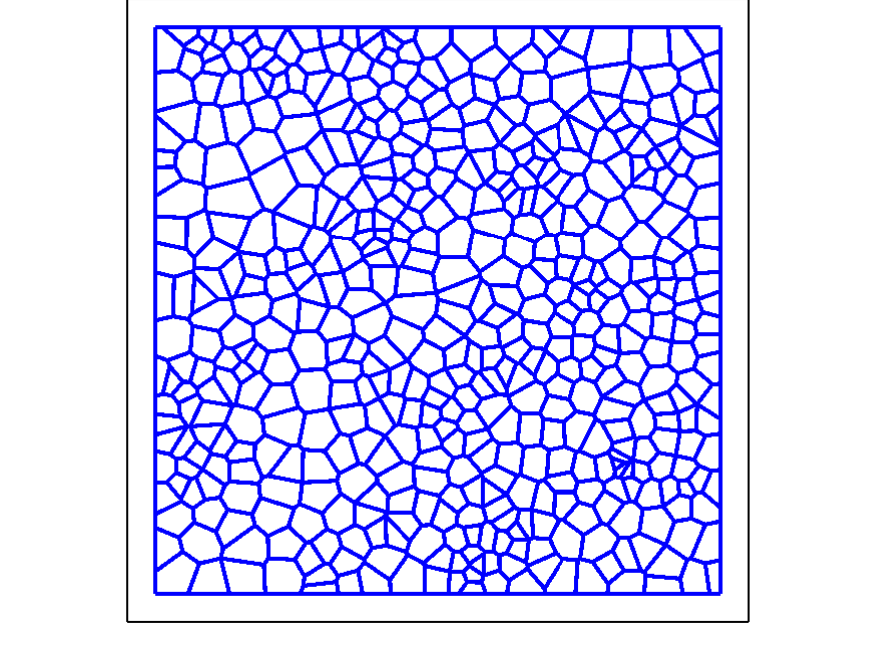}
	\includegraphics[width=0.32\linewidth]{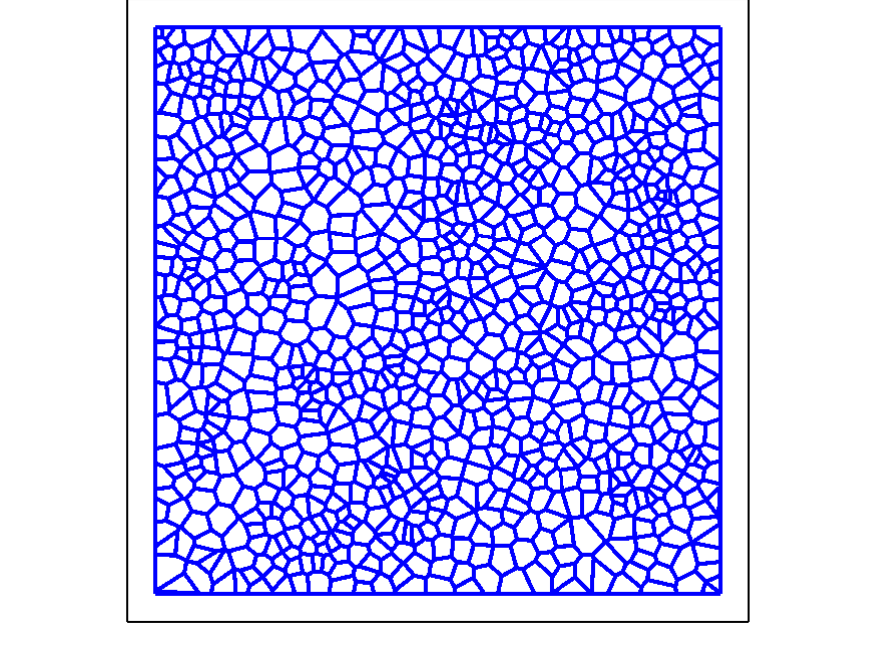}
	\caption{Example of Voronoi meshes for the unit square.}
	\label{fig:mesh_voro}
\end{figure}

\rev{
	
	To assess the performance of the new reduced basis  stabilization, we start by comparing the condition number of the resulting stiffness matrix with the condition number of, on the one hand, the standard VEM stiffness matrix with \emph{dofi--dofi} stabilization, and, on the other hand, the ideal (non computable) stiffness matrix, which we assembled after constructing the virtual basis functions by means of a finite element solution of the local PDE \eqref{truebasisfunctions}.   
	On a unit squared domain, we consider both a Laplace operator, and two anisotropic diffusion problems with diffusivity tensors $\aa_1$ and $\aa_2$ respectively defined as
	\[
	\aa_1 = \begin{bmatrix}
		1 & 0 \\
		0 & 6.25\cdot10^{-4}
	\end{bmatrix}, \qquad \qquad \aa_2 = \begin{bmatrix}
		1 & 10^{-2} \\
		5\cdot10^{-3} & 10^{-4}
	\end{bmatrix}.
	\]
	
	These are the diffusivity tensors that we will use later for the convergence tests. In Figure~\ref{fig:cond}, we plot the evolution of the condition numbers for the different operators discretized on a sequence of Voronoi meshes generated by means of Polymesher~\cite{polymesher}  (see, for instance, Figure~\ref{fig:mesh_voro}). We observe that 
	the condition number of the stiffness matrix obtained through RB stabilization exhibits the same behavior of the condition numbers of both the standard VEM matrix with \emph{dofi--dofi} stabilization, and the ``true'' stiffness.

	\begin{figure}
		\centering
		\textbf{Condition number}
		\subfloat{\includegraphics[width=0.34\linewidth]{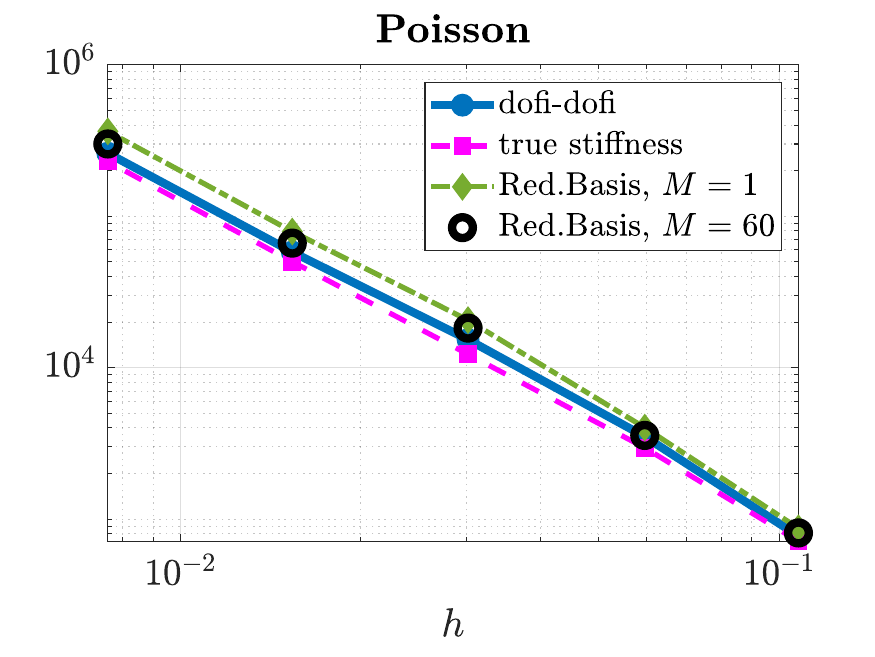}}
		\subfloat{\includegraphics[width=0.34\linewidth]{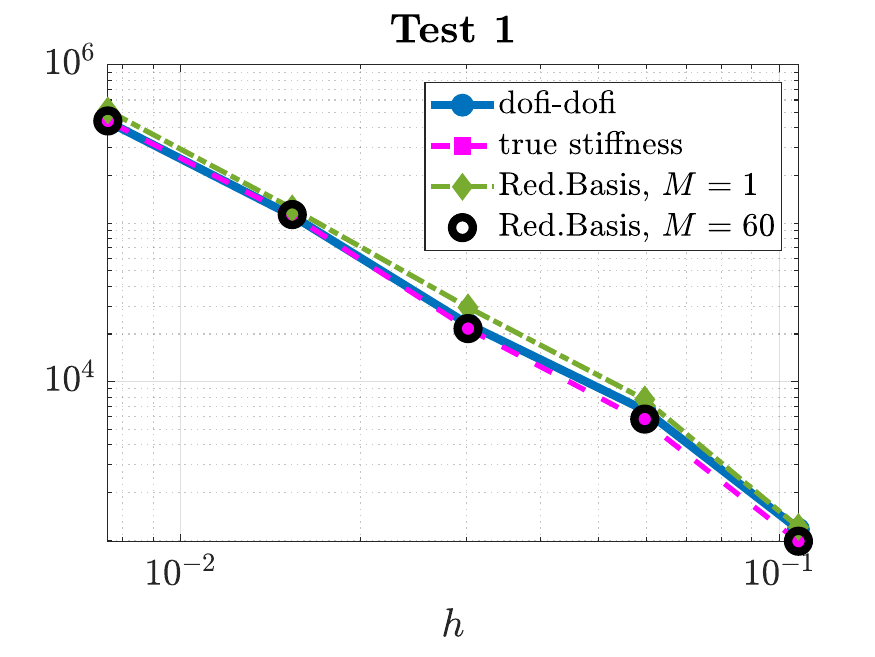}}
		\subfloat{\includegraphics[width=0.34\linewidth]{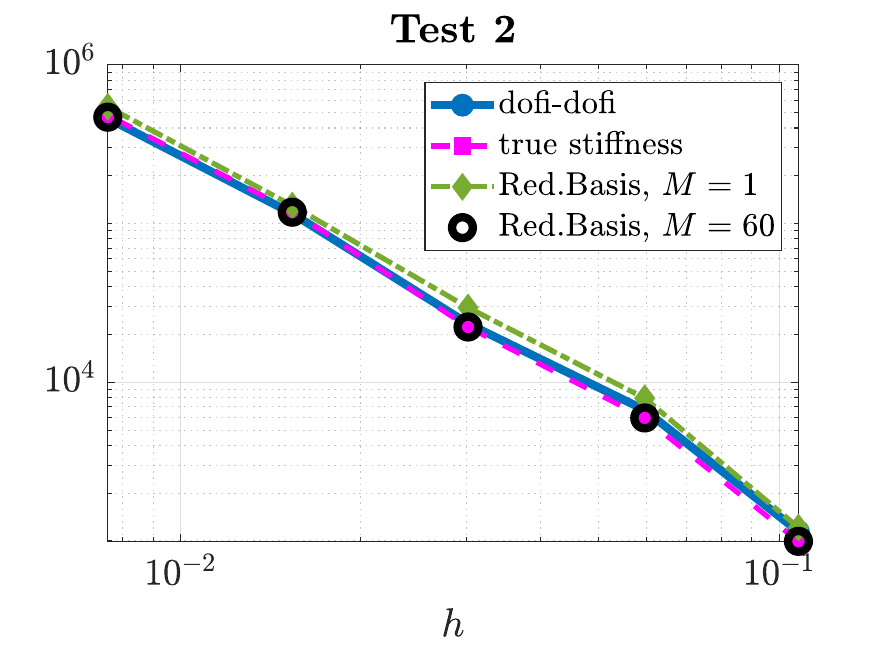}}
		\caption{\rev{From left to right: behavior of the condition number for the Poisson problem, Test~1 ($\aa = \aa_1$) and Test~2 ($\aa=\aa_2$). In the three cases we observe no significant difference between the different versions of the VEM method.}}
		\label{fig:cond}
	\end{figure}

}

\

\rev{
	Focusing on strongly anisotropic problems, where the standard lowest order VEM formulations have been observed to show poor performance \rosso{(see \cite{berrone2022comparison})}, we next compare the convergence properties of the  reduced basis stabilization,  with  the \textit{dofi--dofi} and \textit{D--recipe}, on two test cases. 
	For both, we consider Problem~\ref{pro:cont_poisson} with strongly anisotropic solutions and diffusivity tensors. In both cases the domain $\Omega$ is once again the unit square, discretized by a sequence of Voronoi meshes generated by means of Polymesher~\cite{polymesher}. }
In order to evaluate the performance of the methods, we introduce the following two relative error norms: for $\diamond = u^\text{fe}_h, \PiNabla u_h,\uhRB$ and $\Square = u^\text{fe}_h, \uhRB$, we set
\begin{equation*}
	\begin{aligned}
		&\errstar{\diamond} = \frac{ \bigg( \sum_{\elementVEM\in\mesh} \norm{u-\diamond}^2_{\star,\elementVEM}\bigg)^{\frac12} }{\norm{u}_{\star,\Omega}}&\quad\text{ for }\quad\star=0,1\\
		&\erren{\diamond} = \frac{ \bigg( \sum_{\elementVEM\in\mesh} \norm{ \sqrt{\Kmatr} \grad \big( u-\diamond \big)}^2_{0,\elementVEM}\bigg)^{\frac12} }{\norm{\sqrt{\Kmatr}\grad u}_{0,\Omega}}&\\
		&\errinf{\Square} = \frac{\max_{\x\in\Omega}|u-\Square|}{\max_{\x\in\Omega}|u|}&
	\end{aligned}
\end{equation*}	
As customary for the virtual element method, for the tests in this section, we evaluate $\errstar{\diamond}$ and $\erren{\diamond}$ with  $\diamond = \PiNabla u_h$ (that is, we do not reconstruct the full discrete solution by either the RB method or the finite element method).

\subsubsection*{Test 1.}  
\begin{figure}
	\centering
	\includegraphics[width=0.45\linewidth]{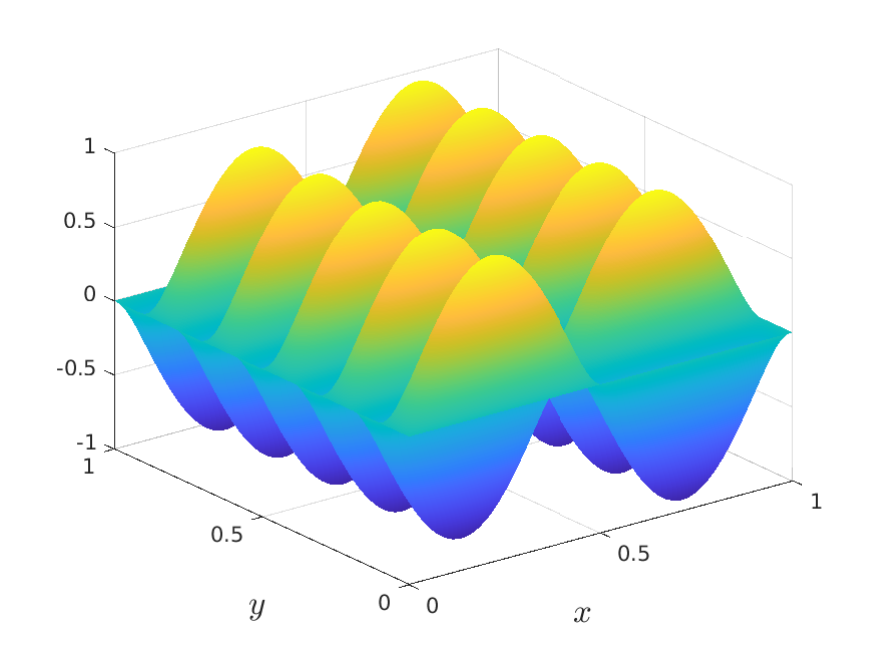}\quad
	\includegraphics[width=0.45\linewidth]{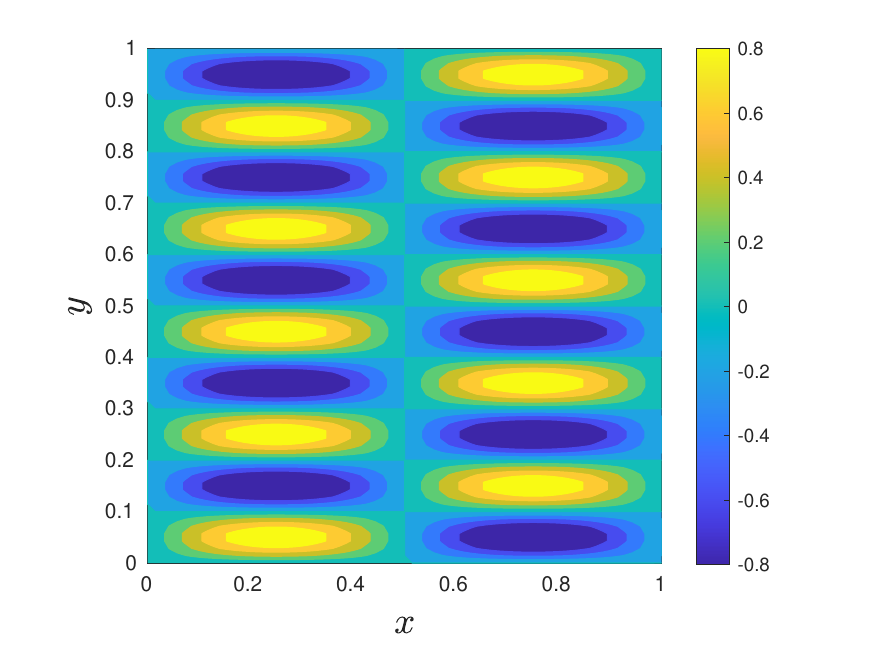}
	\caption{Surface plot and contour plot for the exact solution $u$ defined in \eqref{eq:sol_test_1}. In this case, we set $\parr=10$.}
	\label{fig:exact_u_1}
\end{figure}

\begin{figure}
	\centering
	\textbf{Comparison of stabilization terms - Test 1}
	\subfloat{\includegraphics[width=0.43\linewidth]{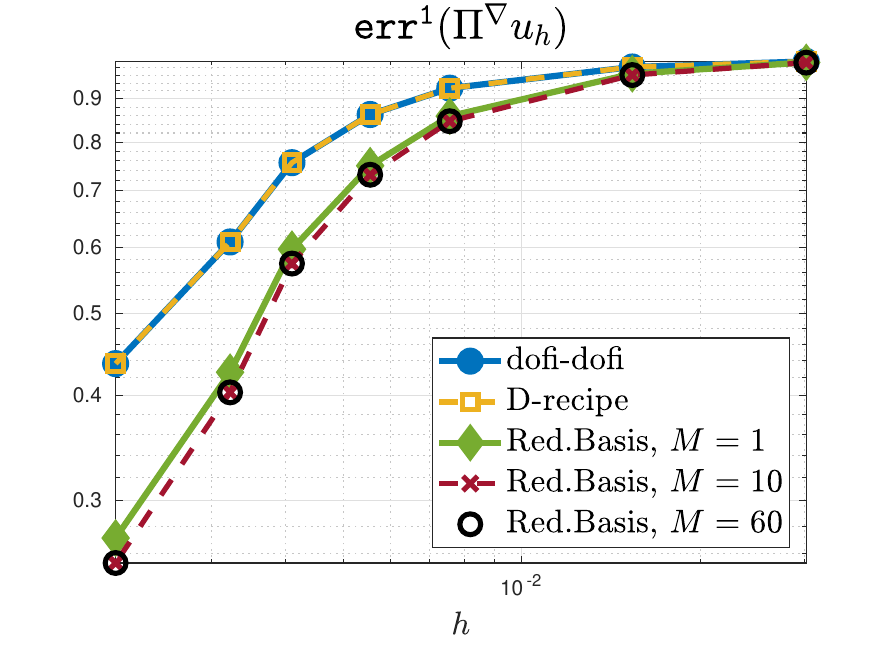}}
	\subfloat{\includegraphics[width=0.43\linewidth]{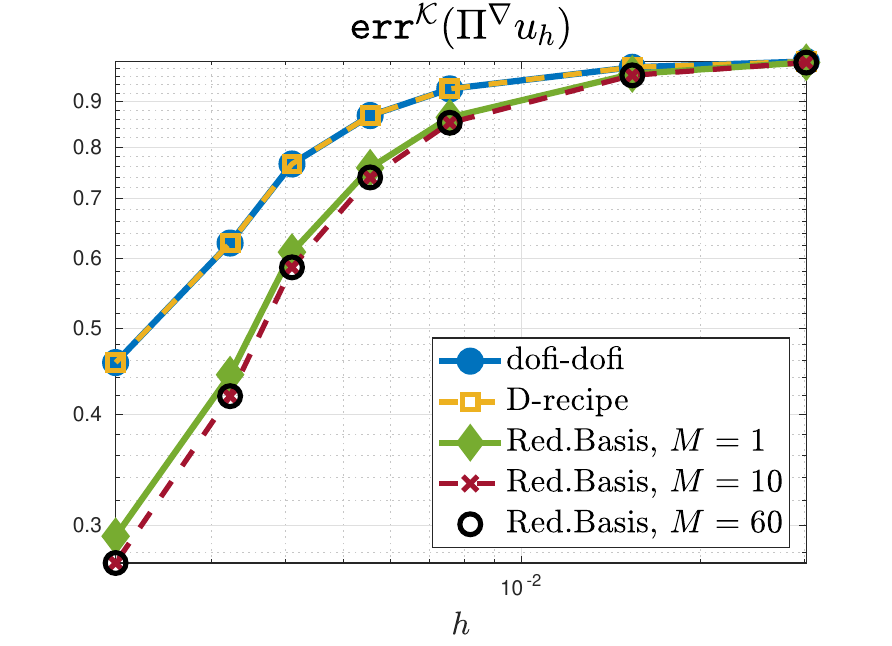}}
	\caption{Convergence plots for Test 1. From left to right, $\erruno{\PiNabla u_h}$ and $\erren{\PiNabla u_h}$. We denote in blue the convergence for \textit{dofi--dofi}, while yellow is used for \textit{D--recipe}. For the reduced basis stabilization, we denote in green the case $M=1$, in red $M=10$ and with black circles $M=60$. In both plots, it is evident that the reduced basis stabilization performs better than the standard choices because it is able to catch the anisotropy of the problem. \rosso{The proposed RB method is effective in the cheapest case since $M=1,10,60$ provide equivalent results.}}
\label{fig:stab_test_1}
\end{figure}

For this test, we choose the right hand side so that the following exact solution is obtained 
\begin{equation}\label{eq:sol_test_1}
	u(x,y) = \sin(2\pi x)\sin(\parr\pi y),\qquad \parr=80.
\end{equation}
The solution is characterized by high frequency oscillations in $y$ direction. We set the diffusivity tensor to be $\aa = \aa_1$.
A simplified plot of the exact solution $u$ when $\parr=10$ in reported in Figure~\ref{fig:exact_u_1}. Convergence plots for this test are collected in Figure~\ref{fig:stab_test_1}, while the convergence history is reported in Table~\ref{tab:errors_test1} for the case of \textit{dofi--dofi} stabilization and RB stabilization built considering $M=1$. We notice that the reduced basis stabilization improves the convergence properties of the method for both $\erruno{\PiNabla u_h}$ and $\erren{\PiNabla u_h}$, with respect to the use of \textit{dofi--dofi} and \textit{D--recipe} (which, in this case, produce equivalent results). This improvement can already be observed  when just one reduced basis function is considered ($M=1$). Moreover, we notice that by increasing the number of reduced bases involved in the stabilization (we consider, in  particular, $M=10,60$), we do not obtain a significant gain with respect to $M=1$, which is then sufficient to cheaply take into account the anisotropy of the problem.

\renewcommand{\arraystretch}{1.4}
\begin{table}
\begin{center}
	\begin{tabular}{c||rr|rr||rr|rr}
		\multicolumn{9}{c}{\textbf{Test 1 - Convergence history}}\\
		\hline
		& \multicolumn{4}{c||}{\em dofi--dofi} & \multicolumn{4}{c}{Reduced Basis, $M=1$} \\
		\hline
		$h$ & $\erruno{\PiNabla u_h}$ & Rate & $\erren{\PiNabla u_h}$ & Rate & $\erruno{\PiNabla u_h}$ & Rate & $\erren{\PiNabla u_h}$ & Rate  \\
		\hline
		3.021e-2	& 9.971e-1	& --	& 9.989e-1	& --	& 9.945e-1	& --	& 9.966e-1	& --   \\
		1.536e-2	& 9.822e-1	& 0.02	& 9.841e-1	& 0.02	& 9.643e-1	& 0.05	& 9.671e-1	& 0.04 \\
		7.569e-3	& 9.269e-1	& 0.08	& 9.304e-1	& 0.08	& 8.586e-1	& 0.16	& 8.641e-1	& 0.16 \\
		5.548e-3	& 8.625e-1	& 0.23	& 8.683e-1	& 0.22	& 7.498e-1	& 0.44	& 7.583e-1	& 0.42 \\
		4.094e-3	& 7.561e-1	& 0.43	& 7.659e-1	& 0.41	& 5.968e-1	& 0.75	& 6.095e-1	& 0.72 \\
		3.219e-3	& 6.088e-1	& 0.90	& 6.239e-1	& 0.85	& 4.264e-1	& 1.40	& 4.436e-1	& 1.32 \\
		2.063e-3	& 4.365e-1	& 0.74	& 4.578e-1	& 0.70	& 2.709e-1	& 1.02	& 2.919e-1	& 0.94 \\
		\hline
	\end{tabular}
\end{center}
\caption{Errors and convergence rates for Test 1. We focus on VEM with \emph{dofi--dofi} stabilization and reduced basis stabilization with $M=1$.}
\label{tab:errors_test1}
\end{table}
\renewcommand{\arraystretch}{1}
\subsubsection*{Test 2.} 
\begin{figure}
\centering
\includegraphics[width=0.45\linewidth]{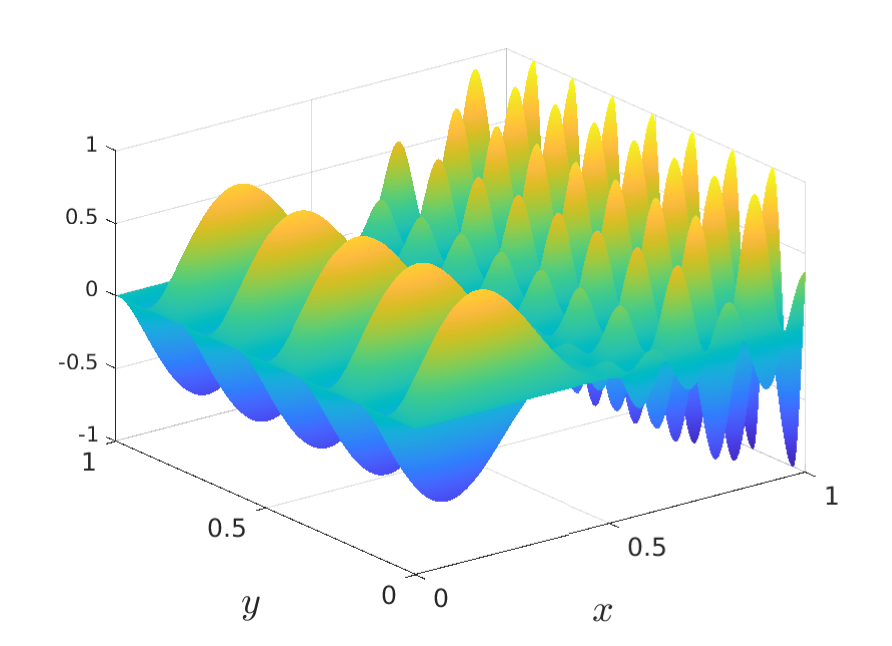}
\includegraphics[width=0.45\linewidth]{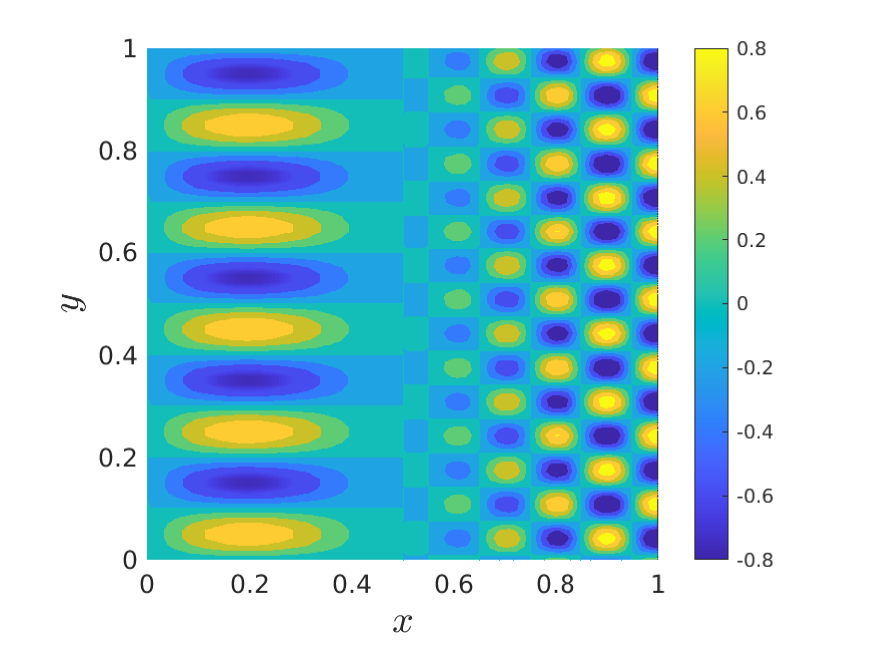}
\caption{Surface plot and contour plot for the exact solution $u$ defined in \eqref{eq:sol_test_2}. In this case, we set $\paru=10$ and $\paruu=15$.}
\label{fig:exact_u_2}
\end{figure}

\begin{figure}
\centering
\textbf{Comparison of stabilization terms - Test 2}
\subfloat{\includegraphics[width=0.43\linewidth]{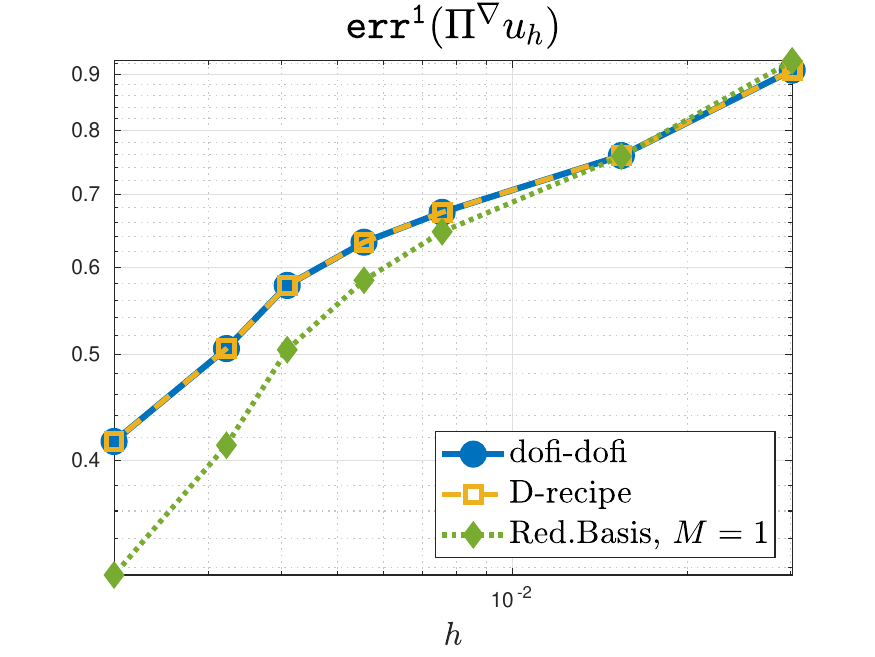}}
\subfloat{\includegraphics[width=0.44\linewidth]{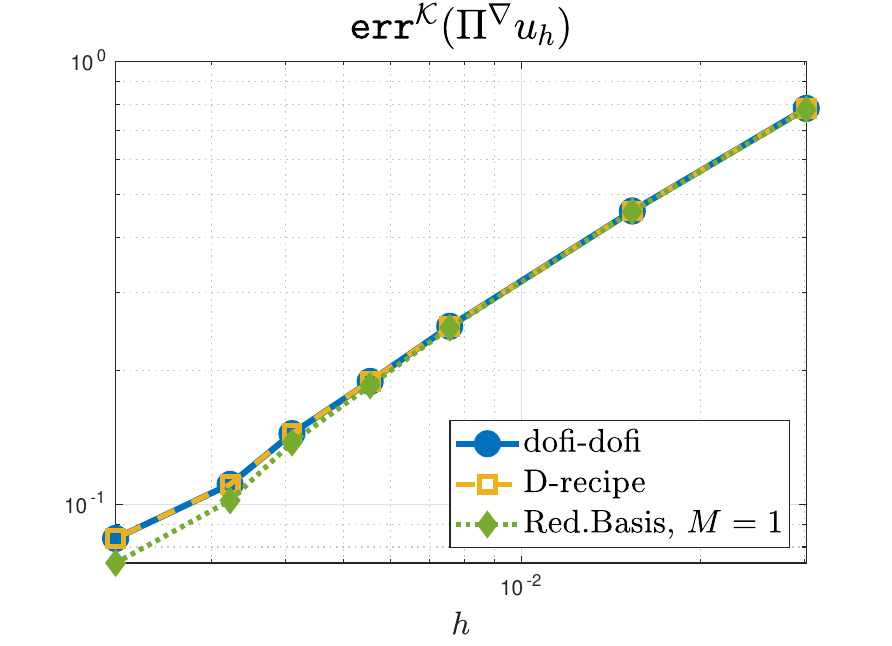}}
\caption{Convergence plots for Test 2. Same format of Figure~\ref{fig:stab_test_1}. Also in this test it is evident that the reduced basis stabilization improves the results of \textit{dofi--dofi} and \textit{D--recipe} in terms of $\erruno{\PiNabla u_h}$. Conversely, all the approaches behave similarly when $\erren{\PiNabla u_h}$ is considered.}
\label{fig:stab_test_2}
\end{figure}

For this test, the right hand side is defined according to the following solution and matrix $\Kmatr$. We choose a solution with discontinuous gradient when $x=1/2$. More precisely, $u$ is the continuous piecewise function 
\begin{equation}\label{eq:sol_test_2}
u(x,y) = \begin{cases}
	\sin(2\pi x)\sin(\paru\pi y)\cos(\pi x) & x\leq \frac{1}{2} \\
	\cos(\paru\pi x) \cos(\pi x) \sin(\paruu(\pi-y)\pi) & x > \frac{1}{2} \\
\end{cases},\qquad \paru=80,\,\paruu=30.
\end{equation}
As shown in Figure~\ref{fig:exact_u_2} for the simplified case with $\paru=10,\,\paruu=15$, this function oscillates in $y$ direction for $0\le x \le 1/2$ and in both $x$ and $y$ directions for $1/2 < x \le 1$. Then, we choose the following nonsymmetric positive definite matrix
\begin{equation*}
\Kmatr = \begin{bmatrix}
	1 & 10^{-2} \\
	5\cdot10^{-3} & 10^{-4}
\end{bmatrix}.
\end{equation*}
In this test, for the reduced basis stabilization, we plot only on the case $M=1$, since for larger values of $M$, \blu{analogously to what happens in the previous case, no significant improvement can be observed.} The results  are collected in Figure~\ref{fig:stab_test_2} and Table~\ref{tab:errors_test2} and are mostly in line with the results of the previous test. Once again, the method with reduced basis stabilization (with $M=1$) performs better than standard cases when $\erruno{\PiNabla u_h}$ is analyzed. Conversely, if we look at $\erren{\PiNabla u_h}$, all the method have same behavior until the last two finest cases, where the use of the reduced basis approach slightly improves.

\renewcommand{\arraystretch}{1.4}
\begin{table}
\begin{center}
	\begin{tabular}{c||rr|rr||rr|rr}
		\multicolumn{9}{c}{\textbf{Test 2 - Convergence history}}\\
		\hline
		& \multicolumn{4}{c||}{\em dofi--dofi} & \multicolumn{4}{c}{Reduced Basis, $M=1$} \\
		\hline
		$h$ & $\erruno{\PiNabla u_h}$ & Rate & $\erren{\PiNabla u_h}$ & Rate & $\erruno{\PiNabla u_h}$ & Rate & $\erren{\PiNabla u_h}$ & Rate  \\
		\hline
		3.021e-2	& 9.077e-1	& --	& 7.833e-1	& --	& 9.255e-1	& --	& 7.787e-1	& -- \\
		1.536e-2	& 7.589e-1	& 0.26	& 4.595e-1	& 0.79	& 7.576e-1	& 0.30	& 4.584e-1	& 0.78 \\
		7.569e-3	& 6.736e-1	& 0.17	& 2.522e-1	& 0.85	& 6.468e-1	& 0.22	& 2.501e-1	& 0.85 \\
		5.548e-3	& 6.327e-1	& 0.20	& 1.893e-1	& 0.93	& 5.842e-1	& 0.33	& 1.853e-1	& 0.97 \\
		4.094e-3	& 5.779e-1	& 0.30	& 1.443e-1	& 0.89	& 5.049e-1	& 0.48	& 1.377e-1	& 0.98 \\
		3.219e-3	& 5.062e-1	& 0.55	& 1.108e-1	& 1.09	& 4.133e-1	& 0.83	& 1.019e-1	& 1.25 \\
		2.063e-3	& 4.165e-1	& 0.44	& 8.364e-2	& 0.63	& 3.148e-1	& 0.61	& 7.360e-2	& 0.73 \\
		\hline
	\end{tabular}
\end{center}
\caption{Errors and convergence rates for Test 2. We focus on VEM with \emph{dofi--dofi} stabilization and reduced basis stabilization with $M=1$.}
\label{tab:errors_test2}
\end{table}
\renewcommand{\arraystretch}{1}

\newcommand{\vemglobRB}{V_h^{\text{rb}}}
\newcommand{\vemlocRB}{V^{\text{rb}}}
\newcommand{\nK}{\nu^K}
\newcommand{\Blocort}{\mathbb{B}^\perp(\bK)}

\begin{remark}
It is interesting to observe that the method resulting from the RB stabilization with $M$ basis functions can be interpreted as a fully conforming method in a suitable discretization space~$\vemglobRB$. Indeed, letting $\Blocort$ be defined as
\[
\Blocort = \{ v \in H^{1/2}(\bK) : \quad \int_{\bK} v \nabla q \cdot \nK \,ds= 0, \ \forall q \in \mathbb{P}_1(K), \ \int_{\bK} v\,ds = 0\}, 
\]
we have that $\PiNabla v = 0$ if and only if $v \in \Blocort$.
We can now let $W^\text{rb}(K) \subseteq H^1(K)$ be defined as
\[
W^{\text{rb}}(K) = \{ v \in \Span\{e_1^\text{rb},\cdots,e_N^{\text{rb}}\} :  
v|_{\bK} \in \Blocort
\},
\]
and let 
\[
\vemlocRB(K) = \mathbb{P}_1(K) \oplus W^{\text{rb}}(K).
\] 
Remark that we have that $\mathbb{P}_1(K) \subseteq 	\vemlocRB(K)$ by construction and that ${\vemlocRB}_{|\bK} = \Bkloc$. Thanks to~\eqref{splitting:2}, it is not difficult to check that, if we define
\[
\vemglobRB = \{v \in V : v_{|K} \in \vemlocRB(K)\ \text{ for all }K \in \mesh\},
\]
and we discretize Problem \ref{pro:cont_poisson} by a conforming Galerkin method, we obtain the same linear system as for the VEM method with RB stabilization.
\end{remark}

\section{Application II: Post-Processing}\label{sec:postproc}
Another possible use for the RB reconstruction of virtual element functions is the design of a post-processing technique where the whole $u_h$ is reconstructed, allowing, in particular, to retrieve a conformal approximation to the true solution out of the virtual element degrees of freedom. This can be used for visualization, for evaluating point values, or, for academic purpose, to compute the actual $H^ 1(\Omega)$ error with respect to a known benchmark solution. Also here the idea is to use the RB reconstruction on the non polynomial part of the solution, while evaluating the polynomial part by the standard VEM approach. 

\begin{algorithm}
	\caption{Reduced basis virtual functions reconstruction}
	\begin{flushleft}
		\textbf{Data:}\\
		\vspace{1mm}
		\begin{itemize}[leftmargin=15pt]
		\item $\elementVEM$: element of $\mesh$\\
		\item $\{u_h(\vv_j)\}_{j=1,\dots,N}$: dofs of numerical solution in $\elementVEM$\\
		\end{itemize}
		\vspace{2.5mm}
%
%
	Go to \textbf{Online phase} (see Algorithm~\ref{alg:rb_online}):\\
	\quad \textbf{Input:} $\elementVEM$\\
	\quad \textbf{Output:} RB approximation of VEM basis functions, $\widehat e_{M,j}^{\,\text{rb}}[K], j=1,\cdots,N$\\
	\vspace{2.5mm}
	
	Pull back $e_j^\elementVEM=\widehat e_{M,j}^{\,\text{rb}}[K] \circ \Bcal_\elementVEM$ on $\meshK=\Bcal^{-1}(\meshfine)$ in $\elementVEM$ for $j=,1\dots,N$\\
	
	Construct projection onto polynomials $\PiNabla u_h$ and evaluate on $\meshK$\\

	Compute $\uhRB = \PiNabla u_h + \sum_{j=1}^N (u_h(\vv_j)-\PiNabla u_h(\vv_j)) \, e_j^\elementVEM$ in $\elementVEM$
	\end{flushleft}
	\label{alg:rb_postproc}
\end{algorithm}

We demonstrate this approach on three examples where we solve Problem~\ref{pro:basis1} with $\Kmatr=I$ on \rev{$\Omega=(0,1)^2$}, once again discretized by Voronoi meshes. 

In the first example, once the equation is solved with VEM with {\em dofi-dofi} stabilization, we use the reduced basis method to reconstruct the solution and visualize it. In the second example we use the RB reconstruction for evaluating the solution along a line. The third example is an academic convergence test, where the true solution is known and the discrete solution is reconstructed element by element in the entire $\Omega$ to compare the convergence of the post-processed solution $\uhRB$ with the standard $\PiNabla u_h$ and the ``exact" reconstruction $u^\text{fe}_h$.

\subsubsection*{Visualization}

For this test, we discretize $\Omega$ with a relatively coarse Voronoi mesh $\mesh$ consisting of 100 elements, as represented in Figure~\ref{fig:mesh_rec}. For visualization purposes, we solve the Poisson problem by means of the standard lowest order virtual element method and then we recostruct a conforming solution via RB method. The right hand side $f$ is chosen in such a way that the continuous solution is 
\begin{equation}\label{eq:sol_post}
u(x,y) = \frac{1}{32\pi^2} \sin(4\pi x)\sin(4\pi y).
\end{equation}

\begin{figure}
\centering
\subfloat[$u$]{\includegraphics[width=0.32\linewidth]{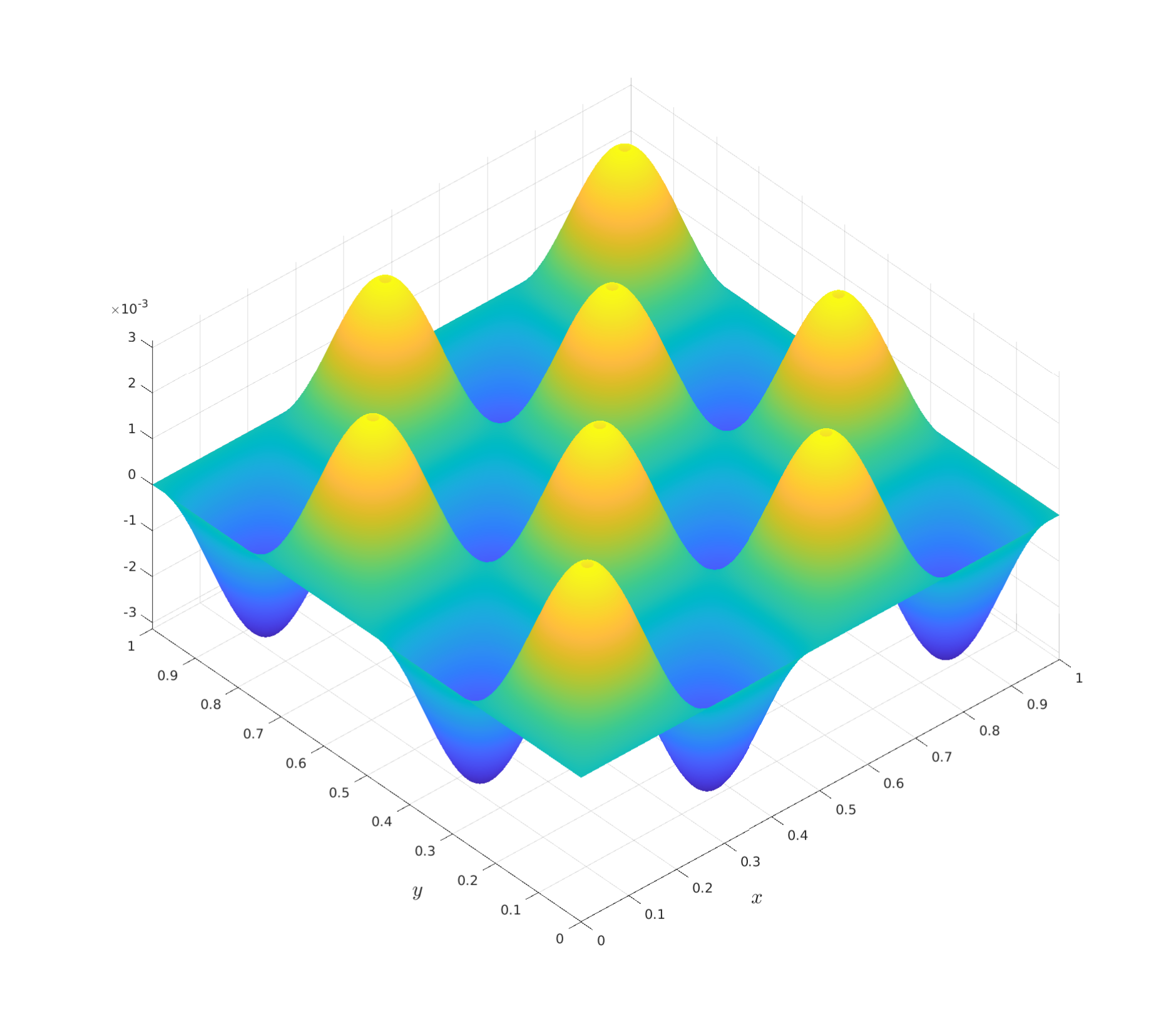}}
\subfloat[$\PiNabla u_h$]{\includegraphics[width=0.32\linewidth]{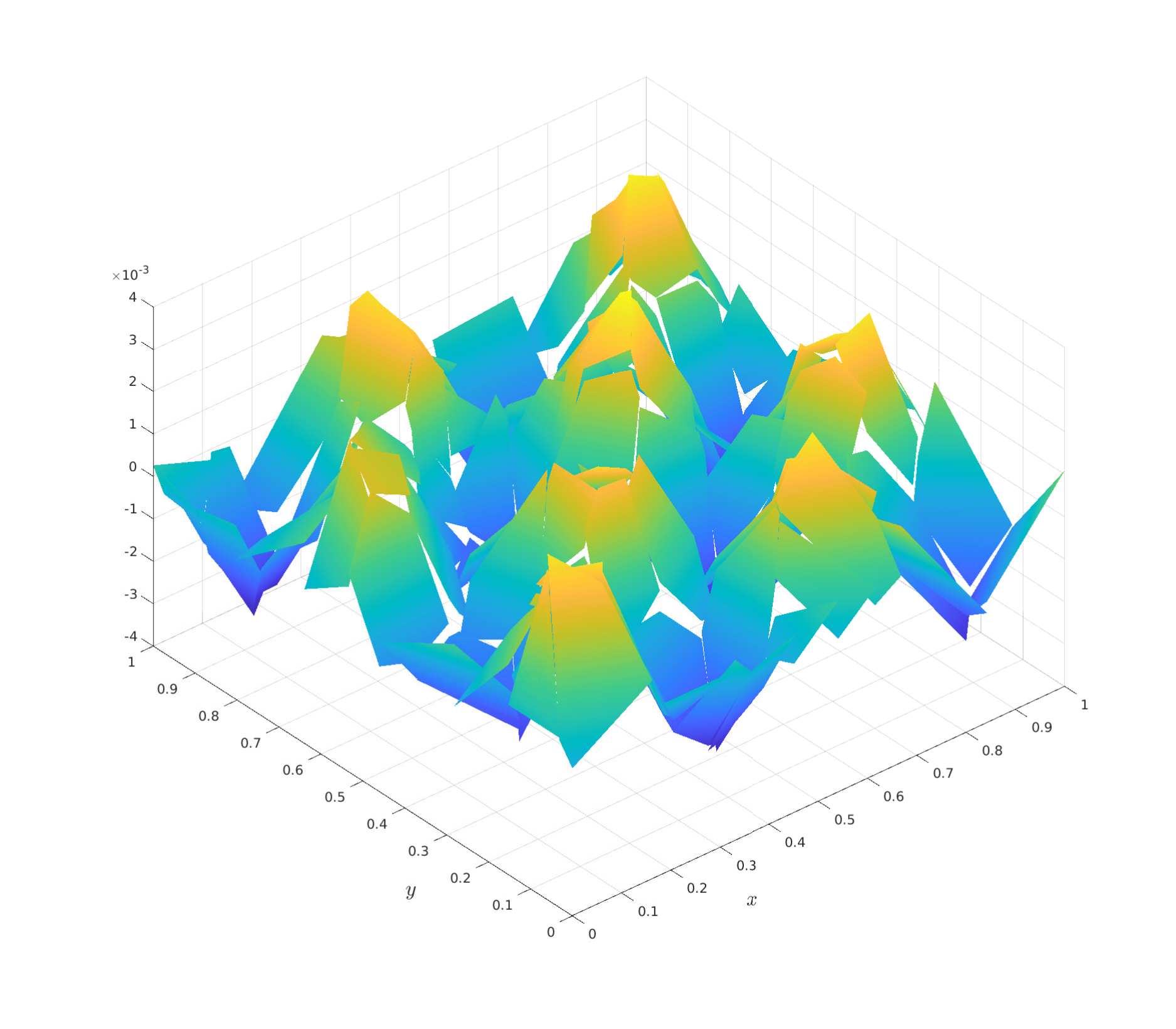}}
\subfloat[$\uhRB,\,M=1$]{\includegraphics[width=0.32\linewidth]{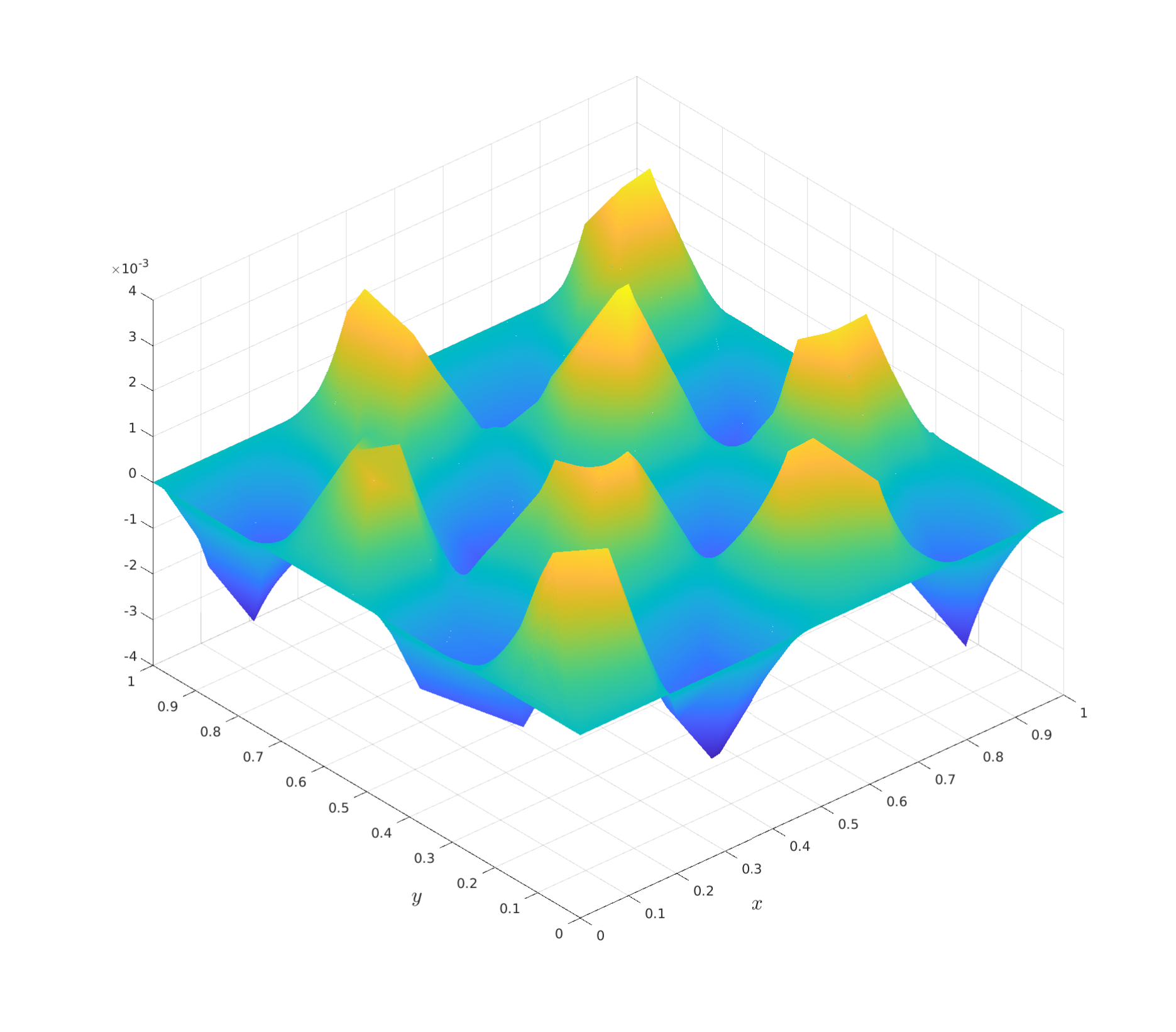}}
\caption{From left to right: the exact solution $u$ defined in \eqref{eq:sol_post}, the post-processed solution $\PiNabla u_h$ computed projecting onto polynomials in each element, the reconstructed solution $\uhRB$, computed with a single reduced basis in each element. The considered mesh is represented in Figure~\ref{fig:mesh_rec}.}
\label{fig:comparison_postproc}
\end{figure}

In Figure~\ref{fig:comparison_postproc}, we plot the reduced basis reconstruction $\uhRB$ (computed in each element of $\mesh$ with precision $M=1$) and we compare it with the standard polynomial projection $\PiNabla u_h$. The reduced basis reconstruction is conforming in the VEM space, while $\PiNabla u_h$ is discontinuous across the elements.

\begin{figure}
	\centering
	\subfloat[$u$]{\includegraphics[width=0.31\linewidth]{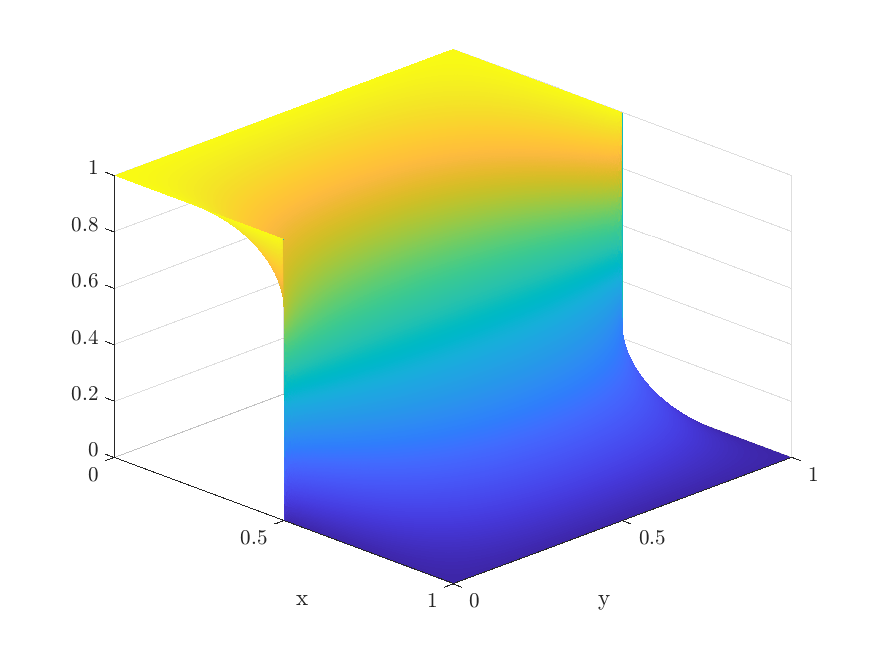}}
	\subfloat[$\PiNabla u_h$]{\includegraphics[width=0.32\linewidth]{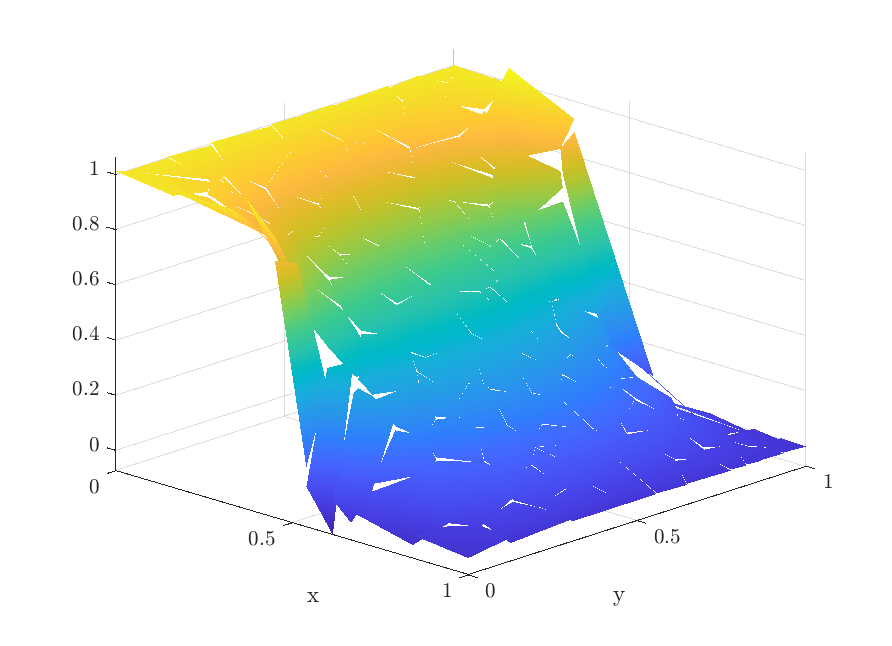}}
	\subfloat[$\uhRB,\,M=1$]{\includegraphics[width=0.32\linewidth]{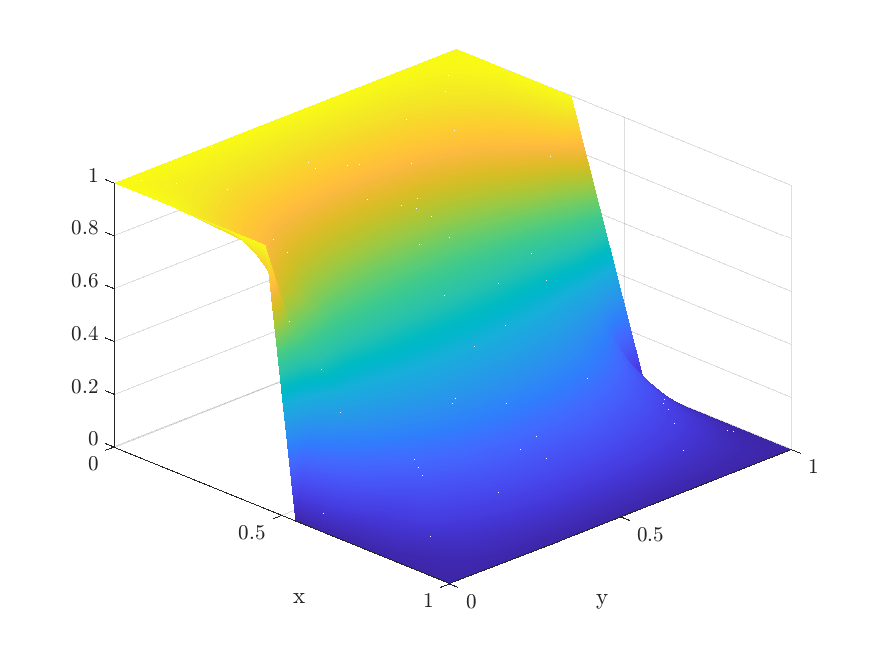}}
	\caption{\rev{From left to right: the ``true'' solution $u$, the post-processed solution $\PiNabla u_h$ computed by projecting onto polynomials in each element, the reconstructed solution $\uhRB$, computed with a single reduced basis in each element. The ``true'' solution $u$ is evaluated by a standard finite element method on a fine triangulation with size $5\cdot10^{-4}$.}}
	\label{fig:comparison_postproc2}
\end{figure}

\rev{We then perform a second visualization test by considering a less regular function: indeed, we solve again the Poisson equation in the unit square with right hand side $f$ being equal to zero and with the following boundary datum:
\begin{equation}
	u = \begin{cases}
		1 & x\le 1/2\\
		0 & x>1/2
	\end{cases}\qquad\text{on }\partial\Omega.
\end{equation}
Observe that this example falls outside of the standard theoretical framework even at continuous level, since the boundary datum does not belong to the space $H^1(\Omega)_{|\partial\Omega} = H^{1/2}(\partial\Omega)$, so that the solution does not belong to $H^1(\Omega)$. The considered mesh is again the one depicted in Figure~\ref{fig:mesh_rec}; we remark that the jump at the boundary does not coincide with any vertex of the mesh elements. We plot the post-processed solutions in Figure~\ref{fig:comparison_postproc2}: we clearly see the improvement resulting from using the RB reconstruction with respect to the standard VEM reconstruction by projection onto discontinuous polynomials.
}

\subsubsection*{Local reconstruction}
\begin{figure}
\centering
\textbf{Local reconstruction test}
\subfloat[Geometry \label{fig:mesh_rec}]{\includegraphics[width=0.4\linewidth]{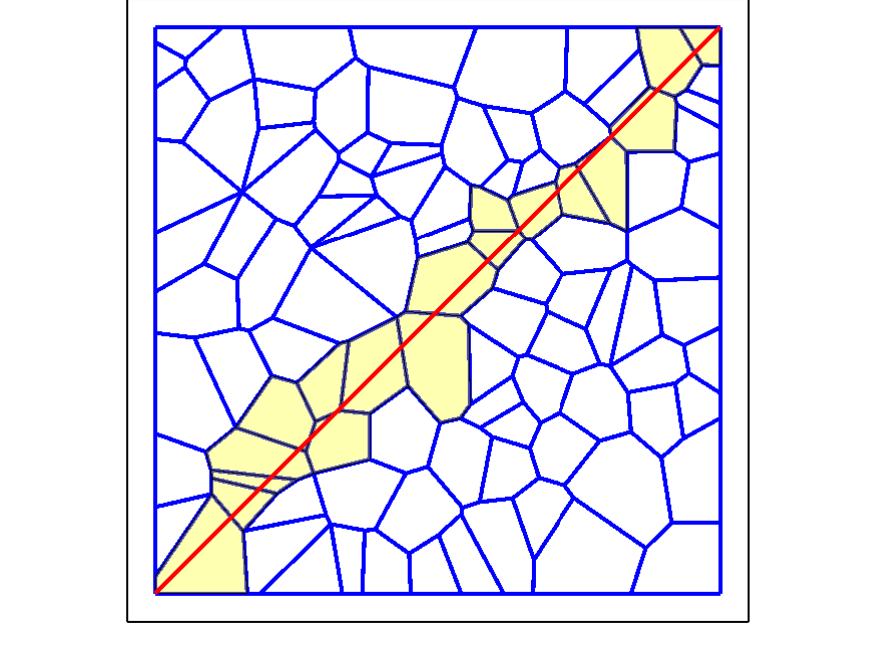}}
\subfloat[$\uhRB$ vs $u^\text{fe}_h$ \label{fig:rec_rb}]{\includegraphics[width=0.45\linewidth]{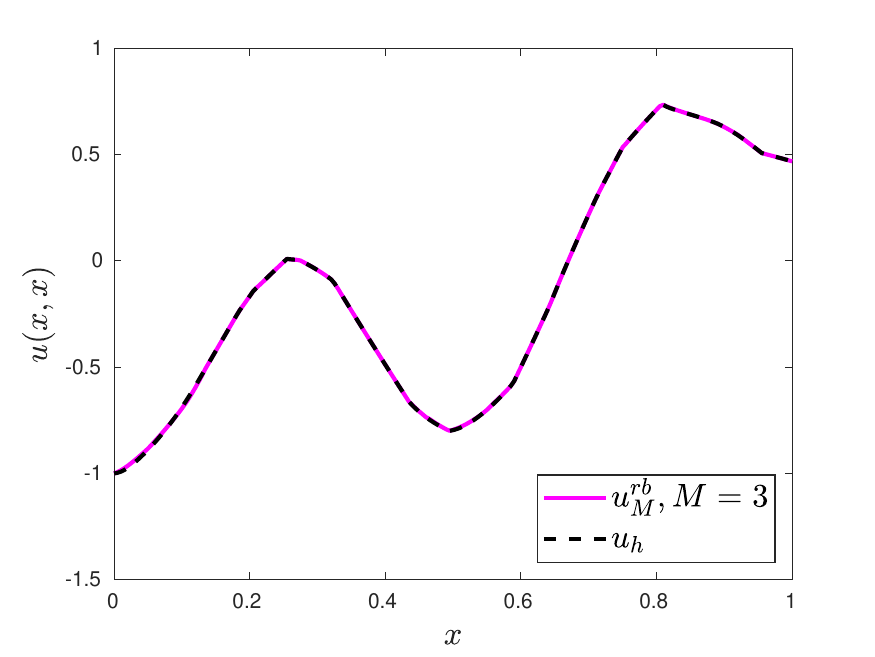}}\\
\subfloat[$\PiNabla u_h$ vs $u^\text{fe}_h$ \label{fig:rec_pi}]{\includegraphics[width=0.45\linewidth]{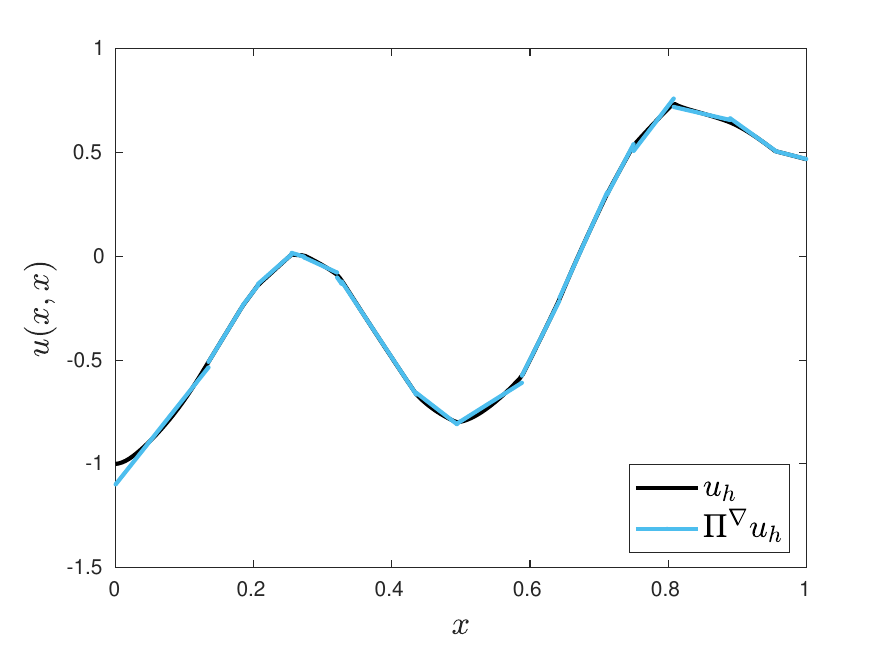}}
\subfloat[Zoom \label{fig:rec_zoom}]{\includegraphics[width=0.45\linewidth]{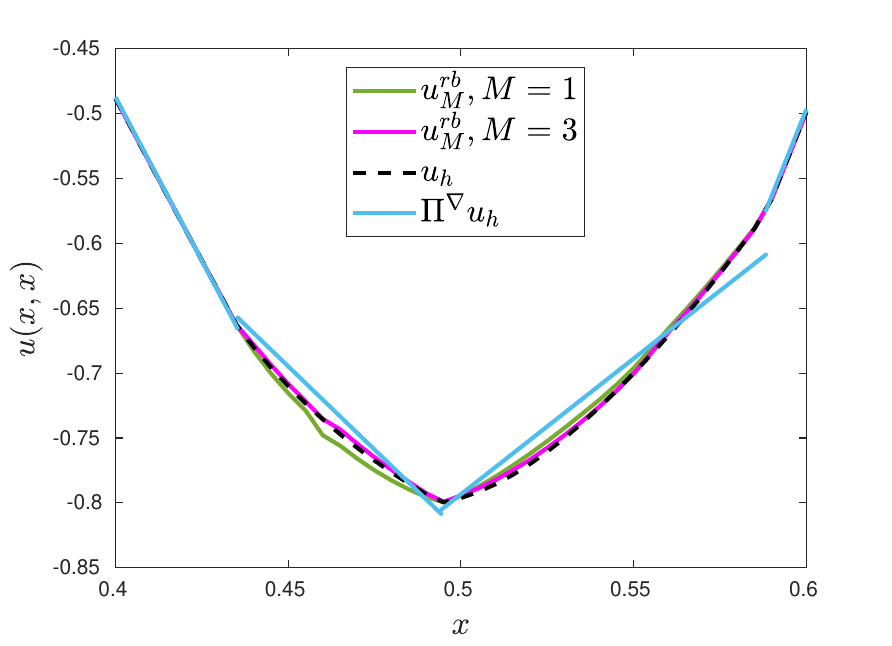}}
\caption{Geometry and results for the local reconstruction test. \textsc{(a)}~The problem is solved on a Voronoi mesh and the solution is reconstructed on the red diagonal. \textsc{(b)}~Comparison between the ``exact" reconstruction $u^\text{fe}_h$ (black line) and the reduced basis approximation $\uhRB,M=3$ (magenta line). \textsc{(c)}~Comparison between $u^\text{fe}_h$ and the projection $\PiNabla u_h$ (blue line). \textsc{(d)}~Zoom in the interval $[0.4,0.6]$ comparing $u^\text{fe}_h$, $\PiNabla u_h$, $\uhRB,M=3$ and in addition $\uhRB,M=1$ (green line). All the reconstructions are globally good: $\uhRB,M=1,3$ are conforming, while $\PiNabla u_h$ is discontinuous across the elements.}
\label{fig:reconstructions}
\end{figure}

In this second post-processing example, we solve again the problem under consideration on the mesh depicted in Figure~\ref{fig:mesh_rec}. The right hand side is chosen in such a way that the continuous solution is 
\begin{equation}
\begin{aligned}
	u(x,y) =\,&x^3-xy^2+yx^2+x^2-xy\\
	&\quad -x+y-1+\sin(5x)\sin(7y)+\log(1+x^2+y^4).
\end{aligned}
\end{equation}

Once the degrees of freedom of the discrete solution are computed, the function is reconstructed on the diagonal $y=x$ of the domain, which is marked in red in Figure~\ref{fig:mesh_rec}. We reconstruct the VEM basis functions in each polygon $\elementVEM$ intersecting the diagonal and then we evaluate the solution on a one dimensional grid. For the explicit finite elements computation of the virtual functions in~$\elementVEM$, we generate a triangulation $\meshK$ with size $\dd = h_\elementVEM/100$, while the reduced basis reconstruction~$\uhRB$ is computed for $M=1,3$.


We plot the post-processed solutions in Figure~\ref{fig:reconstructions}. In particular, in Figure~\ref{fig:rec_rb} we compare~$\uhRB$ built using $M=3$ (magenta) with $u^\text{fe}_h$ (black), whereas in Figure~\ref{fig:rec_pi} we compare $\PiNabla u_h$ (blue) with~$u^\text{fe}_h$. Finally, in Figure~\ref{fig:rec_zoom}, we zoom in the interval $[0.4,0.6]$ to highlight the different behaviors between the reconstructed solutions, considering also $\uhRB,M=1$ (green). \rosso{As already observed in the previous test, both $\uhRB$ and $\PiNabla u_h$ are good approximations of the solution: $\uhRB$ is conforming in the VEM space, whereas $\PiNabla u_h$ is discontinuous.}

%

\subsubsection*{Convergence test}

In this example, the reduced basis post-processing technique is analyzed in terms of convergence properties. We solve the Poisson problem on a sequence of Voronoi meshes $\meshfamily$. We select the exact solution \eqref{eq:sol_post} already chosen for the visualization test.

Once the problem is solved with VEM, we compute the reduced basis reconstruction $\uhRB$ with $M=1$ and we compare its convergence history with those given by the polynomial projection $\PiNabla u_h$ and the ``exact" reconstruction $u^\text{fe}_h$, constructed as in the previous section.

\begin{figure}
\centering
\textbf{Comparison of post-processing techniques}\par\medskip
\includegraphics[width=0.48\linewidth]{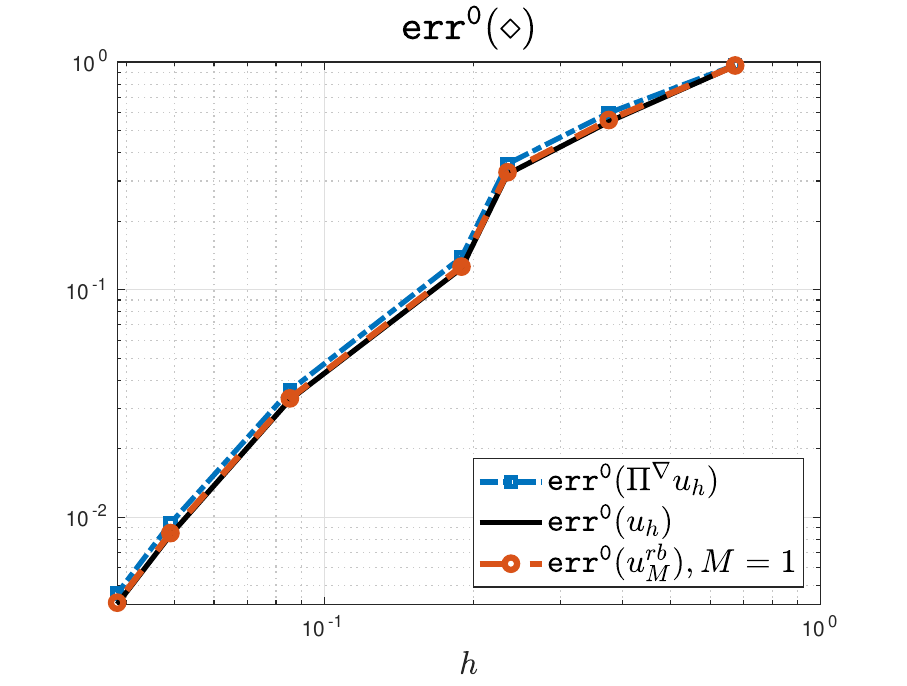}
\includegraphics[width=0.485\linewidth]{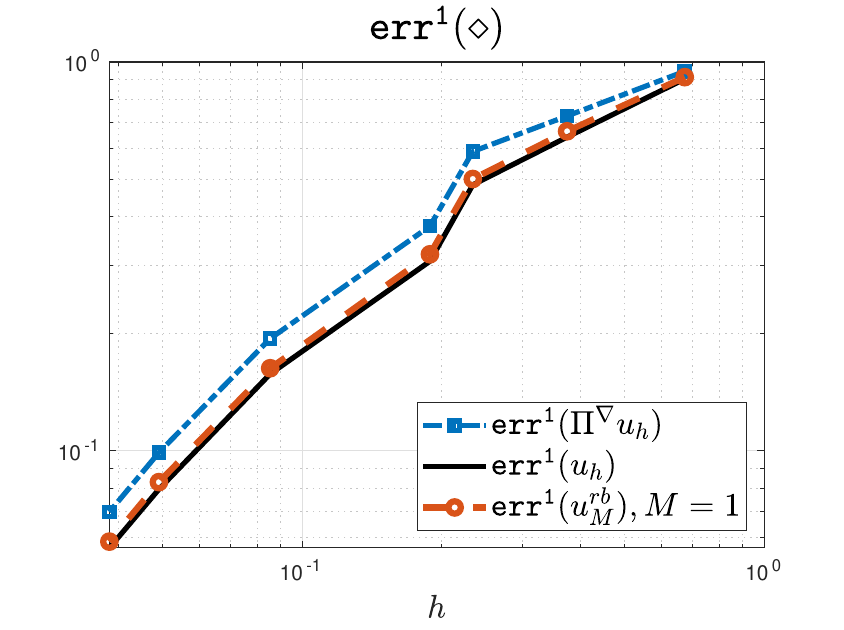}
\includegraphics[width=0.48\linewidth]{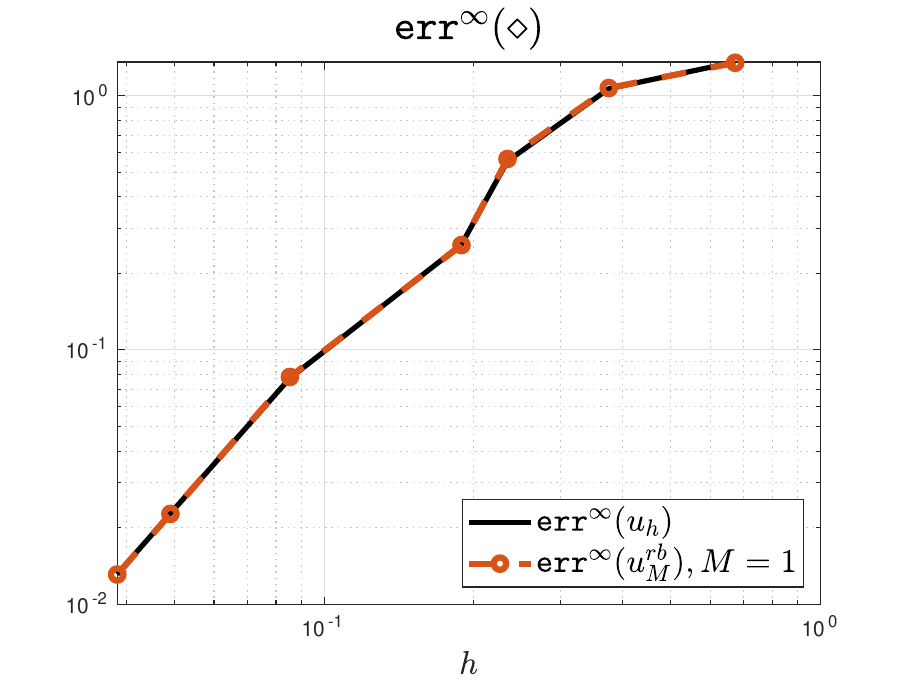}
\caption{Convergence plots comparing the error committed by $\PiNabla u_h$ (blue line), the ``exact'' reconstruction $u^\text{fe}_h$ (black line) and the reduced basis approximation $\uhRB$ (orange line) computed with a single basis. The reduced basis approximation behaves like $u^\text{fe}_h$ for all the considered errors.}
\label{fig:convergence}
\end{figure}

The results are studied in terms of $\errzero{\diamond},\erruno{\diamond},\errinf{\diamond}$ and the convergence plots are collected in Figure~\ref{fig:convergence}. It is clear that the reduced basis approximation (orange line) presents the same convergence history of $u^\text{fe}_h$ (black line), slightly improving the behavior of $\PiNabla u_h$ (blue line). Looking at $\errinf{\diamond}$, we see that the good behavior of $\uhRB$ is confirmed also pointwise, as already observed in the previous test when we carried out a local reconstruction.

	\section{Conclusions}\label{sec:conclusion}
	
We proposed the use of a reduced basis method in support of the virtual element method. By interpreting the partial differential equations involved in the definition of the elemental nodal basis functions as instances of a single parametric equation on a reference element, where the geometry of the physical element plays the role of parameter, it is indeed possible to use a reduced basis method to cheaply reconstruct more or less accurate approximations to such basis functions. These can then be used to reconstruct a conforming approximation to the discrete solution, or to evaluate output quantities. 
It can also be exploited to  evaluate a stabilization term that is closer to the actual contribution of the non polynomial part to the stiffness matrix. 
Depending on the accuracy of the reduced basis reconstruction, the resulting discrete method can be interpreted either as a virtual element method with a suitably designed stabilization, or as a computationally efficient polygonal finite element method based on approximate harmonic coordinates. 
We tested the proposed approach on different test problems, demonstrating its feasibility, effectiveness and computational affordability for two dimensional second order problems. 

Forthcoming work will include the application of the same idea to a larger class of problems, in particular two dimensional elasticity and fourth order problems, and to higher order two dimensional virtual elements, \rev{where we believe that, in  order to get full robustness in $k$,   it is necessary, within the design of the stabilization term, to take the interaction between boundary and interior degrees of freedom into account.} \rev{The main difficulty here stems from the high number of degrees of freedom, and, consequently, of basis functions to be constructed, that risks to entail a significant increase of the resulting  computational overhead. Moreover, as the degree $k$ increases, the solution of the partial differential equation defining the basis function, requires finer and finer meshes, thus significantly increasing the cost of the offline phase. 
} 

The same ideas can, in principle, be applied also in three dimensions, \rev{provided we can define a suitable set of reference polyhedra, playing the role that the $N$-edges regular polygons play in two dimensions. In some cases of interest, such as, for instance, for meshes of 8 nodes bricks with flat triangular faces, this is feasible, and, consequently, the proposed strategy can be quite easily adapted. More generally, however, the a priori choice of an exhaustive set of reference polyhedra, onto which all the elements of any mesh of a given class of meshes can be piecewise linearly mapped, is far from being trivial. It might then be necessary to resort to some kind of domain decomposition strategy, in the spirit of the Reduced Basis Element method \cite{M2AN_2013__47_1_213_0}, by decomposing the given element into the union of pyramid like elements, each with a polygonal basis (coinciding with one of the faces of the polyhedron under consideration), and with the opposite vertex in the element centroid. Whether the resulting computational overhead will be acceptable will have to be ascertained.
}

	\section{Acknowledgments}
	
	The authors are member of the INdAM -- GNCS research group and F. Credali is partially supported by CNR -- IMATI ``E. Magenes", Pavia (Italy). \medskip

	This paper has been realized in the framework of ERC Project CHANGE, which has received funding from the European Research Council (ERC) under the European Union’s Horizon 2020 research and innovation programme (grant agreement No 694515), and was co--funded by the MIUR Progetti di Ricerca di Rilevante Interesse Nazionale (PRIN) Bando 2017 (grant 201744KLJL) and Bando 2020 (grant 20204LN5N5).
	
	\bibliographystyle{abbrv}
	\bibliography{biblio}

\begin{thebibliography}{10}

\bibitem{antonietti2014stream}
P.~F. Antonietti, L.~Beir{\~a}o~da Veiga, D.~Mora, and M.~Verani.
\newblock A stream virtual element formulation of the {S}tokes problem on
  polygonal meshes.
\newblock {\em SIAM Journal on Numerical Analysis}, 52(1):386--404, 2014.

\bibitem{antonietti2023agglomeration}
P.~F. Antonietti, S.~Berrone, M.~Busetto, and M.~Verani.
\newblock Agglomeration-based geometric multigrid schemes for the virtual
  element method.
\newblock {\em SIAM Journal on Numerical Analysis}, 61(1):223--249, 2023.

\bibitem{multi}
P.~F. Antonietti, L.~Mascotto, and M.~Verani.
\newblock A multigrid algorithm for the p-version of the virtual element
  method.
\newblock {\em ESAIM: M2AN}, 52(1):337--364, 2018.

\bibitem{beirao2016virtual}
L.~Beir\~{a}o~da Veiga, F.~Brezzi, L.~D. Marini, and A.~Russo.
\newblock Virtual element method for general second-order elliptic problems on
  polygonal meshes.
\newblock {\em Mathematical Models and Methods in Applied Sciences},
  26(04):729--750, 2016.

\bibitem{beirao2013basic}
L.~Beir{\~a}o~da Veiga, F.~Brezzi, A.~Cangiani, G.~Manzini, L.~D. Marini, and
  A.~Russo.
\newblock Basic principles of virtual element methods.
\newblock {\em Mathematical Models and Methods in Applied Sciences},
  23(01):199--214, 2013.

\bibitem{da2013virtual}
L.~Beir{\~a}o Da~Veiga, F.~Brezzi, and L.~D. Marini.
\newblock Virtual elements for linear elasticity problems.
\newblock {\em SIAM Journal on Numerical Analysis}, 51(2):794--812, 2013.

\bibitem{beirao2014hitchhiker}
L.~Beir{\~a}o~da Veiga, F.~Brezzi, L.~D. Marini, and A.~Russo.
\newblock The hitchhiker's guide to the virtual element method.
\newblock {\em Mathematical models and methods in applied sciences},
  24(08):1541--1573, 2014.

\bibitem{da2016serendipity}
L.~Beir{\~a}o~da Veiga, F.~Brezzi, L.~D. Marini, and A.~Russo.
\newblock Serendipity nodal {VEM} spaces.
\newblock {\em Computers \& Fluids}, 141:2--12, 2016.

\bibitem{beir2016basic}
L.~Beir{\~a}o~da Veiga, A.~Chernov, L.~Mascotto, and A.~Russo.
\newblock Basic principles of hp virtual elements on quasiuniform meshes.
\newblock {\em Mathematical Models and Methods in Applied Sciences},
  26(08):1567--1598, 2016.

\bibitem{beirao2017stability}
L.~Beir{\~a}o~da Veiga, C.~Lovadina, and A.~Russo.
\newblock Stability analysis for the virtual element method.
\newblock {\em Mathematical Models and Methods in Applied Sciences},
  27(13):2557--2594, 2017.

\bibitem{da2017divergence}
L.~Beir{\~a}o~da Veiga, C.~Lovadina, and G.~Vacca.
\newblock Divergence free virtual elements for the {S}tokes problem on
  polygonal meshes.
\newblock {\em ESAIM: Mathematical Modelling and Numerical Analysis},
  51(2):509--535, 2017.

\bibitem{benedetto2016hybrid}
M.~F. Benedetto, S.~Berrone, A.~Borio, S.~Pieraccini, and S.~Scialo.
\newblock A hybrid mortar virtual element method for discrete fracture network
  simulations.
\newblock {\em Journal of Computational Physics}, 306:148--166, 2016.

\bibitem{berrone2021lowest}
S.~Berrone, A.~Borio, and F.~Marcon.
\newblock Lowest order stabilization free virtual element method for the
  {P}oisson equation.
\newblock {\em arXiv preprint arXiv:2103.16896}, 2021.

\bibitem{berrone2022comparison}
S.~Berrone, A.~Borio, and F.~Marcon.
\newblock Comparison of standard and stabilization free virtual elements on
  anisotropic elliptic problems.
\newblock {\em Applied Mathematics Letters}, 129:107971, 2022.

\bibitem{berrone2022virtual}
S.~Berrone and M.~Busetto.
\newblock A virtual element method for the two-phase flow of immiscible fluids
  in porous media.
\newblock {\em Computational Geosciences}, 26(1):195--216, 2022.

\bibitem{berrone2023virtual}
S.~Berrone, M.~Busetto, and F.~Vicini.
\newblock Virtual element simulation of two-phase flow of immiscible fluids in
  discrete fracture networks.
\newblock {\em Journal of Computational Physics}, 473:111735, 2023.

\bibitem{berrone2022efficient}
S.~Berrone and A.~Raeli.
\newblock Efficient partitioning of conforming virtual element discretizations
  for large scale discrete fracture network flow parallel solvers.
\newblock {\em Engineering Geology}, 306:106747, 2022.

\bibitem{bertoluzza2022stabilization}
S.~Bertoluzza, G.~Manzini, M.~Pennacchio, and D.~Prada.
\newblock Stabilization of the nonconforming virtual element method.
\newblock {\em Computers \& Mathematics with Applications}, 116:25--47, 2022.

\bibitem{bertoluzza2017bddc}
S.~Bertoluzza, M.~Pennacchio, and D.~Prada.
\newblock {BDDC} and {FETI}-{DP} for the virtual element method.
\newblock {\em Calcolo}, 54:1565--1593, 2017.

\bibitem{bertoluzza2020feti}
S.~Bertoluzza, M.~Pennacchio, and D.~Prada.
\newblock {FETI}-{DP} for the three dimensional virtual element method.
\newblock {\em SIAM Journal on Numerical Analysis}, 58(3):1556--1591, 2020.

\bibitem{bertoluzza2022interior}
S.~Bertoluzza, M.~Pennacchio, and D.~Prada.
\newblock Interior estimates for the virtual element method.
\newblock {\em arXiv preprint arXiv:2204.09955}, 2022.

\bibitem{bertoluzza2022weakly}
S.~Bertoluzza, M.~Pennacchio, and D.~Prada.
\newblock Weakly imposed dirichlet boundary conditions for 2d and 3d virtual
  elements.
\newblock {\em Computer Methods in Applied Mechanics and Engineering},
  400:115454, 2022.

\bibitem{brenner2017some}
S.~C. Brenner, Q.~Guan, and L.~Sung.
\newblock Some estimates for virtual element methods.
\newblock {\em Computational Methods in Applied Mathematics}, 17(4):553--574,
  2017.

\bibitem{brezzi2014basic}
F.~Brezzi, R.~S. Falk, and L.~D. Marini.
\newblock Basic principles of mixed virtual element methods.
\newblock {\em ESAIM: Mathematical Modelling and Numerical Analysis},
  48(4):1227--1240, 2014.

\bibitem{brezzi2013virtual}
F.~Brezzi and L.~D. Marini.
\newblock Virtual element methods for plate bending problems.
\newblock {\em Computer Methods in Applied Mechanics and Engineering},
  253:455--462, 2013.

\bibitem{cangiani2017hp}
A.~Cangiani, Z.~Dong, E.~Georgoulis, and P.~Houston.
\newblock hp-version discontinuous galerkin methods on polytopic meshes.
\newblock {\em to appear}, 2017.

\bibitem{cangiani2017posteriori}
A.~Cangiani, E.~H. Georgoulis, T.~Pryer, and O.~J. Sutton.
\newblock A posteriori error estimates for the virtual element method.
\newblock {\em Numerische mathematik}, 137:857--893, 2017.

\bibitem{cangiani2015hourglass}
A.~Cangiani, G.~Manzini, A.~Russo, and N.~Sukumar.
\newblock Hourglass stabilization and the virtual element method.
\newblock {\em International Journal for Numerical Methods in Engineering},
  102(3-4):404--436, 2015.

\bibitem{chen2018some}
L.~Chen and J.~Huang.
\newblock Some error analysis on virtual element methods.
\newblock {\em Calcolo}, 55:1--23, 2018.

\bibitem{cihan2022virtual}
M.~Cihan, B.~Hudobivnik, J.~Korelc, and P.~Wriggers.
\newblock A virtual element method for 3{D} contact problems with
  non-conforming meshes.
\newblock {\em Computer Methods in Applied Mechanics and Engineering},
  402:115385, 2022.

\bibitem{da2014mimetic}
L.~B. da~Veiga, K.~Lipnikov, and G.~Manzini.
\newblock {\em The mimetic finite difference method for elliptic problems},
  volume~11.
\newblock Springer, 2014.

\bibitem{DASSI20221}
F.~Dassi, A.~Fumagalli, A.~Scotti, and G.~Vacca.
\newblock Bend 3{D} mixed virtual element method for {D}arcy problems.
\newblock {\em Computers \& Mathematics with Applications}, 119:1--12, 2022.

\bibitem{dassi2020parallel}
F.~Dassi and S.~Scacchi.
\newblock Parallel block preconditioners for three-dimensional virtual element
  discretizations of saddle-point problems.
\newblock {\em Computer Methods in Applied Mechanics and Engineering},
  372:113424, 2020.

\bibitem{dassi2020parallel2}
F.~Dassi and S.~Scacchi.
\newblock Parallel solvers for virtual element discretizations of elliptic
  equations in mixed form.
\newblock {\em Computers \& Mathematics with Applications}, 79(7):1972--1989,
  2020.

\bibitem{dassi2022robust}
F.~Dassi, S.~Zampini, and S.~Scacchi.
\newblock Robust and scalable adaptive {BDDC} preconditioners for virtual
  element discretizations of elliptic partial differential equations in mixed
  form.
\newblock {\em Computer Methods in Applied Mechanics and Engineering},
  391:114620, 2022.

\bibitem{de2016nonconforming}
B.~A. de~Dios, K.~Lipnikov, and G.~Manzini.
\newblock The nonconforming virtual element method.
\newblock {\em ESAIM: Mathematical Modelling and Numerical Analysis},
  50(3):879--904, 2016.

\bibitem{DeRose2006HarmonicC}
T.~DeRose and M.~Meyer.
\newblock Harmonic coordinates.
\newblock 2006.

\bibitem{di2020hybrid}
D.~A. Di~Pietro and J.~Droniou.
\newblock The hybrid high-order method for polytopal meshes.
\newblock {\em Number 19 in Modeling, Simulation and Application}, 2020.

\bibitem{floater2015generalized}
M.~S. Floater.
\newblock Generalized barycentric coordinates and applications.
\newblock {\em Acta Numerica}, 24:161--214, 2015.

\bibitem{gain2014virtual}
A.~L. Gain, C.~Talischi, and G.~H. Paulino.
\newblock On the virtual element method for three-dimensional linear elasticity
  problems on arbitrary polyhedral meshes.
\newblock {\em Computer Methods in Applied Mechanics and Engineering},
  282:132--160, 2014.

\bibitem{hesthaven2016certified}
J.~S. Hesthaven, G.~Rozza, B.~Stamm, et~al.
\newblock {\em Certified reduced basis methods for parametrized partial
  differential equations}, volume 590.
\newblock Springer, 2016.

\bibitem{Hormann:2017:GBC}
K.~Hormann and N.~Sukumar, editors.
\newblock {\em Generalized Barycentric Coordinates in Computer Graphics and
  Computational Mechanics}.
\newblock CRC Press, Boca Raton, FL, 2017.

\bibitem{lions2012non}
J.~L. Lions and E.~Magenes.
\newblock {\em Non-homogeneous boundary value problems and applications: Vol.
  1}, volume 181.
\newblock Springer Science \& Business Media, 2012.

\bibitem{manzini2014new}
G.~Manzini, A.~Russo, and N.~Sukumar.
\newblock New perspectives on polygonal and polyhedral finite element methods.
\newblock {\em Mathematical Models and Methods in Applied Sciences},
  24(08):1665--1699, 2014.

\bibitem{martin2008polyhedral}
S.~Martin, P.~Kaufmann, M.~Botsch, M.~Wicke, and M.~Gross.
\newblock Polyhedral finite elements using harmonic basis functions.
\newblock In {\em Computer graphics forum}, volume~27, pages 1521--1529. Wiley
  Online Library, 2008.

\bibitem{mascotto2018ill}
L.~Mascotto.
\newblock Ill-conditioning in the virtual element method: Stabilizations and
  bases.
\newblock {\em Numerical Methods for Partial Differential Equations},
  34(4):1258--1281, 2018.

\bibitem{mascotto2019nonconforming}
L.~Mascotto, I.~Perugia, and A.~Pichler.
\newblock A nonconforming {T}refftz virtual element method for the {H}elmholtz
  problem.
\newblock {\em Mathematical Models and Methods in Applied Sciences},
  29(09):1619--1656, 2019.

\bibitem{mascotto2019nonconforming2}
L.~Mascotto, I.~Perugia, and A.~Pichler.
\newblock A nonconforming {T}refftz virtual element method for the {H}elmholtz
  problem: numerical aspects.
\newblock {\em Computer Methods in Applied Mechanics and Engineering},
  347:445--476, 2019.

\bibitem{mascotto2022nonconforming}
L.~Mascotto, I.~Perugia, and A.~Pichler.
\newblock The nonconforming {T}refftz virtual element method: general setting,
  applications, and dispersion analysis for the {H}elmholtz equation.
\newblock In {\em The Virtual Element Method and its Applications}, pages
  363--410. Springer, 2022.

\bibitem{mengolini2019engineering}
M.~Mengolini, M.~F. Benedetto, and A.~M. Arag{\'o}n.
\newblock An engineering perspective to the virtual element method and its
  interplay with the standard finite element method.
\newblock {\em Computer Methods in Applied Mechanics and Engineering},
  350:995--1023, 2019.

\bibitem{nguyen2018virtual}
V.~M. Nguyen-Thanh, X.~Zhuang, H.~Nguyen-Xuan, T.~Rabczuk, and P.~Wriggers.
\newblock A virtual element method for {2D} linear elastic fracture analysis.
\newblock {\em Computer Methods in Applied Mechanics and Engineering},
  340:366--395, 2018.

\bibitem{M2AN_2013__47_1_213_0}
D.~B. Phuong~Huynh, D.~J. Knezevic, and A.~T. Patera.
\newblock A {Static} condensation {Reduced} {Basis} {Element} method :
  approximation and \protect\emph{a posteriori }error estimation.
\newblock {\em ESAIM: Mathematical Modelling and Numerical Analysis},
  47(1):213--251, 2013.

\bibitem{prada2019feti}
D.~Prada, S.~Bertoluzza, M.~Pennacchio, and M.~Livesu.
\newblock {FETI}-{DP} preconditioners for the virtual element method on general
  2{D} meshes.
\newblock In {\em Numerical Mathematics and Advanced Applications ENUMATH
  2017}, pages 157--164. Springer, 2019.

\bibitem{sukumar2006recent}
N.~Sukumar and E.~Malsch.
\newblock Recent advances in the construction of polygonal finite element
  interpolants.
\newblock {\em Archives of Computational Methods in Engineering}, 13:129--163,
  2006.

\bibitem{sutton2017virtual}
O.~J. Sutton.
\newblock The virtual element method in 50 lines of {MATLAB}.
\newblock {\em Numerical Algorithms}, 75(4):1141--1159, 2017.

\bibitem{polymesher}
C.~Talischi, G.~H. Paulino, A.~Pereira, and I.~Menezes.
\newblock {P}oly{M}esher: a general-purpose mesh generator for polygonal
  elements written in {M}atlab.
\newblock {\em Structural and Multidisciplinary Optimization}, 45:309--328,
  2012.

\bibitem{vacca2018h}
G.~Vacca.
\newblock An ${H}^1$-conforming virtual element for {D}arcy and {B}rinkman
  equations.
\newblock {\em Mathematical Models and Methods in Applied Sciences},
  28(01):159--194, 2018.

\bibitem{valtr1994probability}
P.~Valtr.
\newblock Probability that n random points are in convex position.
\newblock 1994.

\bibitem{Vendorschot}
S.~Vendoschot.
\newblock Generating random convex polygons, 2017.

\bibitem{wriggers2017low}
P.~Wriggers and B.~Hudobivnik.
\newblock A low order virtual element formulation for finite elasto-plastic
  deformations.
\newblock {\em Computer Methods in Applied Mechanics and Engineering},
  327:459--477, 2017.

\bibitem{wriggers2017efficient}
P.~Wriggers, B.~D. Reddy, W.~Rust, and B.~Hudobivnik.
\newblock Efficient virtual element formulations for compressible and
  incompressible finite deformations.
\newblock {\em Computational Mechanics}, 60:253--268, 2017.

\bibitem{yu2022mvem}
Y.~Yu.
\newblock m{VEM}: A {MATLAB} software package for the virtual element methods.
\newblock {\em arXiv preprint arXiv:2204.01339}, 2022.

\end{thebibliography}

\end{document}